\documentclass[a4,12pt]{scrartcl}

\usepackage{xcolor} 
\definecolor{ocre}{RGB}{52,177,201} 

\usepackage{diagbox}
\usepackage{ stmaryrd }

\usepackage{listings}
\usepackage{afterpage}

\usepackage{avant} 
\usepackage{mathptmx} 

\usepackage{microtype} 
\usepackage[utf8]{inputenc} 
\usepackage[T1]{fontenc} 

\usepackage{float}
\usepackage{caption}
\usepackage{subcaption}
\usepackage[final]{graphicx}

\usepackage{tikz}
\usepackage{pgfplots} 
\usepackage{pgfplotstable}
\usetikzlibrary{math}
\usetikzlibrary{calc}
\usepackage{tkz-euclide}
\usetikzlibrary{patterns}
\usetikzlibrary{decorations.pathmorphing,shapes.geometric,arrows}
\usepgfplotslibrary{external} 
\usepackage{mdwlist}
\usepackage{multirow}
\usepackage{siunitx}

\protected\def\psverb#1{\def\innerpsverb##1#1{\texttt{##1}}\innerpsverb}

\usepackage[edges]{forest}

\def\Size{4pt}
\tikzset{
	folder/.pic={
		\filldraw[draw=folderborder,top color=folderbg!50,bottom color=folderbg]
		(-1.05*\Size,0.2\Size+5pt) rectangle ++(.75*\Size,-0.2\Size-5pt);  
		\filldraw[draw=folderborder,top color=folderbg!50,bottom color=folderbg]
		(-1.15*\Size,-\Size) rectangle (1.15*\Size,\Size);
	}
}

\usepackage{graphicx} 

\usepackage{tikz} 

\usepackage[english]{babel} 

\usepackage{enumitem} 
\setlist{nolistsep} 

\usepackage{booktabs} 

\usepackage{eso-pic} 

\usepackage[titletoc]{appendix} 

\definecolor{folderbg}{RGB}{124,166,198}
\definecolor{folderborder}{RGB}{110,144,169}
\definecolor{foldercolor}{RGB}{124,166,198}

\definecolor{folderbg}{RGB}{124,166,198}
\definecolor{folderborder}{RGB}{110,144,169}

\definecolor{Navy}{RGB}{000,000,128}
\definecolor{DarkRed}{RGB}{139,000,000}
\definecolor{DarkOrange}{RGB}{205,102,000}
\definecolor{DarkGreen}{RGB}{000,100,000}
\definecolor{Beige}{RGB}{250,250,240}
\definecolor{applegreen}{rgb}{0.55, 0.71, 0.0}
\definecolor{palecopper}{rgb}{0.85, 0.54, 0.4}
\definecolor{platinum}{rgb}{0.9, 0.89, 0.89}
\definecolor{CAUpurple}{RGB}{153,000,153}
\definecolor{CAUblue}{RGB}{051,051,153}
\definecolor{DarkPurple}{rgb}{0.07, 0.04, 0.56}

\definecolor{U0}{RGB}{183,028,028}
\definecolor{U1}{RGB}{136,014,079}
\definecolor{U2}{RGB}{074,020,140}
\definecolor{U3}{RGB}{049,027,146}
\definecolor{U4}{RGB}{026,035,126}
\definecolor{U5}{RGB}{013,071,161}
\definecolor{U6}{RGB}{001,087,155}
\definecolor{U7}{RGB}{000,096,100}
\definecolor{U8}{RGB}{000,077,064}
\definecolor{U9}{RGB}{027,094,032}

\definecolor{codegreen}{rgb}{0,0.6,0}
\definecolor{codegray}{rgb}{0.5,0.5,0.5}
\definecolor{codepurple}{rgb}{0.58,0,0.82}
\definecolor{backcolour}{rgb}{0.95,0.95,0.92}

\lstdefinestyle{mystyle}{
	backgroundcolor=\color{backcolour},   
	commentstyle=\color{codegreen},
	keywordstyle=\color{magenta},
	numberstyle=\tiny\color{codegray},
	stringstyle=\color{codepurple},
	basicstyle=\ttfamily\footnotesize,
	breakatwhitespace=false,         
	breaklines=true,                 
	captionpos=b,                    
	keepspaces=true,                 
	numbers=left,                    
	numbersep=5pt,                  
	showspaces=false,                
	showstringspaces=false,
	showtabs=false,                  
	tabsize=2
}

\usepackage{fancyhdr} 

\usepackage{amsmath,amsfonts,amssymb,amsthm} 

\usepackage[backend=bibtex,
citestyle=numeric]{biblatex}

\bibliography{sources/sources}

\usepackage{hyperref}

\hypersetup{
    bookmarks=true,
    unicode=false,
    pdftoolbar=true,
    pdfmenubar=true,
    pdffitwindow=false,
    pdfstartview={FitH},
    pdftitle={My title},
    pdfauthor={B. Philippi},
    pdfsubject={Subject},
    pdfcreator={Creator},
    pdfproducer={Producer},
    pdfkeywords={keyword1} {key2} {key3},
    pdfnewwindow=true,
    colorlinks=false,
    linkbordercolor={0 0 1},
    linkcolor=blue,
    citecolor=green,
    filecolor=magenta,
    urlcolor=cyan
}

\usepackage[left=3cm,right=3cm,top=2cm,bottom=2cm,includeheadfoot]{geometry}

\title{A Micro-Macro Parareal Implementation for the Ocean-Circulation Model FESOM2}
\author{B. $\text{Philippi}^1$, T. $\text{Slawig}^2$} 

\usepackage{csquotes}

\begin{document}

\maketitle

\begin{abstract}
	A micro-macro variant of the parallel-in-time algorithm Parareal has been applied to the ocean-circulation and sea-ice model model FESOM2. The state-of-the-art software in climate research has been developed by the Alfred-Wegener-Institut (AWI) in Bremen, Germany. The algorithm requires two meshes of low and high spatial resolution to define the coarse and fine propagator. As a first assessment we refined the PI mesh, increasing its resolution by factor 4. The main objective of this study was to demonstrate that micro-macro Parareal can provide convergence in diagnostic variables in complex climate research problems. After the introduction to FESOM2 we show how to generate the refined mesh and which interpolation methods were chosen. With the convergence results presented we discuss the success of this attempt and which steps have to be taken to extend the approach to current research problems.
\end{abstract}

\vspace{0.25cm}
\small
\textbf{${}^1$ \textit{Christian-Albrecht-Universität Kiel, Dept. of Computer Science, b.k.philippi@gmail.com}}

\textbf{${}^2$ \textit{Christian-Albrecht-Universität Kiel, Dept. of Computer Science, ts@informatik.uni-kiel.de}}
\normalsize

\newpage

\section{Introduction}

Predicting the earths future climate by numerical simulation represents one of the most urgent scientific tasks of our time. The urgency for understanding global climate change is matched by the complexity of this undertaking. In order to obtain a holistic representation a multitude of geophysical models are coupled in order to account for the different contributions to climate change, e.g. atmosphere, ocean-circulation, bio-geochemistry and human influence. An introduction into the complexity of geophysical science is given in the working group 1 contribution to the IPCC 6 report \cite{IPCC}. Numerical simulation of coupled climate systems requires considerable amount of time, even with the advancements made in parallel computing on high-performance-clusters (HPC). The necessity to cover large time intervals from decades to centuries, and even longer, extend the possibilities for run-time reductions by classical domain decomposition techniques. In order to be able to further exploit the resources for parallelization on modern HPC, the concept of an parallel-in-time approach offers a possibility to further reduce computational time. In this paper a micro-macro variant of the time-parallel algorithm Parareal is applied to the ocean-circulation and sea-ice model FESOM2. The algorithm was first introduced in \cite{Maday2001} and has been applied to vast variety of problems since then. There have been studies for climate problems reduced in complexity \cite{CaldasEtAl2023}, \cite{HamonEtAl2020}, \cite{Samuel2012} that successfully implemented Parareal and modifications to it. But, to our knowledge, no application to a state-of-the-art climate model software has been achieved so far. In an attempt to reduce the overwhelming complexity of coupled climate simulations \cite{FESOM2echam6} we focused on the ocean model FESOM2, developed at the Alfred-Wegener-Institut (AWI) \cite{FESOM2standardSim},\cite{FESOM2oce}. In a former study we applied Parareal to FESOM2 on a test mesh configuration in order to approximate prognostic variables, like temperature distributions local in time and space \cite{philippi2022parareal}. From the findings we learned that for climate research diagnostic variables, like annual mean temperature, have to be targeted and that a micro-macro approach to the problem is needed. For direct numerical simulation (DNS) of turbulent flows the convergence to diagnostic, statistical quantities by micro-macro Parareal was able to generate speedups \cite{Lunet2018}. The micro-macro Parareal variant in this study was introduced and tested for an energy balance model in \cite{Slawig2018}, where a 0-D model for macroscopic scale and an 1-D model for the microscopic scale were used. For the macro-scale we used the low resolution mesh PI, that is used for test purposes only. For the micro-scale we build a new mesh with higher resolution on basis of the PI mesh. By doing so we preserved the geometry, which is necessary for conserving energy and momentum for the interpolation between the meshes. \\
In this paper we will introduce the Parareal variant according to \cite{Slawig2018}. Afterwards a brief introduction to FESOM2 is given before we give an overview to FESOM2 meshes and how the PI mesh was refined. Since interpolation methods are indispensable for the proposed time-parallel approach, we present the possibilities for interpolation between unstructured meshes how we applied them. We conclude the preliminary work with an evaluation of the new mesh and compare the results to FESOM2 simulations on the PI mesh and the high resolution mesh CORE. The numerical experiments will be evaluated with respect to convergence to a selection of diagnostic variables. Based on the findings we will discuss the promising results and how to improve the approach in order to approximate standard test cases of FESOM2 in climate research. By giving insight to the process of preparing the setup for a Parareal application to FESOM2, we hope to give a glimpse into the complexity into climate research and the great amount of work necessary to make a test setting work.

\section{Micro-Macro Parareal Algorithm}

The Parareal algorithm as it was introduced by \cite{Maday2001} is concerned with the time-parallel approximation to solutions of initial value problems of the form:
\begin{equation}
	u_t(t) \; = \; f(u(t),t) \; , \quad t \in [0,T] \; , \quad u(0) = u_0 \; .
	\label{EQ:IVP}
\end{equation}
The algorithm exploits the decomposition of the time domain $T$ into $N_t$ time slices $\Delta T = [t_n,t_{n+1}]$ by iteratively computing approximations to the solution $u(t)$ of Eq.\ref{EQ:IVP} for each time slice. The algorithm requires two time integration methods of different order. The fine propagator $F$ is commonly defined as an accurate high order method with high computational cost. In contrast a coarse solver of low order is required, for which numerical accuracy is secondary to being fast. Both solver provide approximations to the solution $u(t_{n+1})$ from a given initial value $U_n$ at $_tn$ based on their properties:
\begin{equation}
	\begin{split}
		F_{n+1} \; &:= \; F(U_n) \; \approx \; u(t_{n+1}) \; , \\
		G_{n+1} \; &:= \; G(U_n) \; \approx \; u(t_{n+1}) \; .
	\end{split} 
\end{equation}

The $F$ propagator is generally applied in serial over the entire time interval $[0,T]$ and therefore represents the reference solution to the IVP in Eq.\ref{EQ:IVP}. The concept of Parareal is to iterate an approximate solution on each time slice $\Delta T$ in parallel. To allow for time-parallel execution every time slice needs initial values provided by the fast propagator $G$ over the entire time domain:
\begin{equation}
	\begin{split}
		U^0_0 \; &= \; u_0 \; , \\
		U^0_{n+1} \; &= \; G^0_{n+1} \; := \; G^0(U^0_n) \; , \; n \; = \; 0,1,\dots,N_t-1 \; .
	\end{split}
\end{equation}
Here, the superscript indicates the initial iteration for the iteration process. With the initial values provided in every time slice the fine solver $F$ can be executed in parallel in the k-th iteration by:
\begin{equation}
	F^k_{n+1} \; = \; F(U^k_n) \; , \; n \; = \; 0,1,\dots,N_t-1 \; .
\end{equation}
The initial value for all iterations at $t=0$ is given by $U^k_0 = u_0$. With the computationally expensive solver executed in parallel a serial correction is started with the coarse solver $G$. The differences $F^k_{n+1} - G^k_{n+1}$ for every time slice are added to the new serial execution of $G$ as updates:
\begin{equation}
	\begin{split}
		U^{k+1}_0 \; &= \; u_0 \; , \\
		G^{k+1}_{n+1} \; &= \; G(U^{k+1}_n) \; , \\
		U^{k+1}_{n+1} \; &= \; G^{k+1}_{n+1} \; + \; F^k_{n+1} \; - \; G^k_{n+1} \; .
	\end{split}
	\label{EQ:PAR}
\end{equation}
The algorithm is executed until iteration $K$, when a stopping criterion $U^{k+1}_{n+1} - U^k_{n+1} \leq \varepsilon$ is fulfilled. Iterating for $K=N_t$ recomputes the reference solution and therefore, Parareal necessarily converges. On the contrary, in order to generate speedups the algorithm is required to converge with $K \ll N_t$ iterations. \\

The Parareal algorithm presented above allows for different time integration methods or time step sizes as choices for the coarse and fine propagator. In the context of FESOM2 these options are very limited, since the code only offers one time integration method. A variation of time step sizes to distinguish coarse and fast solver was attempted in \cite{philippi2022parareal}. It was observed that the impact of temporal resolution on the relevant diagnostic variables was negligible. In order to define appropriate fine and coarse methods, we decided to execute FESOM2 on two meshes of high and low resolution, respectively. \\
The micro-macro variant to the Parareal algorithm has been introduced and tested for an energy balance climate model in \cite{Slawig2018}. The correction procedure is split into two sub-steps, the first on the macro (low resolution) level and the finalizing step on the micro (high resolution) level:
\begin{equation}
	\begin{split}
		\tilde{U}^{k+1}_{n+1} \; &= \; G^{k+1}_{n+1} \; + \; \mathbf{R} (F^k_{n+1}) \; - \; G^k_{n+1} \; , \\
		U^{k+1}_{n+1} \; &= \; \mathbf{L} (\tilde{U}^{k+1}_{n+1}) \; + \; F^k_{n+1} \; - \; \mathbf{L} (\mathbf{R} (F^k_{n+1})) \; .
	\end{split}
	\label{EQ:MMPAR}
\end{equation}
In addition to the two meshes, interpolation methods must be used to transfer data between the two resolutions in form of a lifting operator $\mathbf{L}(\cdot)$ and restriction operator $\mathbf{R}(\cdot)$. The update on the microscopic level in Eq.\ref{EQ:MMPAR} contains the terms $F^k_{n+1} \; - \; \mathbf{L}$, which are supposed to add fine scale information to the iterative solution $U^{k+1}_{n+1}$ that is lost during the restriction process to the macroscopic level.

\subsection{A priori Speedup Estimate}

The speedup estimate presented in this section is assumed to indicate a priori whether a run-time reduction can be achieved with the respective admissible amount of iterations $K$. The speedup estimate in dependency of the iterations $K$, the run-time ratio $m$ and the amount of time slices $N_t$ is given by:
\begin{equation}
	S_K \; = \; \min \left[ \frac{m}{K+1} , \frac{N_t}{K} \right] \; , \; m \; = \; \frac{\tau_F}{\tau_C} \; .
	\label{EQ:SpeedUP}
\end{equation}
The run-time ratio $m$ is defined by the estimated wall-clock times of the fine and coarse propagator $\tau_F$ and $\tau_C$ for one time slice. The estimate in Eq.\ref{EQ:SpeedUP} requires the run-time ratio to be $m \gg 1$ in order to generate substantial speedups. It must be made clear that the a priori assessment cannot be more than an indicator. Accurate runtime reductions should in any case be estimated by wall-clock time measurements to account for overhead during the execution of Parareal. 

\section{FESOM 2}

In this study, version 2.0 of the Finite-volumE Sea ice-Ocean circulation Model (FESOM2) is chosen to evaluate the parallelization-in-time concept of Parareal. It is formulated on unstructured finite volume meshes. The major advantage over its Finite-Element-Method based predecessor FESOM 1.4 is a computational efficiency comparable to ocean models based on structured meshes. The ocean mesh is split in horizontal planes that are aligned in vertical direction throughout the ocean. By doing so, the numerical efficiency and scalability for parallel computation on HPCs was significantly improved, see \cite{FESOM2performance}. 

FESOM2 provides linear and non-linear free-surface models, turbulent eddy modeling by the Gent-McWilliams parameterization and an isoneutral Redi diffusion. For a detailed description of the solver and the implementation of all its numerical models, see \cite{FESOM2oce}, \cite{FESOM2standardSim}. All variables of ocean and sea ice are intertwined with each other and need to be investigated as an entity to understand the complex physical behavior of the system. However, a evaluation of Parareal with respect to all variables would go beyond the scope of this study. And hence, during the course of the numerical experiments the focus is set on the investigation of diagnostic variables. FESOM2 provides the time-averages of the zonal velocity $u$, vertical velocity $w$, temperature $T$ and salinity $S$, which are computed during runtime and provided at the end of the simulation. We figured, the feasibility of Parareal with spatial coarsening for FESOM2 as a first assessment is best estimated by these diagnostics.  Evaluation of the ocean's global annual surface temperature and the meridional overturning circulation derived from the velocity components are important diagnostic variables in maritime geoscience. In case of interest in the sea-ice model the reader is referred to \cite{FESOM2ice},\cite{FESOM2icepack}.

\section{Mesh Generation}

\subsection{Geometry of FESOM2 meshes}

\begin{figure}[H]
	\centering
	\includegraphics{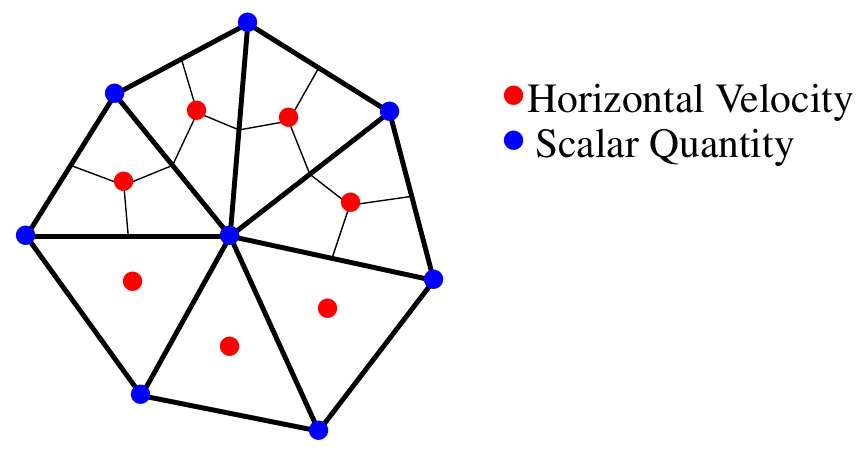}
	\caption{Illustration of variable placement on the horizontal plane. The horizontal velocities are located at the cell centers (red) and scalar quantities at the nodes (blue).}
	\label{FIG:FVMcell}
\end{figure}

Before diving into the mesh generation process an introduction to FESOM2 meshes is given. The mesh is divided in horizontal layers and vertical levels. The ocean surface is modeled by a triangular mesh and the depth created by using multiple layers. This creates a creates a triangular column from the ocean surface to the bottom, which is divided into prisms by the different horizontal levels. The finite volume prisms for vector valued variables are based on the surface triangles, while the scalar quantities are discretized by median-dual control volumes, see Fig. \ref{FIG:FVMcell}. The terms vertex and node will be used as synonyms for the remainder of this document. Scalar quantities, like temperature and salinity, are placed at the nodes (blue) and vector valued variables at the cell centers (red), compare Fig. \ref{FIG:FVMcell}. \\

\begin{figure}[H]
	\centering
	\includegraphics{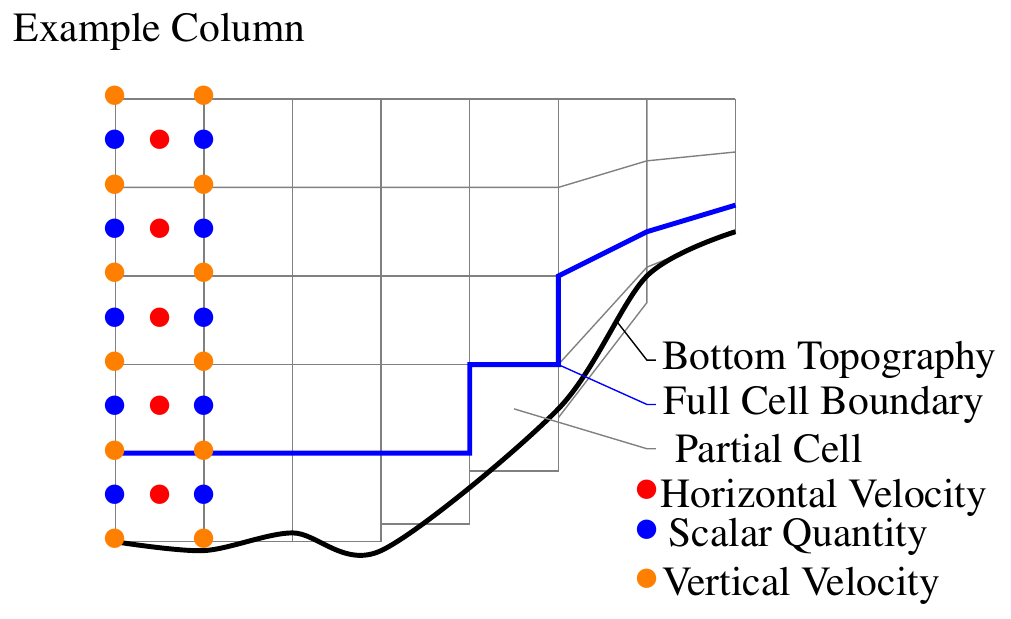}
	\caption{Illustration of the vertical discretization of the ocean. The staggered placement of variables is demonstrated for one top-to-bottom column. FESOM2 allows for use of partial cells and variation in thickness and amount of layers by applying the Arbitrary Lagrangian Euler (ALE) method.}
	\label{FIG:FVMvertical}
\end{figure}

In the vertical discretization scheme the placement is staggered, as shown in Fig. \ref{FIG:FVMvertical}. The horizontal velocities and scalars are placed at mid-level, while the vertical velocity for interaction between layers is stored at full-level. The staggered placement of vertical velocities at the nodes and the horizontal variables at mid-level corresponds to the arrangement of an Arakawa B-grid \cite{FESOM2oce}.
FESOM2 applies the Arbitrary Lagrangian Euler (ALE) method for the vertical coordinates while the horizontal layers remain fixed, allowing the ocean surface and the vertical interfaces to move. Further, the levels follow the bottom topography and the amount of horizontal layers can vary, which lead the developers of FESOM2 to introduce partial cells, as depicted in the right part of Fig. \ref{FIG:FVMvertical}. The advantage of this full nonlinear free-surface model is that fresh water fluxes, by e.g. rain or ice freezing and melting, can be added to the surface layer. By doing so, the salinity at the surface changes with the increase of surface prisms volume. FESOM2 offers two options: \verb|zlevel| and \verb|zstar|. The \verb|zstar| option allows only the surface layer to move, keeping all cells and interfaces below fixed. \verb|zstar| provides the possibility to move all layers except for the bottom level, that remains predefined in space. 

\begin{figure}[H]
	\centering
	\includegraphics{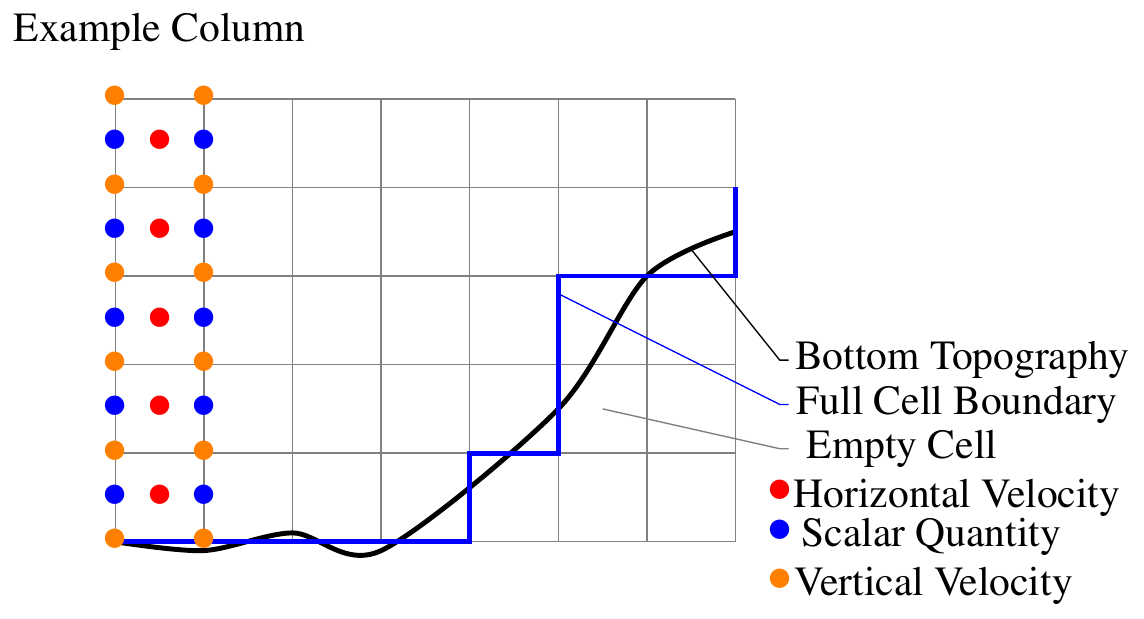}
	\caption{Illustration of the vertical discretization with the linear free-surface ({\psverb+linfs+}) approach. The layer thickness and the volumes of the cells remain constant during the simulation. Cells and nodes outside the boundary (blue line) are part of the mesh but neglected during computations.}
	\label{FIG:FVMvertSim}
\end{figure}

Beside these two options the user is given the choice to run simulations without the ALE approach by keeping the mesh a priori fixed in space, as shown in Fig. \ref{FIG:FVMvertSim}. In the linear free-surface (\verb|linfs|) approach all prism volumes and faces are constant. With no adaption of layer thicknesses the surface requires virtual salinity fluxes for the boundary condition to take the impact of fresh water fluxes into account. Nevertheless, this method means less computational effort and allows for a less complex approach in terms of grid modification. With respect to the feasibility of creating a new mesh, as discussed in the upcoming section, the \verb|linfs| method offers itself as an suitable way to keep the process as simple as possible.

\subsection{On the construction of a new mesh}

Before discussing the step by step approach for the construction of a new domain, it should be noted that with each new mesh evaluation is necessary. Hence, it makes sense to use an already existing mesh provided by the FESOM2 community in order to refine or coarsen into a new mesh. Altering the resolution of an unstructured mesh should at best not affect the quality. To estimate the quality and suitability of a mesh the skewness of the prisms was used in this work. For the prisms in the horizontal layers of the mesh the triangular faces have to examined for skewness, to which the estimate in Eq.\ref{EQ:SKEW} will refer to throughout this chapter. 
Skewness decsribes the grid distortion with respect to the angles between grid lines. It is defined as:
\begin{equation}
	SKEW \; = \; \max \left[ \frac{\theta_\text{max} - \theta_e}{180 - \theta_e} \; , \; \frac{\theta_e - \theta_\text{min} }{\theta_e} \right] \; ,
	\label{EQ:SKEW}
\end{equation} 
with $\theta_\text{max}$ and $\theta_\text{min}$ representing the largest and smallest angle, respectively. $\theta_e = 60$ denotes the angle for an equiangular triangle. The equiangular skew of a prism is $SKEW \in [0,1[$. With increasing skewness the triangle becomes more and more distorted.
\begin{figure}[H]
	\centering
	\includegraphics{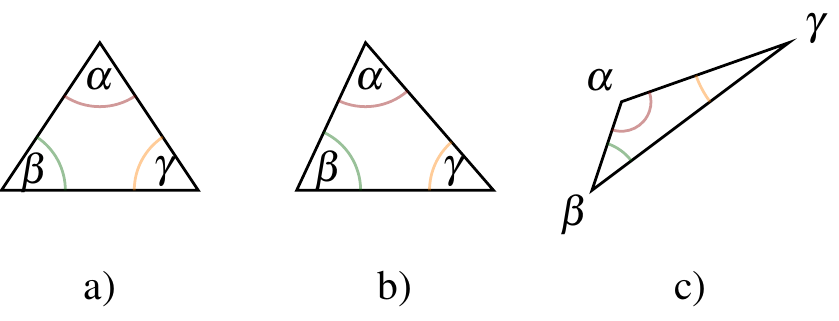}
	\caption{Examples of triangluar skewness: a) triangle with equal angles. b) slightly distorted triangle $SKEW = 0.182$ c) distorted triangle $SKEW = 0.583$}
	\label{FIG:skew}
\end{figure}
In Fig. \ref{FIG:skew} three examples for angular skewness are given: an equilateral triangle, a slightly skewed and a heavily distorted one. Since skewness has significant impact on the numerical stability triangles with no or small distortion are preferred \cite{Ferziger}. 
\begin{figure}[H]
	\centering
	\includegraphics{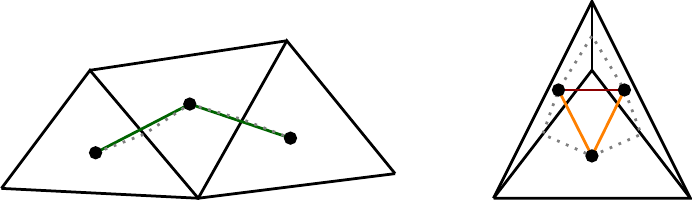}
	\caption{Illustration of skewness in triangles and their impact on the mesh quality. The dotted, gray lines represent connections between cell centers and edge midpoints, colored lines the direct link between centers. Optimally, those lines should overlap or be as close as possible. }
\end{figure}

\begin{figure}[H]
	\centering
	\includegraphics{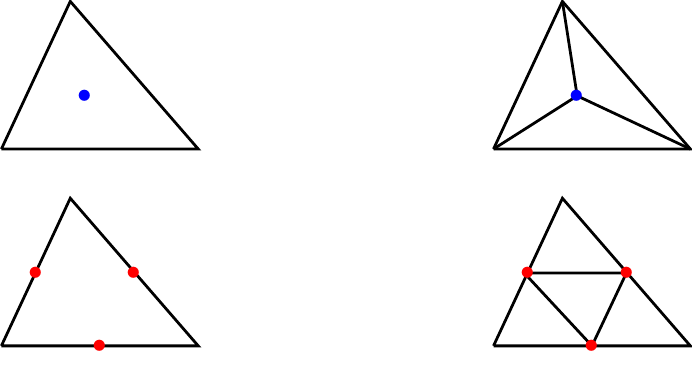}
	\caption{Spatial refinement of unstructured meshes: a) refinement by constructing new triangles with centroids. b) refinement via edge midpoints maintaining the cell skewness.}
	\label{FIG:MESHref}
\end{figure}

Therefore, when altering an existing FESOM2 mesh to obtain a new one, it makes sense to preserve the skewness of the finite volume prisms. For the new mesh in this work, different possibilities for refinement and coarsening of the grid were considered. The refinement strategy concerns adding nodes to the existing mesh that form additional triangles within. In Fig. \ref{FIG:MESHref} two examples for refinement were given: centroid based in Fig. \ref{FIG:MESHref} a) and edge based in Fig. \ref{FIG:MESHref} b). The centroid based approach replaces the existing base triangle with three smaller ones, which are constructed with the bases centroid. By doing so, one increases $\theta_\text{max}$ and reduces $\theta_\text{min}$, and hence the skewness increases. For the edge based approach, three new nodes are added to the base triangle's edge midpoints, over which four new triangles are constructed. Although, this method requires more effort, it can be assured that the angles and skewness are preserved. Both approaches enable refinement by a loop over the existing elements, not affecting respective neighbors. In terms of numerical stability the approach depicted in Fig. \ref{FIG:MESHref} b) should be preferred.

\begin{figure}[H]
	\centering
	\includegraphics{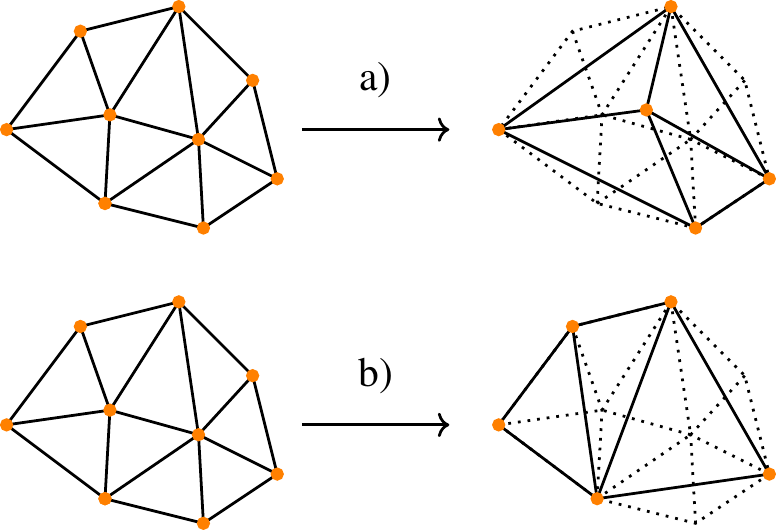}
	\caption{Spatial coarsening of unstructured meshes: a) new triangulation after removing and adding nodes. b) reconnection to triangles by removing nodes only.}
	\label{FIG:Meshcoar}
\end{figure}

The coarsening of unstructured meshes is more complicated. In Fig. \ref{FIG:Meshcoar} depicts two approaches for removing nodes from the mesh to reduce the resolution. In Fig. \ref{FIG:Meshcoar} a) nodes are removed and a new node is inserted in to the domain to construct the coarse mesh. In Fig. \ref{FIG:Meshcoar} b) the amount of nodes is reduced and the remaining nodes are reconnected. Either way, the choice of which nodes to remove, and if necessary to insert, is not unique. The challenge lies in finding optimal nodes for the coarsening process, especially since the triangles' skewness should not be increased significantly. In order to meet these requirements an optimization and triangulation algorithm would be a suitable choice to find those optimal nodes. Without further investigation and knowledge of the base mesh it is not clear, which of the two coarsening concepts would perform better.  
In addition to these considerations to the basic procedure of altering the spatial resolution, the specific requirements in the discretization of the oceans have to be taken into account. Most importantly, one has to preserve the oceans' topology when creating a new computational grid. With identical coast lines and bottom topography for both meshes issues with mass and momentum conservation during the interpolation steps of the Parareal algorithm are avoided. For instance, the PI mesh does not contain the Baltic sea, which is part of the spatial domain of the CORE mesh. Accordingly, an interpolation scheme would fail to generate data for the Baltic sea when mapping from PI to CORE mesh. Which excludes using these meshes for the application of the micro-macro Parareal algorithm. Analogously, this problem does not only occur on the ocean's surface but also in the depth, and so the same considerations are necessary for the bottom topology. Deviations in the meshes ocean depth will result in the same conflict for interpolation attempts. 

For this work the PI mesh was chosen as a basis for refinement by edge midpoints. The advantages are twofold in terms of suitability as a coarse solver and for mesh refinement. Regarding the former, simulations on the PI mesh require few computational resources and produce simulation results in short run-times, which are excellent properties to be applied as a coarse solver. With respect to the latter, it became apparent during the considerations for construction of the new mesh, that refining an existing mesh is more straight forward to accomplish while naturally preserving the topology. In conclusion, three key steps for refining the PI mesh were identified:
\vspace{.25cm}
\begin{itemize}
	\setlength{\itemindent}{1cm}
	\item[1.] Compile all triangles and their edges in arrays.
	\item[2.] Congruent refinement of elements into four triangles of identical skewness.
	\item[3.] Interpolate bottom topography at the added nodes.
\end{itemize}
\vspace{.25cm}
The refined PI mesh will be denoted as FPI mesh for the rest of the report.

\subsection{Mesh Construction}

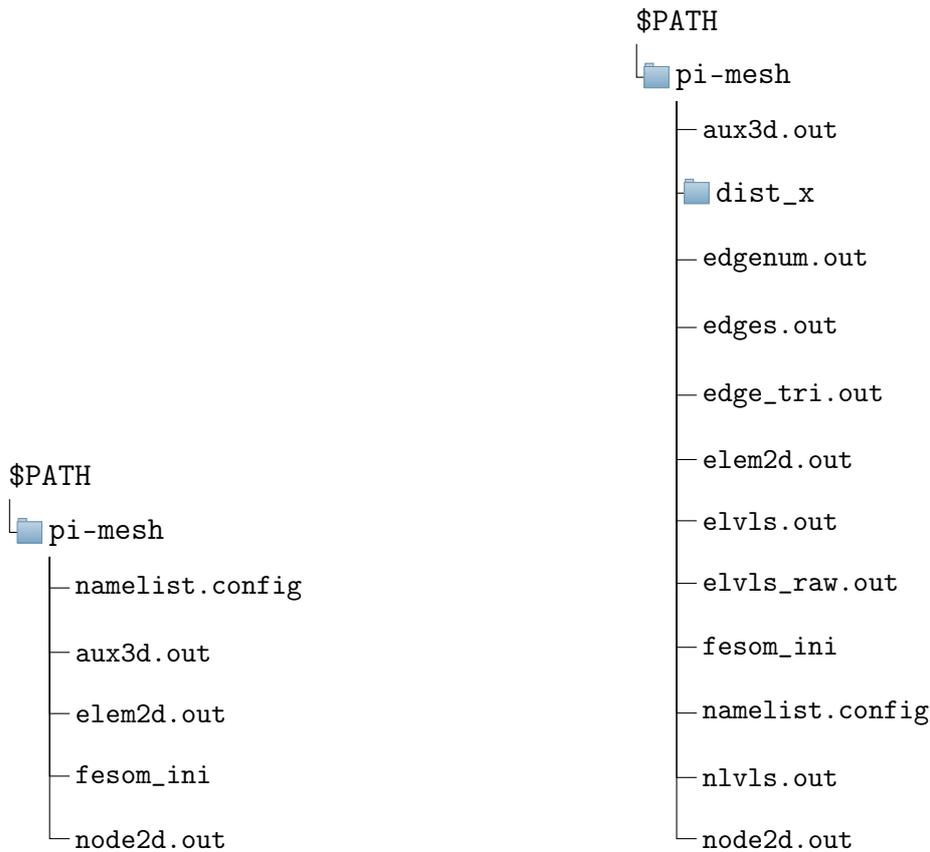
\begin{figure}[H]
	\centering
	\begin{subfigure}{.45\textwidth}
			\begin{forest}
  for tree={
	font=\ttfamily,
	grow'=0,
	child anchor=west,
	parent anchor=south,
	anchor=west,
	calign=first,
	inner xsep=7pt,
	edge path={
		\noexpand\path [draw, \forestoption{edge}]
		(!u.south west) +(7.5pt,0) |- (.child anchor) pic {folder} \forestoption{edge label};
	},
	before typesetting nodes={
		if n=1
		{insert before={[,phantom]}}
		{}
	},
	fit=band,
	before computing xy={l=15pt},
	file/.style={edge path={\noexpand\path [draw, \forestoption{edge}]
			(!u.south west) +(7.5pt,0) |- (.child anchor) \forestoption{edge label};},
		inner xsep=2pt,font=\small\ttfamily
	}
}  
	[\$PATH
		[pi-mesh
		[\textbf{namelist.config}, file]
		[\textbf{aux3d.out}, file]
		[\textbf{elem2d.out}, file]
		[fesom$\_$ini, file]
		[\textbf{node2d.out}, file]
		]
	]
	\end{forest}		
	\end{subfigure}
	\hfill
	\begin{subfigure}{.45\textwidth}
			\tikzsetnextfilename{MeshFolderIni}
\begin{forest}
  for tree={
	font=\ttfamily,
	grow'=0,
	child anchor=west,
	parent anchor=south,
	anchor=west,
	calign=first,
	inner xsep=7pt,
	edge path={
		\noexpand\path [draw, \forestoption{edge}]
		(!u.south west) +(7.5pt,0) |- (.child anchor) pic {folder} \forestoption{edge label};
	},
	before typesetting nodes={
		if n=1
		{insert before={[,phantom]}}
		{}
	},
	fit=band,
	before computing xy={l=15pt},
	file/.style={edge path={\noexpand\path [draw, \forestoption{edge}]
			(!u.south west) +(7.5pt,0) |- (.child anchor) \forestoption{edge label};},
		inner xsep=2pt,font=\small\ttfamily
	}
}  
	[\$PATH
		[pi-mesh
			[\textbf{aux3d.out}, file]
			[dist$\_$x]
			[edgenum.out, file]
			[edges.out, file]
			[edge$\_$tri.out, file]
			[\textbf{elem2d.out}, file]
			[elvls.out, file]
			[elvls$\_$raw.out, file]
			[fesom$\_$ini, file]
			[\textbf{namelist.config}, file]
			[nlvls.out, file]
			[\textbf{node2d.out}, file]
			]
	]
	\end{forest}		
	\end{subfigure}
	\caption{Left: Structure of the mesh folder before executing the FESOM2 mesh builder. Right: The mesh folder after building contains all mesh files required to run FESOM2 simulations. The folder {\psverb+dist$\_$x+} contains the decomposed domain, where x denotes a placeholder (e.g. 2,128,...) for the number of partitions for spatial parallelization.}
	\label{FIG:meshfolder}
\end{figure}

With the identification of a general way to refine meshes the approach used in this work is given in detail. FESOM2 mesh information are stored as ASCII files within a folder. The minimum amount of information to define a mesh is provided by three files: \verb|nod2d.out|, \verb|elem2d.out| and \verb|aux3d.out|. With these files present, the mesh builder of FESOM2 \verb|fesom.ini| can be executed within the working folder to generate additional auxiliary files, see Fig. \ref{FIG:meshfolder}. \\

The created files provide crucial information about triangle edges and the amount of prisms between ocean bottom and surface, but are otherwise of no significance for the mesh creation process. Partitioning settings given in the \verb|namelist.config| file are used for simulations and the mesh building process. FESOM2 provides two parameters for its hierarchical decomposition, \verb|n_levels| and \verb|n_part|. \verb|n_levels| denotes the amount of hierarchical levels, where for each level the number of processors is given by \verb|n_part|. \\
FESOM2 meshes are defined by identical horizontal layers of triangles distinguished by of the amount of vertical layers and the bottom topography. \verb|nod2d.out| contains the horizontal nodes (vertices) that form the triangles (elements) in \verb|elem2d.out|. The amount of layers and the bottom topography at nodes are given in \verb|aux3d.out|.

The \verb|nod2d.out| contains a list of index entries and is followed by coordinates in longitude $\lambda$ and latitude $\phi$. At the end of each line a classification is given, whether the node lies on the boundary or the internal field. For the distinction a "1" is used for the boundary and a "0" for the internal field, respectively. In the first line the total amount of nodes is given. The first lines of \verb|nod2d.out| for the PI mesh are given in Lst. \ref{LST:nod2dpi}.

\begin{lstlisting}[caption={Excerpt of {\psverb+nod2d.out+} for the PI mesh.},label=LST:nod2dpi]
  3140                          # totol amount of nodes
  1 267.4665  84.5252        0 	# indexing starts with "1"
  2 270.9702  84.3700        0	# index long lat classification
  3 268.3207  84.0503        0
  4 263.8669  84.5900        0
  5 264.3323  84.2140        0
  6 261.3389  84.8389        0
             ...
\end{lstlisting}

The  \verb|elem2d.out| files contains the indices from \verb|nod2d.out| and denotes the arrangement of the triangles. As in  \verb|nod2d.out|, the total number of elements is recorded in the first line, see Lst. \ref{LST:elem2pi}. There is no explicit indexing of the triangles, but all variables stored on elements follow the implicit numbering of the file.

\begin{lstlisting}[caption={Excerpt of {\psverb+elem2d.out+} for the PI mesh.},label=LST:elem2pi]
  5839                  # total amount of elements
  1       12        2   # node indices, triangle no. 1
  2       12       10   # node indices, triangle no. 2
  2       10        9   # node indices, triangle no. 3
  3        1        2   # node indices, triangle no. 4
  3        5        1   # node indices, triangle no. 5
  3        9        7   # node indices, triangle no. 6
          ...
\end{lstlisting}
In the \verb|aux3d.out| file the amount of vertical layers, the depth of each level and bottom topography at each node are defined, see Lst. \ref{LST:aux3dpi}. As the file name suggests, the information within contain the instructions to describe the three dimensional domain that consists of two dimensional horizontal layers.

\begin{figure}[H]
	\centering
	\includegraphics{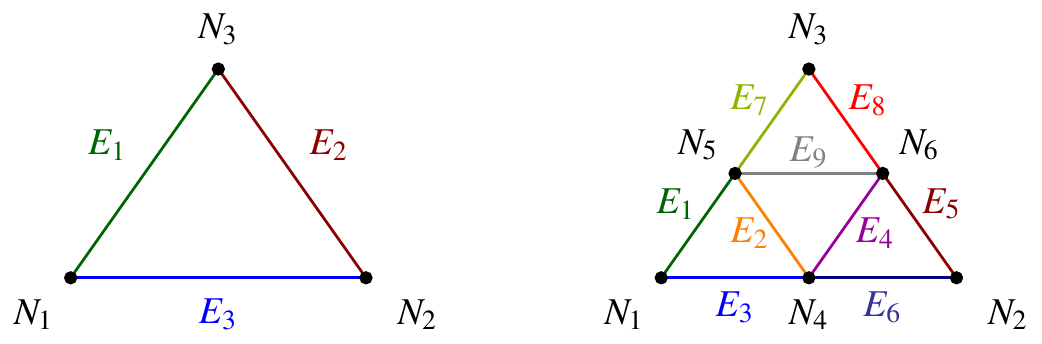}
	\caption{Example mesh consisting of one surface triangle refined with the skew preserving approach. New nodes are inserted at the edge midpoints. Triangle edges are replaced by two new edges forming 4 surface triangles after the refinement process. }
	\label{FIG:EXMESH}
\end{figure}

\begin{lstlisting}[caption=Excerpt of {\psverb+aux3d.out+} for the PI mesh.,label=LST:aux3dpi]
  48        # total amount of horizontal layers
  0.0       # depth horizontal layer 1 (surface)
  -5.0      # depth horizontal layer 2
  -10.0     # depth horizontal layer 3
  ...
  -6000.0   # depth horizontal layer 47
  -6250.0   # depth horizontal layer 48
  -672      # bottom depth at node 1
  -534      # bottom depth at node 2
  -621      # bottom depth at node 3
  -744      # bottom depth at node 4
  ...
\end{lstlisting}
From the ocean depth defined at each node the \verb|fesom.ini| tool derives the amount vertical layers for internal cells (within the ocean) at the respective location. In Fig.\ref{FIG:EXMESH} an example mesh consisting of one surface triangle is shown to illustrate the the refinement process. 

In Lst.\ref{LST:remesh} the \verb|nod2d.out| files before and after the refinement of the example in Fig.\ref{FIG:EXMESH} are given. The edge midpoints of the old mesh are simply attached to the already existing nodes. 
\begin{lstlisting}[caption={\psverb+nod2d.out+} for the example triangle.,label=LST:remesh]
	nod2d_example.out           new_nod2d_example.out
	3                           6
	1 x1 y1 class1              1 x1 y1 class1
	2 x2 y2 class2              2 x2 y2 class2
	3 x3 y3 class3              3 x3 y3 class3
	                            4 x4 y4 class4
	                            5 x5 y5 class5
	                            6 x6 y6 class6
\end{lstlisting}
The definition of elements in \verb|elem2d.out| is entirely replaced by the new elements, compare \ref{LST:remeshELEM}. The orientation of the triangles has to be counter-clockwise. 
\begin{lstlisting}[caption={\psverb+elem2d.out+} for the example refinement.,label=LST:remeshELEM]
	elem2d_example.out          new_elem2d_example.out
	1                           4
	1 2 3                       1 4 5
	                            4 2 6
	                            4 6 5
	                            5 6 3
\end{lstlisting}

\begin{figure}[H]
	\centering
	\includegraphics{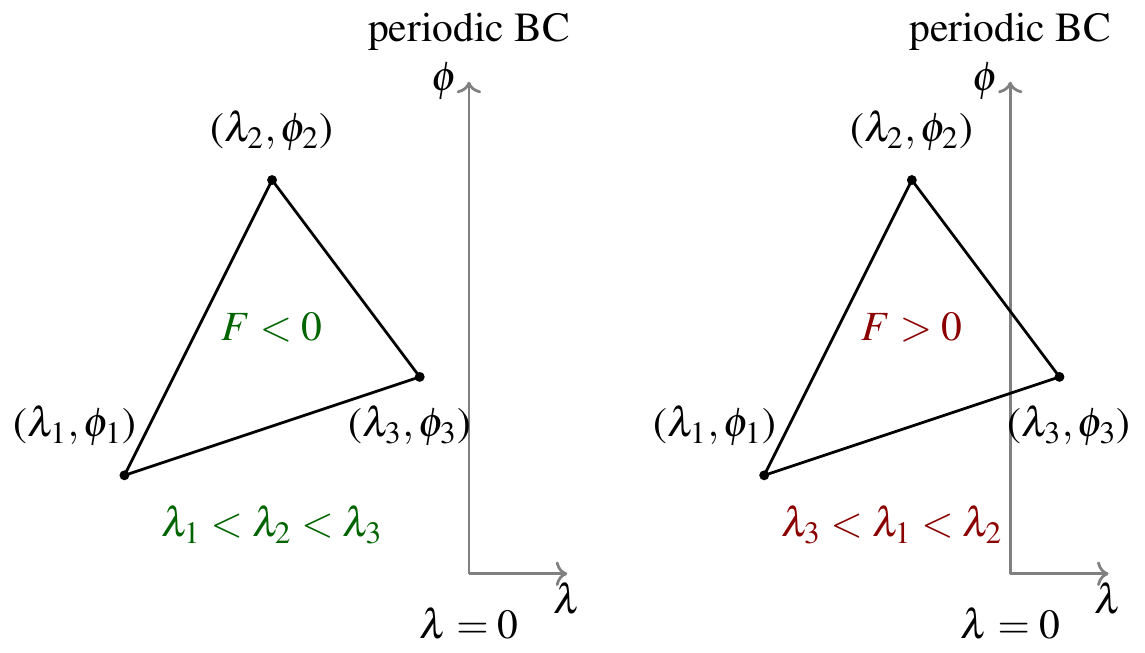}
	\caption{Change in the sign of the orientation of the surface triangles across periodic boundaries. When the change occurs, the $\lambda_i$ coordinates have to be shifted along the $\lambda$-axis such that $F<0$.}
	\label{FIG:triOrient}
\end{figure}

\begin{figure}[H]
	\centering
	\begin{subfigure}[b]{.49\textwidth}
			\includegraphics{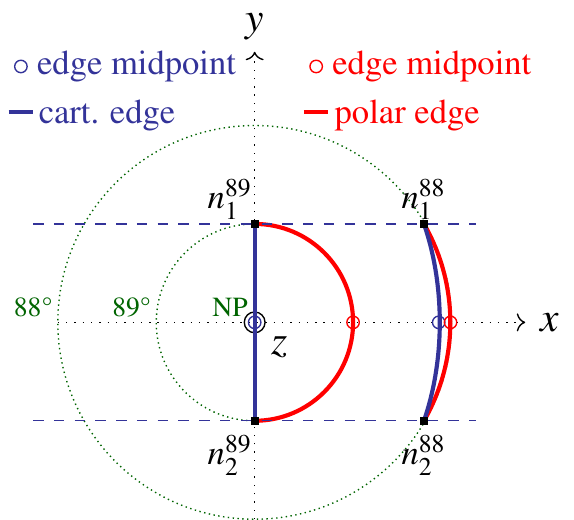}
			\caption{Example edges and nodes.}
		\end{subfigure}
	\hfill
	\begin{subfigure}[b]{.49\textwidth}
			\includegraphics{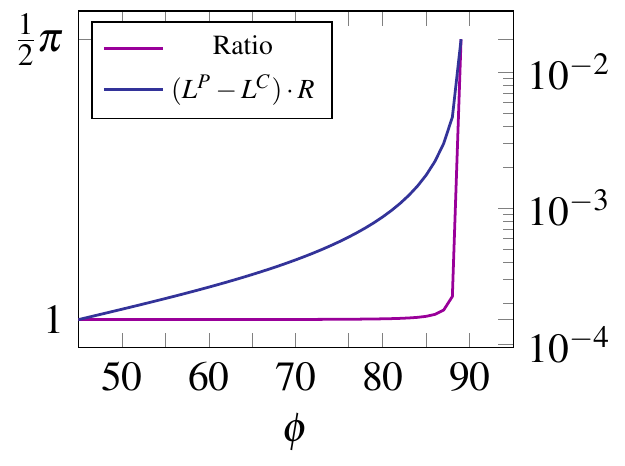}
			\caption{Ratio in dependency of $\phi$.}
		\end{subfigure}
	\caption{Illustration of the difference in edge lengths when interpolating in polar and Cartesian coordinates. a) The Cartesian edge length is chosen to be constant and independent of $\phi$.  b) The length ratio with respect to $\phi$. The ratio approaches $r \rightarrow \frac{1}{2} \pi$ with $\phi \rightarrow 90$° and equals 1 with $\phi=0$ at the equator. In blue is the difference in segment length for both coordinates on a unit sphere given. }
	\label{FIG:PoleComp}
\end{figure}

\begin{lstlisting}[caption={\psverb+aux3d.out+} for the example refinement.,label=LST:remeshAUX]
	aux3d_example.out           new_aux3d_example.out	                            
	48                          48     # amount layers
	0.0                         0.0    # depth layer 1 (surface)
	...                         ...
	-6250                       -6250  # depth layer 48
	-672                        -672   # depth at node 1
	-534                        -534   # depth at node 2
	-621                        -621   # depth at node 3
	                            -603   # depth at node 4
	                            -646.5 # depth at node 5
	                            -577.5 # depth at node 6
\end{lstlisting}
The amount of layers in the refinement process was not altered and therefore, the 49 lines in Lst.\ref{LST:remeshAUX} remain unchanged. In dependency of the sequence of listed nodes in \verb|nod2d.out| the bottom depth is defined. With the new files generated the \verb|fesom.ini| can build the mesh and perform the requested domain decomposition for FESOM2 to work with. \\
There were issues during the refinement of the PI mesh, when the computation of edge midpoints is carried out for triangles crossing the periodic boundary and in the vicinity of the North Pole. In order to identify triangles across the periodic boundary, the formula for computing the triangle area is used:
\begin{equation}
\begin{split} 
F \; &= \; \frac{1}{2} \left( (\lambda_2-\lambda_1)(\phi_3-\phi_1) \; - \; (\lambda_3-\lambda_1)(\phi_2-\phi_1) \right) \; , \\
N_1 \; &= \; (\lambda_1 , \phi_1) \; , \\
N_2 \; &= \; (\lambda_2 , \phi_2 ) \; , \\ 
N_3 \; &= \; (\lambda_3 , \phi_3 ) \; .
\end{split}
\end{equation}
If $F > 0$ the vertices of the triangle are oriented counterclockwise and if $F < 0$ clockwise, respectively. In FESOM2 meshes the triangles are defined in a clockwise manner. Definition of triangles across the periodic boundary will change the orientation to counterclockwise, see Fig.\ref{FIG:triOrient}.
The computation of edge midpoints for triangles across the periodic boundary has to be modified with respect to the $\lambda$ coordinates. As demonstrated in Fig.\ref{FIG:triOrient} the coordinates need to be shifted along the $\lambda$ axis to avoid errors in computing the midpoints. 
\begin{figure}[H]
	\centering
	\includegraphics{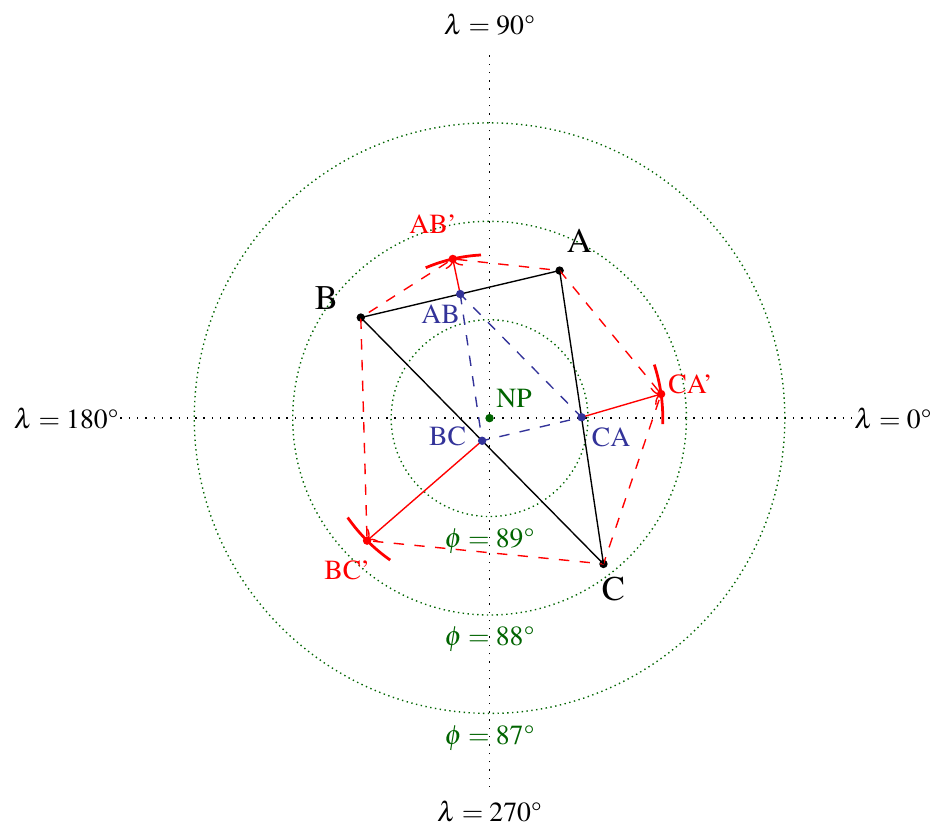}
	\caption{The pole problem for the linear interpolation of edge midpoints. Red lines demonstrate the shift of interpolated nodes in the longitude-latitude representation of coordinates. The edge midpoints interpolated on the spheres surface in Cartesian coordinates are colored in blue.}
	\label{FIG:PoleProblem}
\end{figure}
The nodes of FESOM2 meshes are provided in polar coordinates and linear interpolation of to obtain new nodes in the vicinity of the North Pole will lead to errors in the topology of the prisms. 
All meridians, or lines of constant longitude $\lambda$, merge at the poles. The South Pole is excluded in this case because it is not part the oceans, unlike the North Pole. In Fig.\ref{FIG:PoleProblem} the difference for interpolation in polar and Cartesian coordinates is illustrated. 
For the example triangle the edge midpoints are located outside when interpolated in polar coordinates, leading to a wrong refinement of the parent triangle. We figured two ways to overcome this problem. Transforming the nodes into Cartesian coordinates or rotate the mesh by Euler angles. We chose the latter for our refinement process fo the PI mesh. In older versions of FESOM2 the meshes were rotated such that the coordinate North Pole was not part of the ocean. Hence, the required rotation angle was already known to us. Later, we investigated the differences between the approaches in a test case, shown in Fig.\ref{FIG:PoleComp}. We chose an constant segment across the pole and its correct edge midpoint and compared it to the segment that is generated by linear interpolation in polar coordinates. The ratio of segments approaches $r \rightarrow \frac{1}{2}$ towards the North Pole. The difference in the segment length on a unit sphere is shown in blue in the right panel of Fig.\ref{FIG:PoleComp}. To account for the difference on FESOM2 meshes the value has to be multiplied with the earth' radius, denoting $R=6.317\cdot10^{6}$m. 

\subsection{Interpolation Schemes}

The section's headline suggests multiple possibilities for the interpolation between the PI and FPI mesh. However, we only found the first order conservative scheme for horizontal velocities and the bilinear interpolation for nodal variables successful during first test runs. The first order conservative interpolation scheme is crucial for conserving momentum between the two meshes and was provided by CDO (Climate Data Operators) \cite{CDO}. Any other attempt of interpolating horizontal velocities would fail to have FESOM2 continuing from interpolated restart files. For variables located at nodes, like temperature and salinity, we tested three different interpolation types: the nearest neighbor, bilinear and trilinear method for unstructured grids. The methods were were provided by SciPy \cite{SciPy} and all Python scripts were build upon Numpy \cite{NumPy}. Of these methods, the bilinear scheme was applicable only during test runs and the experimental studies. The nearest neighbor scheme would allow for further serial propagation from interpolated  restart files, but lead to immediate breakdowns during the application of the Parareal algorithm. For the trilinear interpolation we could not interpolate suitable restart files from either mesh to the respective other. Hence, during this section the conservative first order scheme and the bilinear approach are presented. \\

\begin{figure}[H]
	\centering
	\includegraphics{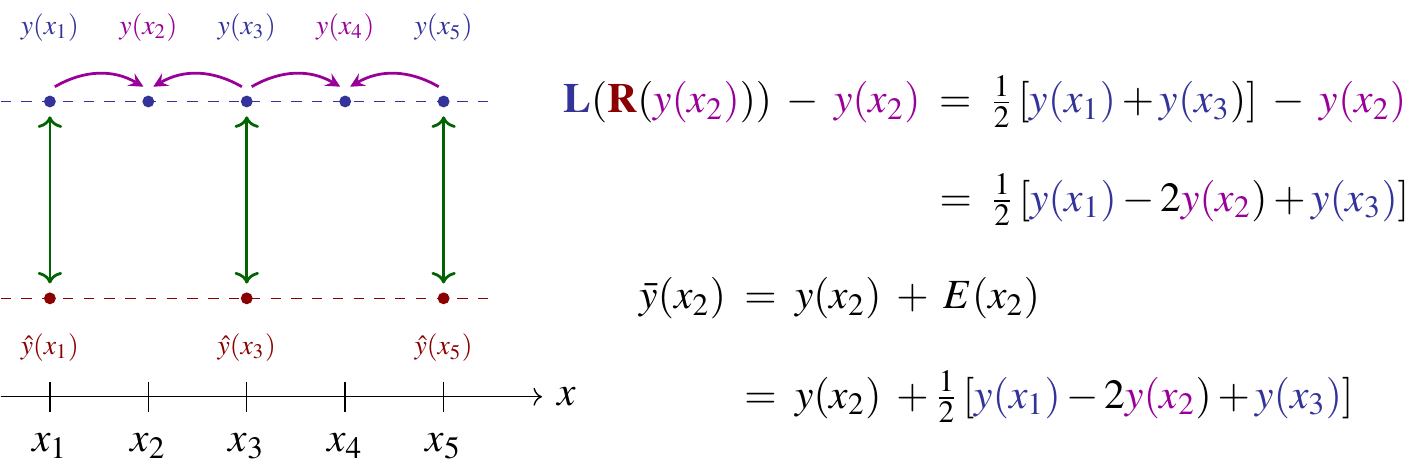}
	\caption{Simplified example for linear interpolation between to meshes. The nodes $x_1,x_3,x_5$ are shared by both example meshes, while $x_2$ and $x_4$ are exclusive to the fine grid only. $y(x_j)$ are the variable values of the fine mesh and $\hat{y}(x_j)$ of the coarse mesh, respectively. The interpolated value $\bar{y}(x_2)$ is given explicitly in this example.}
	\label{FIG:BILIN}
\end{figure}

\begin{figure}[H]
	\centering
	\begin{subfigure}[t]{.48\textwidth}
		\includegraphics{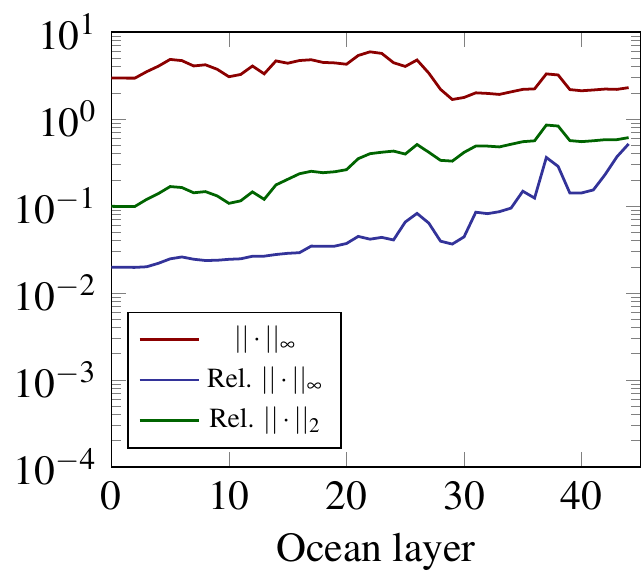}
		\caption{$E^N_I(F)$ for temperature.}
	\end{subfigure}
	\begin{subfigure}[t]{.48\textwidth}
		\includegraphics{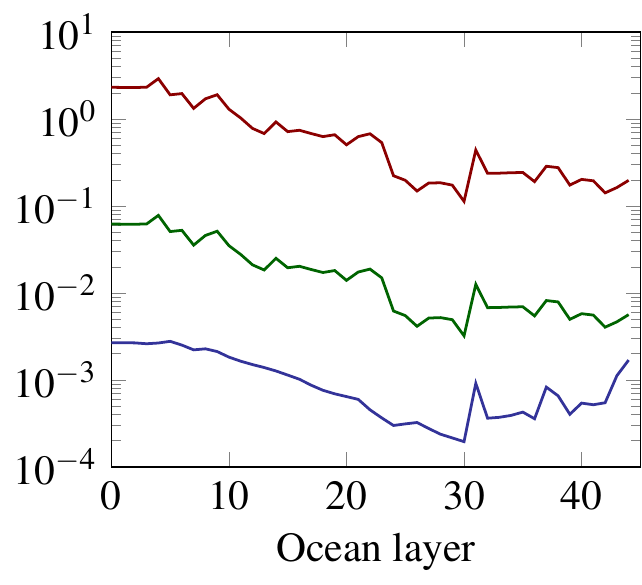}
		\caption{$E^N_I(F)$ for salinity.}
	\end{subfigure}
	\caption{Node interpolation error norms for temperature and salinity.}
	\label{FIG:NINTERROR}
\end{figure}

For the errors $E_I(\cdot)$ in node and element variables when interpolating to the respective other mesh and back are given by:
\begin{equation}
	\begin{split}
		\text{Node interpolation :} &\quad \begin{cases}
			E^N_I(G^k_{n}) \; = \; \mathbf{R} (\mathbf{L} (G^k_{n})) \; - \; G^k_{n} \; &= \; \mathbf{0} \; , \\
			E^N_I(F^k_{n}) \; = \; \mathbf{L} (\mathbf{R} (F^k_{n})) \; - \; F^k_{n} \; &\neq \; \mathbf{0} \; .
		\end{cases} \\
		\text{Element interpolation :} &\quad \begin{cases}
			E^E_I(G^k_{n}) \; = \; \mathbf{R} (\mathbf{L} (G^k_{n})) \; - \; G^k_{n} \; &\neq \; \mathbf{0} \; , \\
			E^E_I(F^k_{n}) \; = \; \mathbf{L} (\mathbf{R} (F^k_{n})) \; - \; F^k_{n} \; &\neq \; \mathbf{0} \; .
		\end{cases}
	\end{split}	
	\label{EQ:IntpCases}
\end{equation}
In Fig.\ref{FIG:BILIN} an illustration of the bilinear interpolation for variables located on nodes is given. Interpolating coarse mesh variables to the fine grid and back (C-F-C) yields no errors $E^N_I(G^k_{n})$ since the PI mesh nodes are shared by both meshes. For $E^N_I(F^k_{n})$ artificial diffusion is added in the process. We evaluated the errors with respect to the ocean layers for temperature and salinity, see Fig.\ref{FIG:NINTERROR}.  

\begin{figure}[H]
	\centering
	\includegraphics{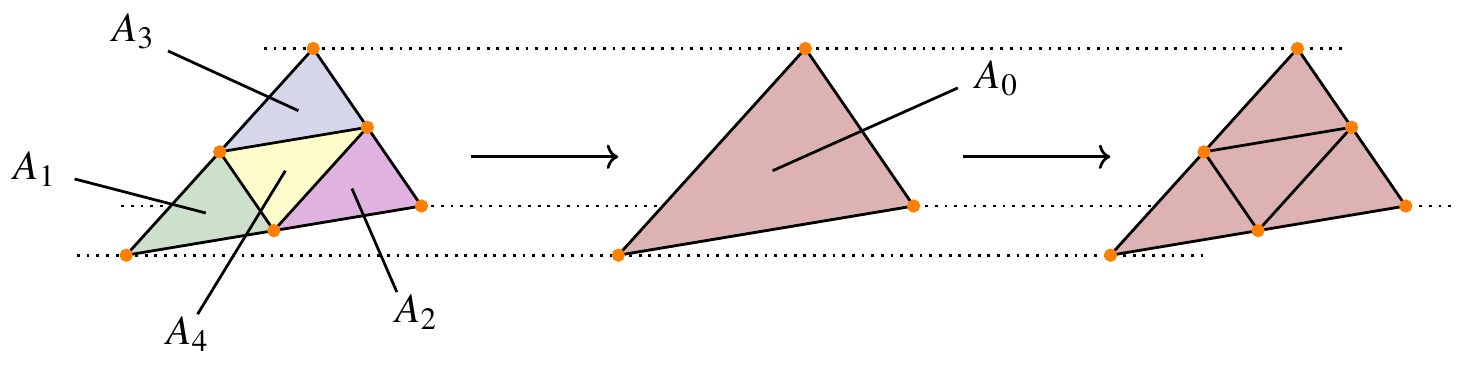}
	\caption{Illustration of the first order conservative interpolation of horizontal velocities located at the centers. Due to the chosen refinement method the union of the triangle areas $A_1,A_2,A_3,A_4$ corresponds to the area $A_0$ of the coarse mesh.}	
	\label{FIG:CONSIN}
\end{figure}

\begin{figure}[H]
	\centering
	\begin{subfigure}[t]{.48\textwidth}
		\centering
		\includegraphics{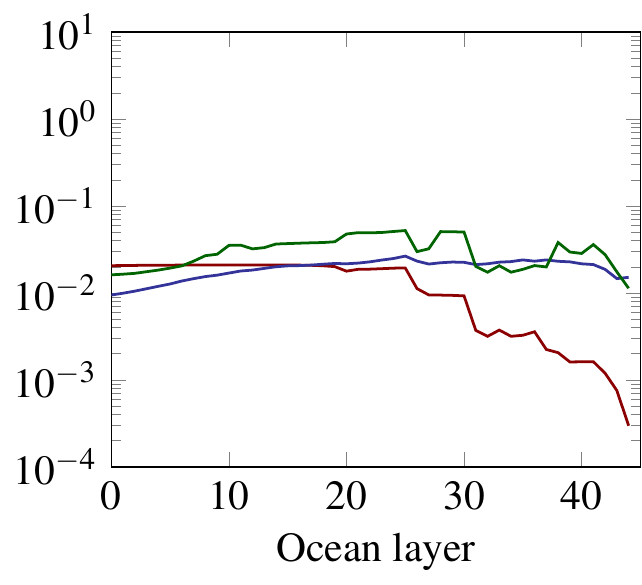}
		\caption[short]{$E^E_I(G)$}
	\end{subfigure}
	\begin{subfigure}[t]{.48\textwidth}
		\centering
		\includegraphics{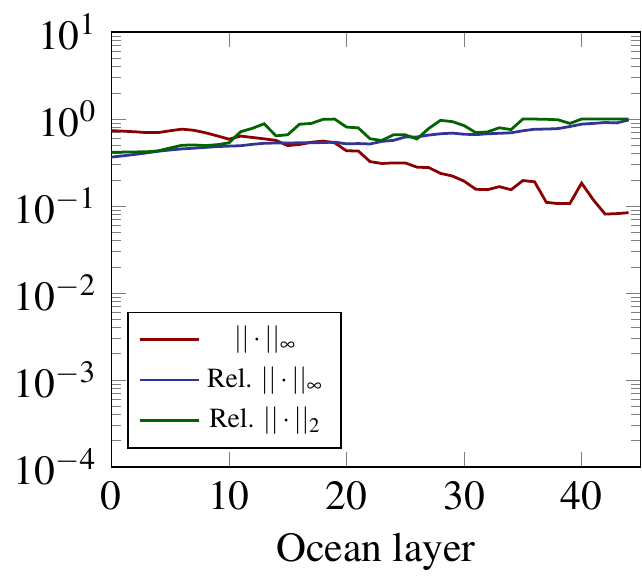}
		\caption[short]{$E^E_I(F)$}
	\end{subfigure}	
	\caption{Element interpolation error norms for the horizontial velocity component $u$.}
	\label{FIG:EINTERROR}
\end{figure}

In Fig.\ref{FIG:CONSIN} the conservative interpolation method of first order is presented. The F-C-F mapping applies an averaging of the fine mesh element values based on their respective volumes. It becomes apparent, that the final interpolation result is filtered. Micro-macro Parareal accounts for this loss of microscopic information during the second update step in Eq.\ref{EQ:MMPAR}. The error estimates in Fig.\ref{FIG:EINTERROR} show that for $E^E_I(G^k_{n})$ (C-F-C) the error is not zero, but stays around $10^{-2}$ for all layers. As expected, the error estimates for $E^E_I(F^k_{n})$ (F-C-F) are significantly impacted by the filtering occurring in the scheme. Nevertheless, we observed that the conservative mapping of the horizontal velocities was crucial to the interpolation. Trying other methods would not succeed in providing restart files, that would allow restarting the computation with FESOM2 on the target mesh. Instead the simulation would suffer immediately from blow-ups in the first time steps.

In Lst.\ref{LST:BILIN} an extract of the bilinear interpolation function within the Python script is given. The function is provided with the node variable to be interpolated and its horizontal layer in the ocean. The latter allows for parallel execution of interpolation in each layer. FESOM2 distinguishes between internal nodes and land nodes. Towards the bottom of the ocean increasingly more nodes are part of the bottom. Although they are stored in the netcdf files, they carry a zero entry. In contrast to the solver, knowing which cells are internal, the interpolation method in Python does not. Thus, in order to avoid including land nodes into the scheme by mistake we provided index masks for each node variable and layer. 

\begin{lstlisting}[language=Python, caption=Extract of intpparscript.py: bilinear interpolation function for node interpolation from coarse to fine mesh.,label=LST:BILIN]
def bilinearInt_tNodeNz1(var,ii):

# i : vertical level
# load layer mask
fname = "$PATH/griddata/coarse/maskNz1/Nz1_layer_"+str(ii)
mask_layer_coar = np.loadtxt(fname)

IDX = tuple(mask_layer_coar.astype(int))

var_coar = SET_in[IDX,ii]

# gather mesh information and variable arrays
z = var_coar
x = pigrid_coar[IDX,0]
y = pigrid_coar[IDX,1]

x = np.append(x,pigrid_coar[IDX,0]-360.0)
y = np.append(y,pigrid_coar[IDX,1])
z = np.append(z,var_coar)

x = np.append(x,pigrid_coar[IDX,0]+360.0)
y = np.append(y,pigrid_coar[IDX,1])
z = np.append(z,var_coar)

# linear interpolator for unstructured meshes
interpolator = LinearNDInterpolator(list(zip(x, y)),z)

fname = "$PATH/griddata/fine/maskNz1/Nz1_layer_"+str(ii)
mask_layer = np.loadtxt(fname)

intp = np.zeros((NODE_F,1))

# interpolation with respect to the masks.
intp[tuple(mask_layer.astype(int)),0] = interpolator(
	pigrid_fine[tuple(mask_layer.astype(int)),0],
	pigrid_fine[tuple(mask_layer.astype(int)),1]
	)

# write interpolated field to coarse restart file
SET_var[:,ii] = intp[:,0]
\end{lstlisting}

The conservative interpolation scheme is provided by CDO. It is a command line tool designed for post-processing data for comparison of simulation results computed on various meshes by different climate model solvers. There is an Python interface that allows for execution of the tool from within the script. Nevertheless, it requires a considerable amount of file handling. The CDO interpolation uses a different mesh description then FESOM2 provides and hence, the element fields have to restructured in a netcdf file, that CDO can use. Although proceeding this way presents a time-consuming process, we considered using it more feasible than programming the method in Python from scratch. The function call in Python is given in Lst.\ref{LST:CONINT}, which can be executed in parallel in each horizontal layer.

\begin{lstlisting}[language=Python, caption=Extract of intpparscript.py: conservative interpolation function for element interpolation from fine to coarse mesh.,label=LST:CONINT]
# gridfiles locations for the cdo remapcon interpolation
intp_grid_coar = "$PATH/griddata/pigrid.txt"
intp_grid_fine = "$PATH/griddata/finepigrid.txt"

def interpolate_elem_FtoC(var,i):
# i :: vertical level 

# initiate cdo library for use in python script
cdo = Cdo()

# SET_in denotes the fine dataset for the variable 'var' of layer 'i'
output = SET_in[:,i]

# write layer data to output file
fname = "$PATH/tmp/FtoC/slice"+args.slice+"/"+str(var)+"/"+str(var)+"_layer_"+str(i)
np.savetxt(fname,output)

# read file into cdo and create netcdf file
foutput = "$PATH/tmp/FtoC/slice"+args.slice+"/intp_"+str(var)+"/fine_"+str(i)+".nc < "+fname
cdo.input(intp_grid_fine,output=foutput,options="-f nc")


fgrid = intp_grid_coar+" -setgrid,"+intp_grid_fine
finput = "$PATH/tmp/FtoC/slice"+args.slice+"/intp_"+str(var)+"/fine_"+str(i)+".nc"
foutput = "$PATH/tmp/FtoC/slice"+args.slice+"/intp_"+str(var)+"/coar_remapcon_"+str(i)+".nc"

# interpolate and write to output file 
cdo.remapcon(fgrid,input = finput,output = foutput)

# read output file 
intp_coar   = netCDF4.Dataset(foutput,'r+')
intp_elem = intp_coar['var1'][0,:]
intp_coar.close()

# make mask check for bottom topography
if i > 3:
 mask = np.loadtxt('$PATH/griddata/coarse/maskElem/mask_'+str(i))

 if mask.size == 0:
  pass
 else:
  for j in range(len(mask)):
   intp_elem[int(mask[j])] = 0.0 
else:
 pass

# write interpolated field to coarse restart file
SET_var[:,i] = intp_elem
\end{lstlisting}

\begin{figure}[H]
	\centering
	\includegraphics[width=.48\textwidth]{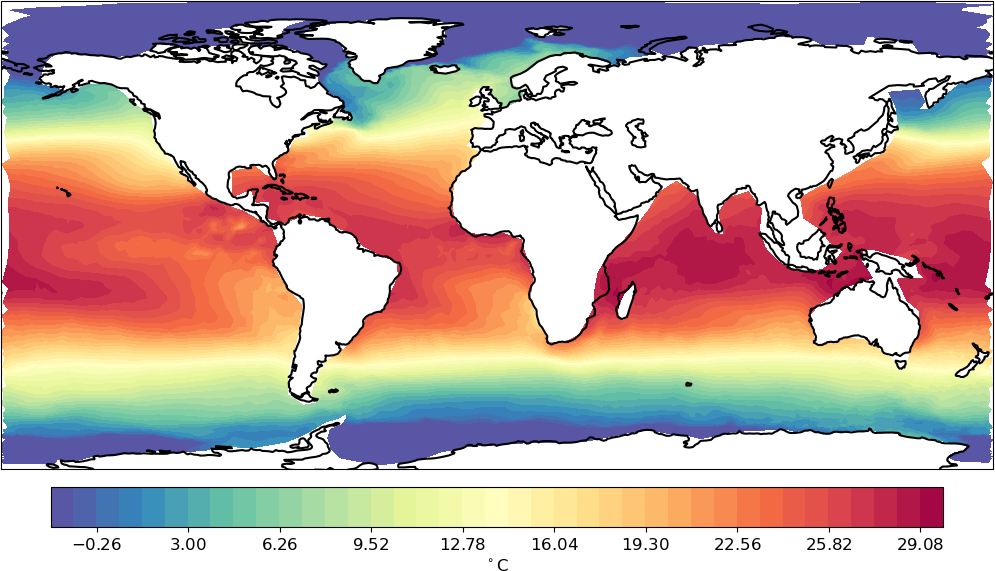}
	\hfill
	\includegraphics[width=.48\textwidth]{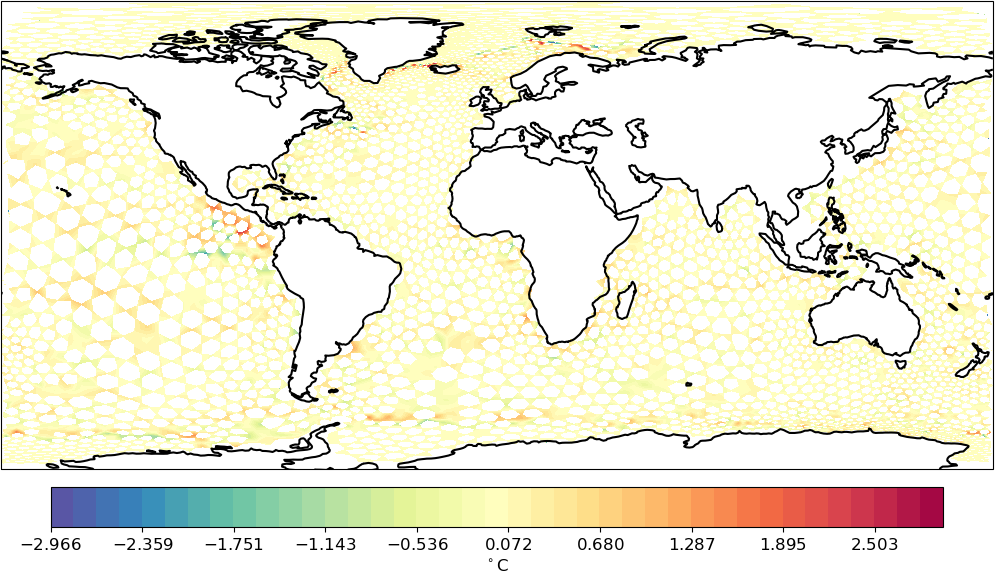}
	\caption{Left: example temperature field on the fine mesh. Right: Difference between reference field and the F-C-F interpolation. The white spots in the figure indicate the shared nodes in both meshes. For those nodes the interpolation is exact and the error is zero.}
	\label{FIG:FCFN}
\end{figure}
\begin{figure}[H]
	\centering
	\includegraphics[width=.48\textwidth]{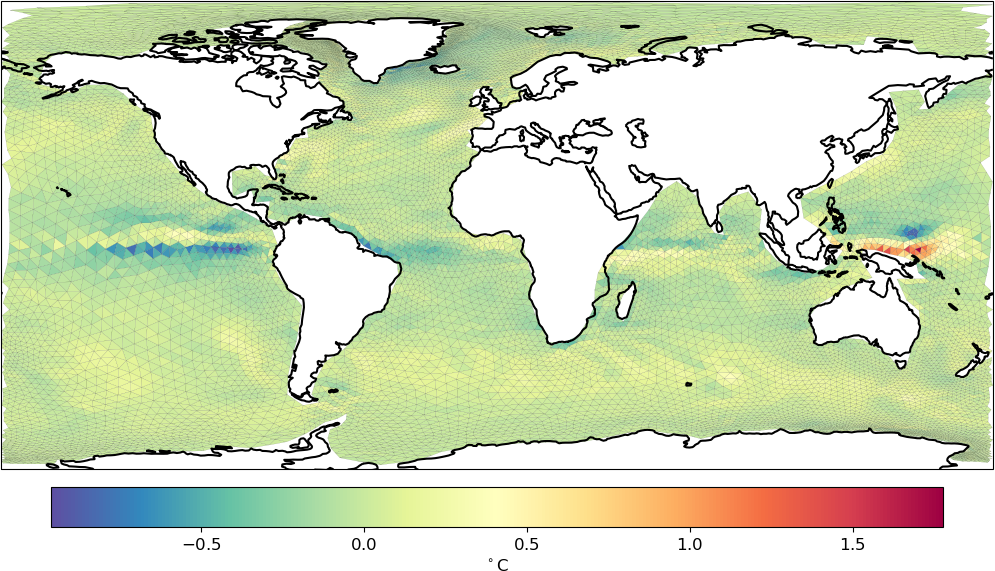}
	\hfill
	\includegraphics[width=.48\textwidth]{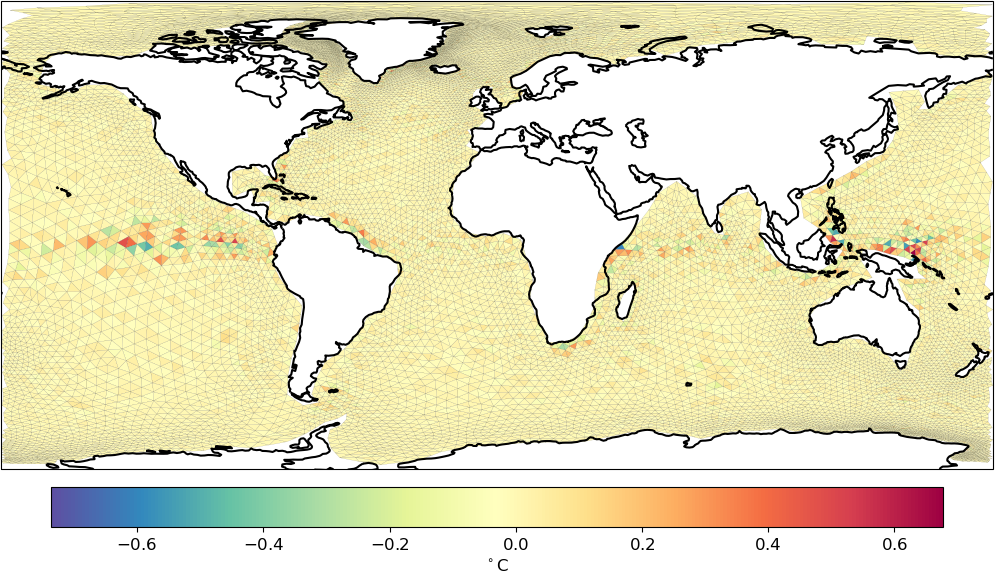}
	\caption{Left: $u$ velocity of the reference field on the fine mesh. Right: Difference between the reference velocity and the F-C-F interpolated field.}
	\label{FIG:FCFE}
\end{figure}

At the end of this section the fields of the fine scale contribution during the second update step in Eq.\ref{EQ:MMPAR} are shown. In Fig.\ref{FIG:FCFN} an example temperature field is given on the FPI mesh. In comparison the $F^k_{n+1} - \mathbf{L} (\mathbf{R} (F^k_{n+1}))$ field is shown in right panel. Both meshes share the nodes of the PI mesh and therefore can be identified by the white spots in the depiction, where the error is zero. The largest errors can be found along the boundaries between warm and cold water fronts at the South and North Pole. 
In Fig.\ref{FIG:FCFE} the $u$ velocity component and its interpolation difference field is shown. Along the equator the $u$ velocity is observed to change its direction with respect north or south. The loss of microscopic dynamics by interpolation in these regions is particularly evident. 

\subsection{The Meshing Tool}

In this section we give a brief introduction to the Python meshing tool and its usage. The main file is shown in Lst.\ref{LST:MTOOL}. The code was created with NumPy \cite{NumPy}. The tool reads in the PI mesh files \verb|nod2d.out|, \verb|elem2d.out| and \verb|aux3d.out|. A class for the mesh is created in which the refinement steps are carried out. In addition to the creation of the FPI mesh, the meshing tool provides gridfiles required for the CDO interpolation of horizontal velocities. As mentioned earlier, the mesh description needs to be changed for CDO and thus, these files are generated along the way of building the new mesh. With completion of the new mesh the \verb|fesom.ini| tool is executed within the output folder of the tool to create auxiliary files for FESOM2 and perform domain decomposition, compare Fig.\ref{FIG:meshfolder}.

\begin{lstlisting}[language=Python, caption=Extract of intpparscript.py: conservative interpolation function for element interpolation from fine to coarse mesh.,label=LST:MTOOL]
#!/usr/bin/env python3

import src.mesher as mesher
import src.netcdf as cdf

def main() -> None:

# paths to the pi-mesh files	
str_elem = "$PATH_TO_ELEM/elem2d.out"
str_node = "$PATH_TO_NODE/nod2d.out"
str_aux  = "$PATH_TO_AUX/aux3d.out"

# read in the mesh files
grid = mesher.convert_str_to_node(str_node)
elements = mesher.convert_str_to_elem(str_elem)
[zlevel, bot_topo] = mesher.convert_str_to_aux(str_aux,grid)

# create mesh class 
pimesh = mesher.MESH(grid,
                     elements,
                     bot_topo,
                     zlevel
                    )
# distinguish nodes into boundary and internal nodes
pimesh.classify_nodes()

# centroids are needed for the auxiliary grid file 
# used in the conservative interpolation of velocities (CDO).
# This step is not necessary for the mesh refinement itself!
pimesh.set_centroids()

# congruent refinement of triangles 
pimesh.refine_congruent()

# refinement of the bottom topography.
# the function allows for three different interpolation types:
#	- nearest neighbor
#	- bilinear 
#	- trilinear 
pimesh.refine_aux3d("linear")

# output of the refined mesh files for the FPI mesh.
pimesh.write_mesh("$PATH_OUTPUT/")

# creating a gridfile that can be utilized by CDO interpolation
# routines. The gridfile contains a different representation of
# the mesh.
pigrid = cdf.DATASET(pimesh.nodes,
                     pimesh.elements,
                     pimesh.centroids
                    )

# output of the PI grid file.
pigrid.write_gridfile("$PATH_TO_OUTPUT/pigrid.txt")

# creating an own mesh class for the FPI mesh.
# The class is needed to produce the gridfile for interpolations
# from PI to FPI and vice versa.
# If desired, the mesh could be further refined.
fineMesh = mesher.MESH(pimesh.nodes_refined,
                       pimesh.elements_refined,
                       pimesh.bottom_refined,
                       zlevel
                      )

fineMesh.set_centroids()

output of the FPI grid file.
finepigrid = cdf.DATASET(fineMesh.nodes,
                         fineMesh.elements,
                         fineMesh.centroids
                        )

finepigrid.write_gridfile("./finepigrid.txt")

if __name__ == "__main__":
main()

\end{lstlisting}

\begin{figure}[H]
	\centering
	\includegraphics[width=\textwidth]{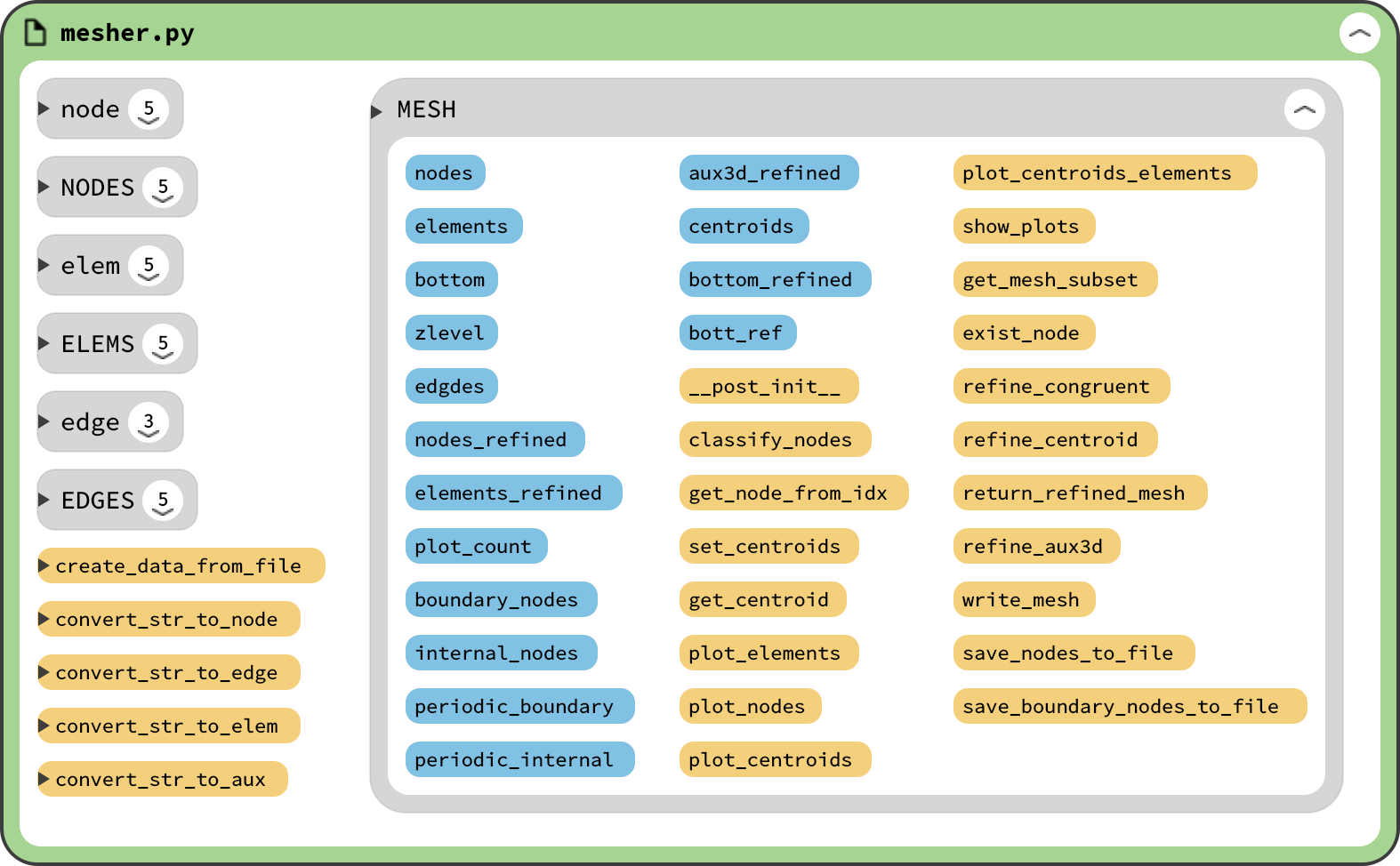}
	\caption{Visualization of the mesher.py file: gray boxes represent classes, blue boxes private members and yellow boxes functions. The MESH class creates is called in pyMesher.py in order to refine the source mesh.}
	\label{FIG:MESHERPY}
\end{figure}

For the sake brevity an overview of the \verb|mesher| class with its functions and members is given in Fig.\ref{FIG:MESHERPY}. In the first attempts to refine a test setting of one triangle, the functions for visualization and auxiliary output files of nodes and elements were implemented in order trace errors. The visualization was made with the open-source tool Sourcetrail, see \url{https://github.com/CoatiSoftware/Sourcetrail}. \\

\begin{figure}[H]
	\centering
	\includegraphics[width=.5\textwidth]{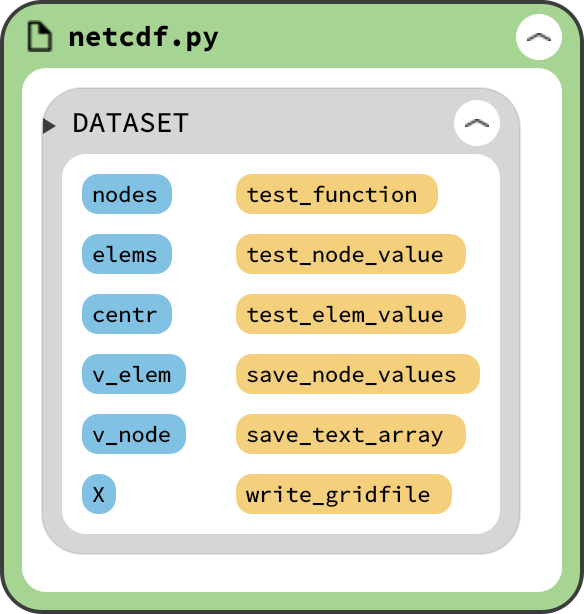}
	\caption{Visualization of the netcdf.py file: gray boxes represent classes, blue boxes private members and yellow boxes functions. The DATASET class is used to generate gridfiles for the interpolation of the velocity fields.}
	\label{FIG:NETCDFPY}
\end{figure}

The \verb|netcdf| class is shown in Fig.\ref{FIG:NETCDFPY}. It is included for the only purpose to generate the gridfiles for the CDO interpolation. It has no effect on the mesh generation, but requires the node and element information from the respective meshes. We considered the integration to the meshing tool reasonable, as it generates the appropriate CDO mesh description along with the creation of a new mesh. 

\begin{figure}[H]
	\centering
	\begin{subfigure}[t]{.49\textwidth}
			\includegraphics[width=\textwidth]{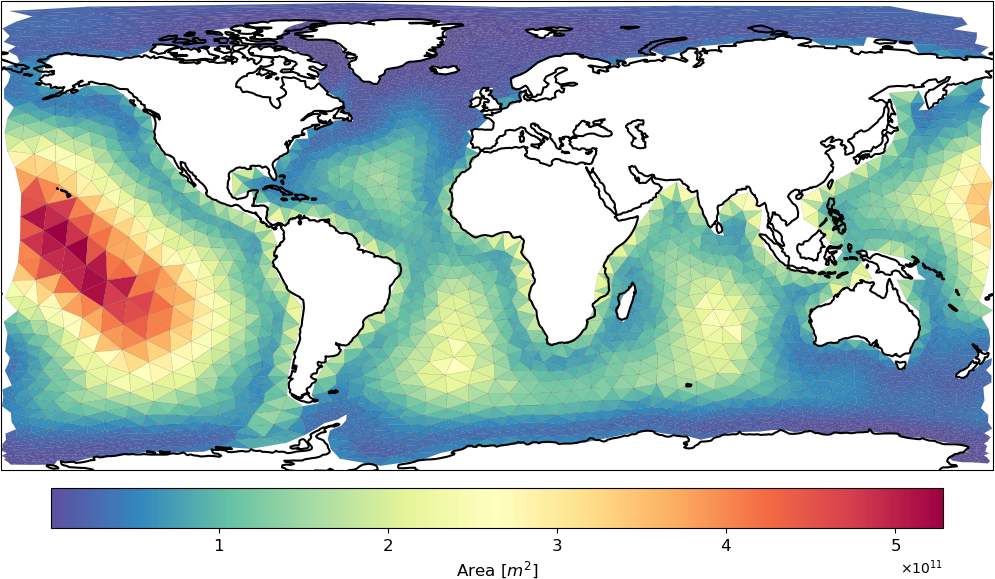}
			\caption{Area of the triangle faces of the PI mesh.}
		\end{subfigure}
	\hfill
	\begin{subfigure}[t]{.49\textwidth}
			\includegraphics[width=\textwidth]{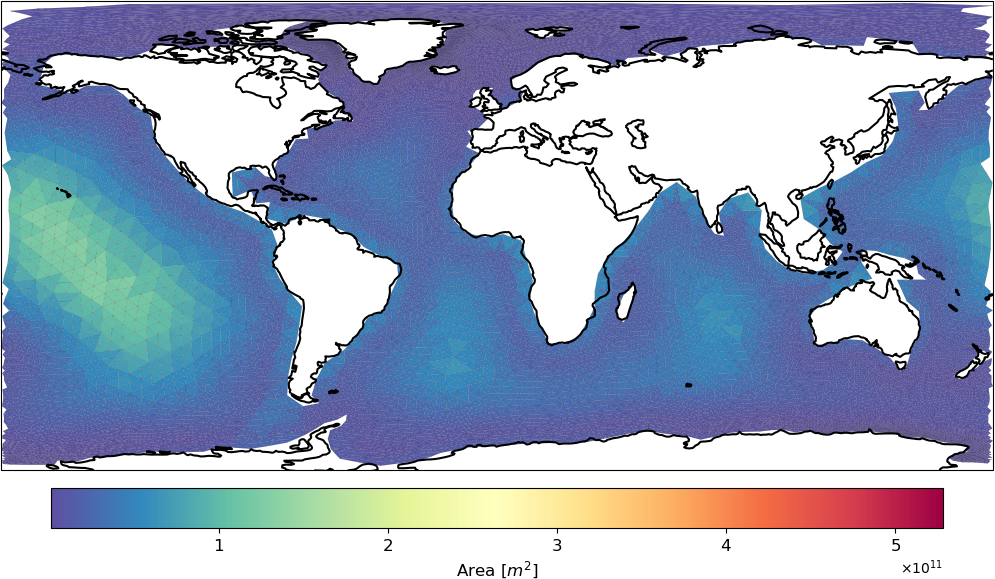}
			\caption{Area of the triangle faces of the FPI mesh.}
		\end{subfigure}
	\caption{Comparison of the PI and FPI resolution by cell area [$m^2$].}
	\label{FIG:COARSEFINE}
\end{figure}

\begin{figure}[H]
	\centering
	\begin{subfigure}[t]{.49\textwidth}
			\includegraphics[width=\textwidth]{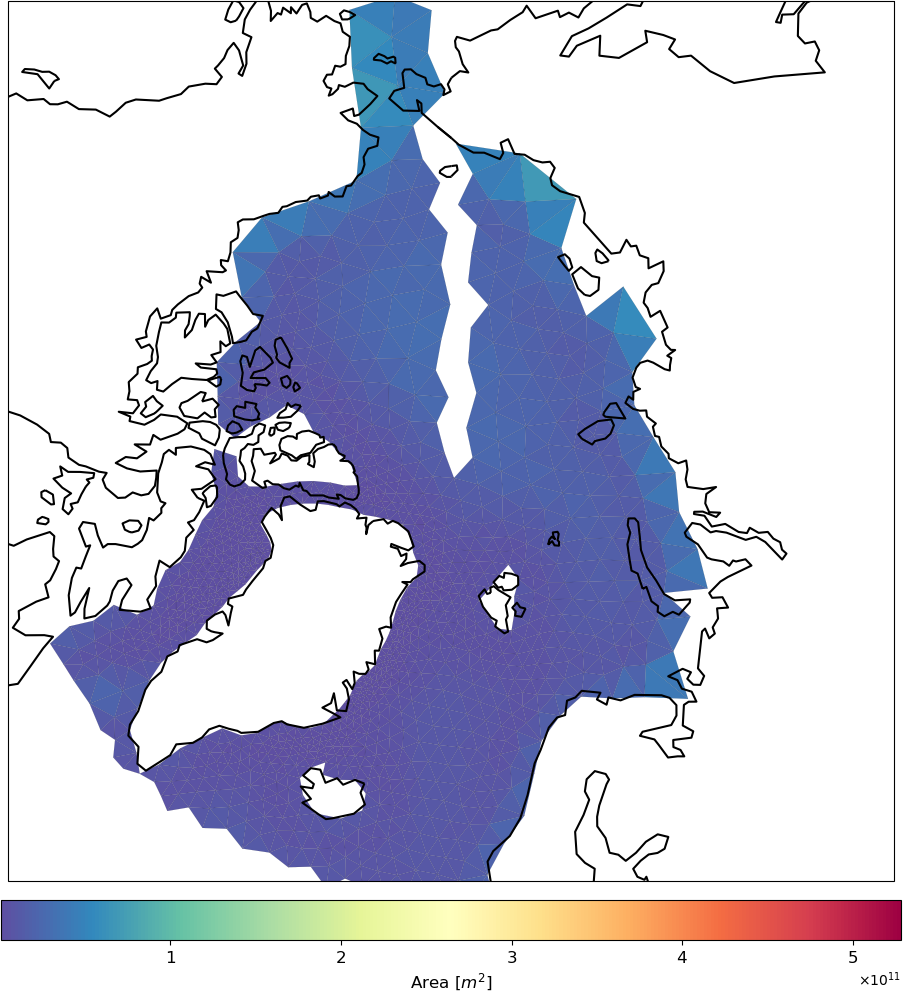}
			\caption{PI mesh}
		\end{subfigure}
	\hfill
	\begin{subfigure}[t]{.49\textwidth}
			\includegraphics[width=\textwidth]{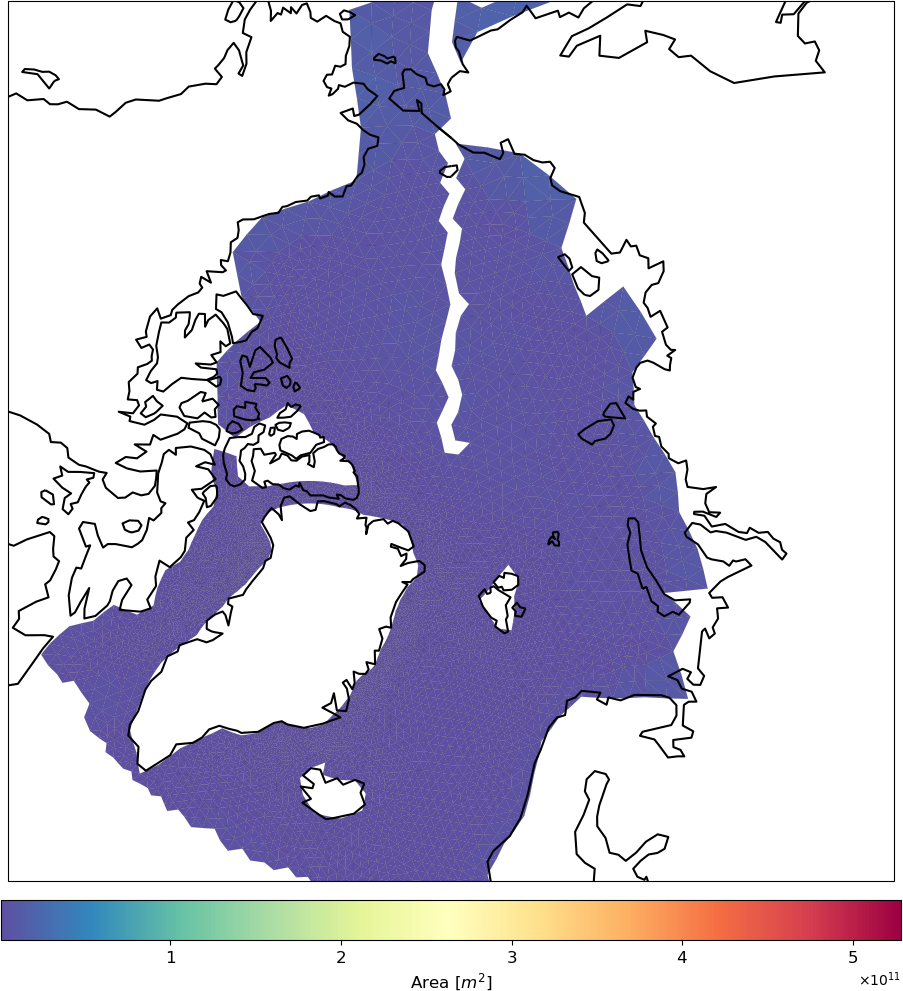}
			\caption{FPI mesh}
		\end{subfigure}
	\caption{Comparison of the PI and FPI resolution at the north pole.}
	\label{FIG:POLECOARSEFINE}
\end{figure}

To conclude this section, we present the resolutions of the PI and FPI mesh in Fig.\ref{FIG:COARSEFINE} and around the North Pole in Fig.\ref{FIG:POLECOARSEFINE}. As demanded, the coast lines and the triangle angles are preserved. In preserving the triangle skewness we observe the smooth gradual transition from large to small surface triangles. The care taken in the creation of the pi mesh could thus be taken over by the proposed refinement method.

\subsection{Background for the Mesh Evaluation}

With the completion of a FESOM2 run the output folder provides up to 43 diagnostic variables. Evaluation of an ocean simulation poses a very demanding and time consuming task. Considering all these variables is beyond the scope of this study and therefore we have chosen to investigate temperature, salinity and velocities. Additionally, we included a comparison of the impact of temporal resolution for both meshes with respect to Global averaged annual temperature at the end of this section. In order to perform an evaluation of the simulation results we used \verb|pyfesom2| provided by the developers of FESOM2. \verb|pyfesom2| is a Python based free software by the AWI to post-process simulation results. In case of further interest in the tool the reader is referred to the \verb|github| repository \url{https://github.com/FESOM/pyfesom2} for detailed functionality and options for data analysis and visualization. In this study we focused on the following diagnostics, which are introduced during this section:
\vspace{.25cm}
\begin{itemize}
	\item[1] Global averaged annual temperature at [0, 500, 0-500] m depth,
	\item[2] Atlantic Meridonal Oveturn Circulation (AMOC),
	\item[3] Atlantic Meridonal Oveturn Circulation (AMOC) at 26.5,
	\item[4] Mean horizontal velocity fields,
	\item[5] Mean salinity fields,
	\item[6] impact of temporal resolutions for PI and FPI mesh.
\end{itemize}
\vspace{.25cm}
In order to access the physical behavior of our fine and coarse settings we computed hindcast simulations with the CORE2 forcing data set, see \cite{LargeCORE2},\cite{core2}. Hindcast simulations rely on data from other models that interact with the ocean. For a complete earth simulation all contributing models (ocean, atmosphere, fresh water fluxes, radiation, bio-chemistry and many more) would have to be coupled and interact with each other. The attempt to unify all models to an earth system model (ESM) poses a highly complex task and is left to state-of-art research efforts. To make this study as feasible as possible, we used CORE2 information to provide the necessary data to approximate an representation of the earths ocean. \\
Since we had no initial values available for our study we decided to perform spin-up simulations. If no initial value is provided for FESOM2 it is derived from a climatology data set. Climatology data sets contain information about the long-time averaged physical conditions of the climate. With these initial values one propagates forward in time until a stable state for the chosen setting is achieved. This process is model and resolution dependent and can take up to several hundred years. Since the CORE2 data set contains forcing information over the period from 1948 to 2006, we were forced to restart the spin-up simulations after each so-called CORE2-cycle of 58 years. In total, we simulated over 600 years for the evaluation part. \\
For the sake of completeness, we compared the results of the PI and FPI spin-ups with a 400 years CORE mesh run. The CORE mesh contains 127k surface nodes and represents a broadly used standard mesh, although there exist several spatially higher resolved meshes.

\subsubsection{Annual Mean Temperature}

The annual mean temperature certainly is the most prominent indicator for climate change. It is commonly evaluated at the oceans surface and in combination with the near-surface atmospheric temperature over land. During the UN Climate Change Conference (COP21) in Paris, on 12 December 2015, the goal to limit the increase in global mean temperature below 2 degrees Celsius was agreed upon by 196 participating countries. At least since then, it became the well-known measure for advancing climate change. \\
During this study we investigated the annual mean temperature at the ocean surface and its vicinity. The diagnostic temperature is averaged annually for each node on the mesh and provided as an own netcdf-file in the output folder. The spatial averaging is performed with \verb|pyfesom2|, where the possibility is given to investigate the changes at 0m (surface) and 500m depth, as well as averaged from 0-500m.

\subsubsection{Atlantic Meridionial Overturn Circulation (AMOC)}

\begin{figure}[H]
	\centering
	\includegraphics{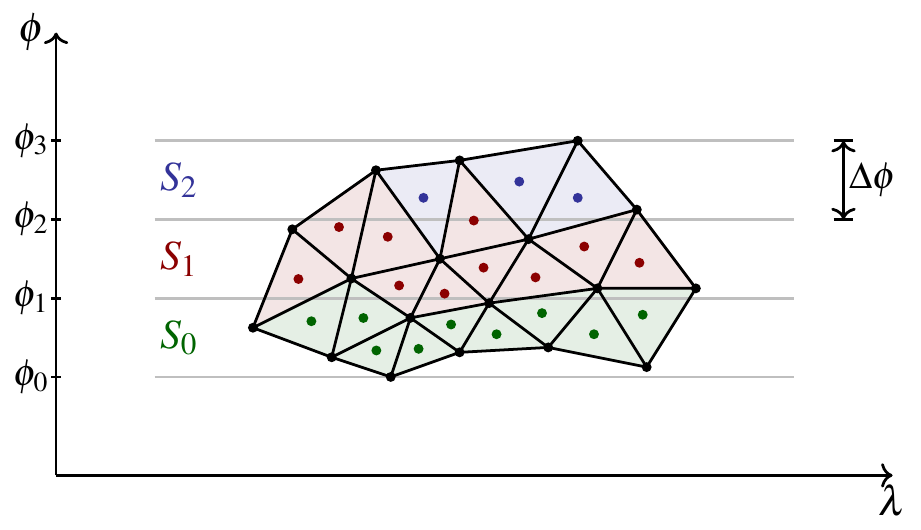}
	\caption{Schematic of the meridional overturn circulation computation. The vertical velocities are interpolated from the nodes to cell centers and subsequently divided into sections along the latitude of length $\Delta \phi$. Triangle faces are colored according to their assigned section, illustrating how the choice of $\Delta \phi$ and the meshes resolution affects the outcome.}
	\label{FIG:AMOCBIN}
\end{figure}

The transport of warm water masses on the Atlantics surface towards the North Pole and cold water along the bottom southwards is defined as the Atlantic Meridionial Overturn Circulation (AMOC). The circulation significantly impacts the climate conditions in Europe and is an important indicator for climate change. \\
According to \cite{FESOM2amoc} there are two ways to compute the meridional overturn circulation for FESOM2 meshes. Both definitions are equivalent due to the assumption of incompressibility for water. With $\nabla \cdot \underline{v} = 0$ and $\underline{v} = (u,v,w)$ the stream function $\Psi(t,\phi)$ can be computed by just taking the vertical velocity $w$ into account:
\begin{equation}
	\Psi(t,\phi) \; = \; \int_{z_{low}}^{z_{up}} R_E d\phi \; \int_{\lambda_W}^{\lambda_E} w d\lambda \; .
\end{equation}
with $R_E$ denoting the earth' radius. To approximate the integral for the stream function three steps are taken: 
\vspace{.25cm}
\begin{itemize}
	\item[1.] Interpolate vertical velocity from nodes to cell centers.
	\item[2.] split the atlantic ocean into latitude bins
	\item[3.] sum all vertical transports for each bin an horizontal layer
\end{itemize}
\vspace{.25cm}
In Fig.\ref{FIG:AMOCBIN} an example for the distribution of cells to latitudinal sections is depicted. The assignment is given by the location of the cell centers. 

The interpolation of node velocities $w_{l,i}$ of an element $e_l$ to its cell center $c_l$ is carried out by:
\begin{equation}
	w(c_l) \; = \; \frac{1}{3} \sum_{i=1}^{3} w_{l,i} \; ,
\end{equation} 
where $l$ denotes the index for the respective element in the FESOM2 mesh.
\begin{figure}[H]
	\centering
	\includegraphics{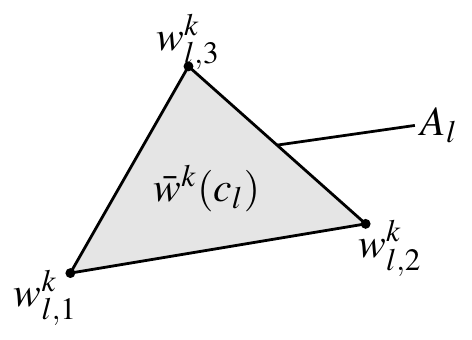}
	\caption{Placement of vertical velocities $w_{l,i}^k$ at the element $e_l$ in horizontal layer $k$. The interpolated velocity $w^k(c_l)$ is located at the cell center $c_l$ of the respective element. Surface triangle face areas $A_l$ are constant through all levels.}
\end{figure}
The sections $S_i$, compare Fig.\ref{FIG:AMOCBIN}, contain all cell center velocities $w(c_l) \in S_i$ with a latitudinal coordinate $c_l^\phi \in [\phi_{i-1},\phi_i]$.
On each latitudinal section the vertical transports of the respective bin are added along the longitude:
\begin{equation}
	\Delta \Psi_{ki} \; = \; \sum_{w(c_l) \in S_i} w_{k,l} A_c \; ,
\end{equation}
where $k$ denotes the horizontal level and $i$ the corresponds to the section $S_i$. The final step denotes a sum over the bins $S_i$ up to the desired longitude:
\begin{equation}
	\Psi_{ki} \; = \; \sum_{j=1}^{i} \Delta \Psi_{kj}
\end{equation}

However, for the two estimates on the discrete level exist differences up to several $Sv$, as emphasized in \cite{FESOM2amoc}. Furthermore, determining the transport by the above algorithm adds interpolation errors to the outcome.  The overall question how to solve for transports on unstructured meshes still poses a challenge \cite{FESOM2amoc}. Despite this level of uncertainty towards the results these approaches allow for an interpretation of the AMOC in simulations.

\subsection{Mesh Evaluation}

The results in this section, in the same way as the numerical experiments later on, are not measured by their ability to represent the physics of the oceans with sufficient accuracy. In the context of Parareal, we found it more important to investigate differences in the solutions than obtaining an acceptable representation of the oceans. Especially for the latter, we knew that the PI mesh is used exclusively for testing purposes and that refining the resolution from the outset would not be sufficient. For the application of Parareal we are predominantly interested if we can converge from the coarse initial solution (PI mesh) to the fine reference solution (FPI mesh). The differences in these results are given in this section. 

\subsubsection{Temperature}

\begin{figure}[H]
	\centering
	\includegraphics{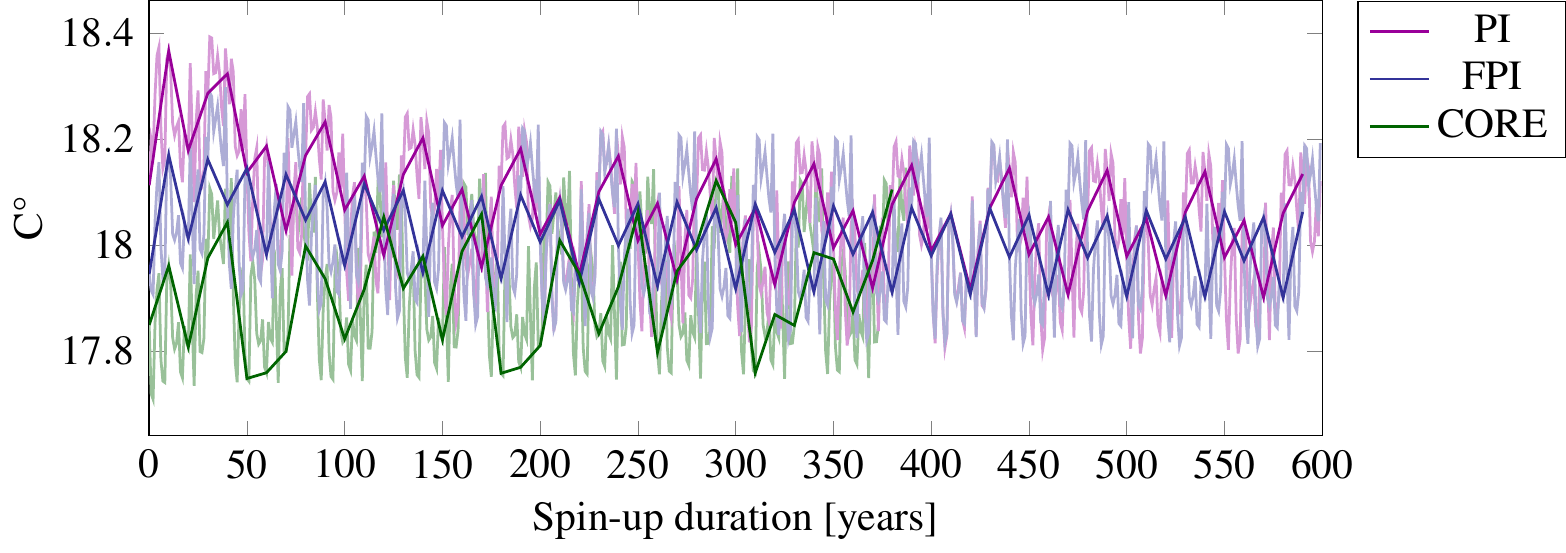}
	\caption{Annual global mean temperature at ocean surface. Solid lines represent every 10th value for the sake of visibility. Transparent lines show the complete temporal evolution of the temperature.}
	\label{FIG:MeanTemp0}
\end{figure}

In a first step the mean temperatures at different ocean depths will be investigated. We chose the annual mean temperature at the ocean's surface, in 500m depth and an average of the layer from 0 to 500m. The CORE mesh was used to put the results of the PI and FPI runs in context. Since there are no analytical solutions available for simulations on both test meshes we decided to present them alongside the FESOM2 standard case. \\
In Fig.\ref{FIG:MeanTemp0} the sea surface temperature is given. At the beginning of the simulation the PI mesh provides the highest temperature, followed by the FPI mesh. With increasing simulation time the profiles seem to settle around 18 degrees and suggest an converged state for the sea surface temperature.  

\begin{figure}[H]
	\centering
	\includegraphics{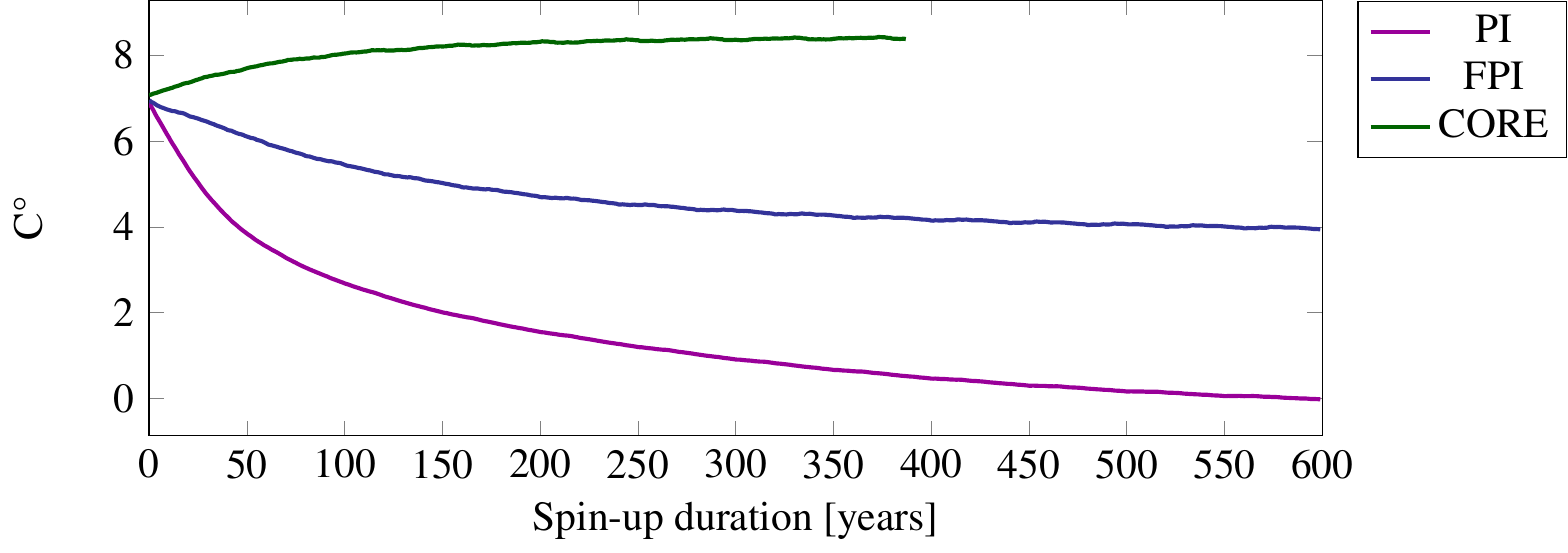}
	\caption{Annual mean temperature at 500m depth.}
	\label{FIG:MeanTemp500}
\end{figure}

\begin{figure}[H]
	\centering
	\includegraphics{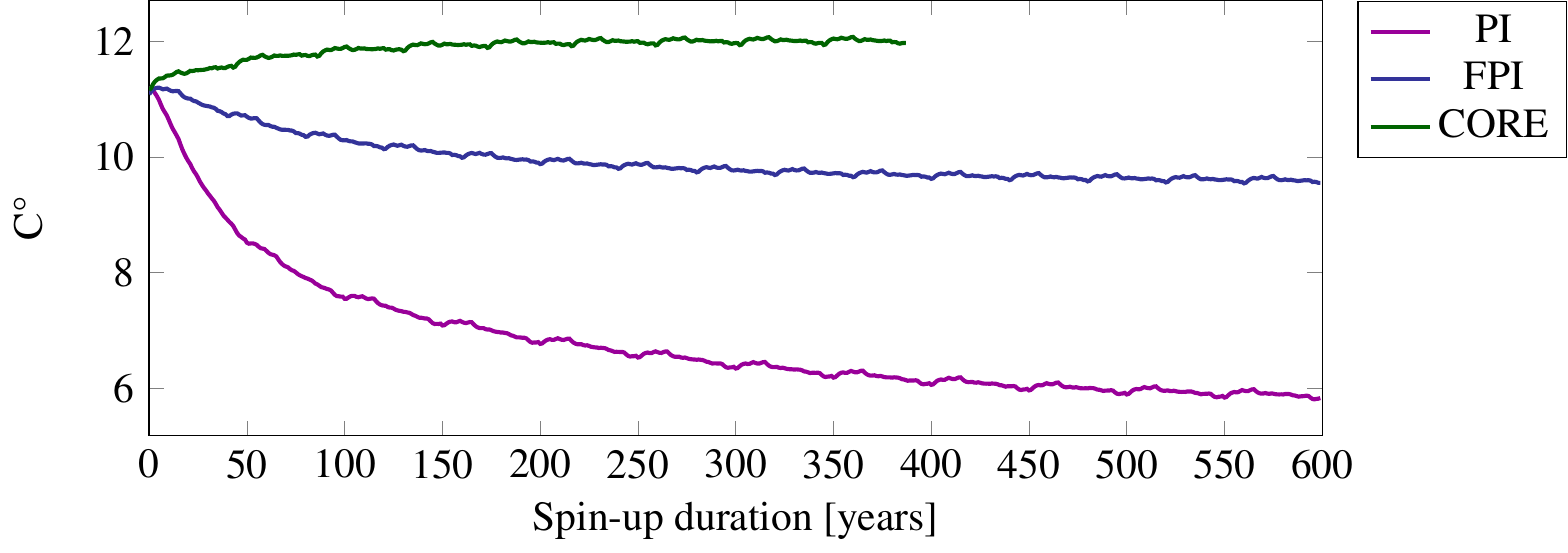}
	\caption{Global mean temperature integrated over 0-500m depth.}
	\label{FIG:MeanTemp0500}
\end{figure}

With the results given in Figs.\ref{FIG:MeanTemp500} and \ref{FIG:MeanTemp0500} we conclude that there is a minimum spatial resolution mandatory in order to obtain sensible results. In Fig.\ref{FIG:MeanTemp500} we find the PI mesh temperature at 500m depth to decrease to zero degrees Celsius after 600 years. For FPI mesh the temperature decline stops at 4 degrees Celsius. The Core mesh shows convergence to ~8°C after 400 years with a slight increase towards the initial condition.

For the averaged temperature over 0 to 500m depth in Fig.\ref{FIG:MeanTemp0500} the same tendency is observed as in Fig.\ref{FIG:MeanTemp500}. Although, the temperature profiles are shifted by +4°C the decline over 600 years is similar to the profiles at 500m depth. The CORE mesh again converges towards 400 years simulation time. The observed minimal oscillations are contributions from the variations at sea surface, where the interface to the forcing data set allows for a reasonable physical representation of the oceans. With decreasing spatial resolution and ocean depth we found the temperature to be further declining. 

\subsubsection{AMOC}

\begin{figure}[H]
	\centering
	\includegraphics{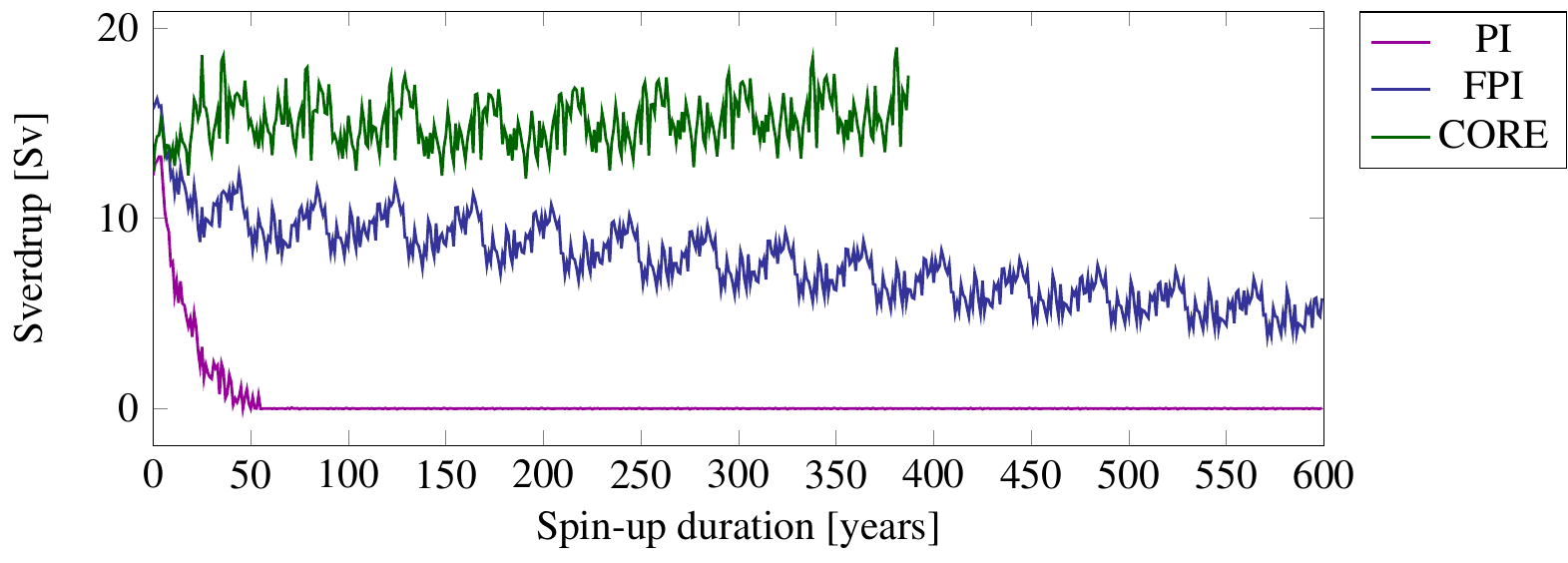}
	\caption{Comparison of the AMOC at 26.5N latitude. The PI mesh AMOC entirely vanishes within in the first CORE2 forcing cycle. The FPI AMOC shows a steady decline over time and reaches no converged state.}
	\label{FIG:AMOC}
\end{figure}

In Fig.\ref{FIG:AMOC} the results of the AMOC evaluation at 26.5N is given. The PI mesh appears to fail in conserving the momentum and therefore, we find the the AMOC collapsing to a standstill of zero Sverdrup. With the increased resolution of the FPI mesh we still observe a steady decline over time, but not as immediate as in the PI mesh simulation. These results suggest strong dependency on the meshes' resolutions for FESOM2 to conserve momentum and approximate the complex ocean dynamics adequately. 
With the exception of the PI grid, the variability of the CORE2 forcing cycles over 58 years during the spin-up is observable. The repetitive pattern is more striking for the FPI mesh results as for the CORE mesh reference. In the later CORE2 cycles the pattern in FPI results becomes less distinctive. The PI mesh profile of the AMOC lacks repetitive variability as the ocean dynamics immediately break down.
\begin{figure}[H]
	\centering
	\begin{subfigure}[t]{.49\textwidth}
			\includegraphics[width=\textwidth]{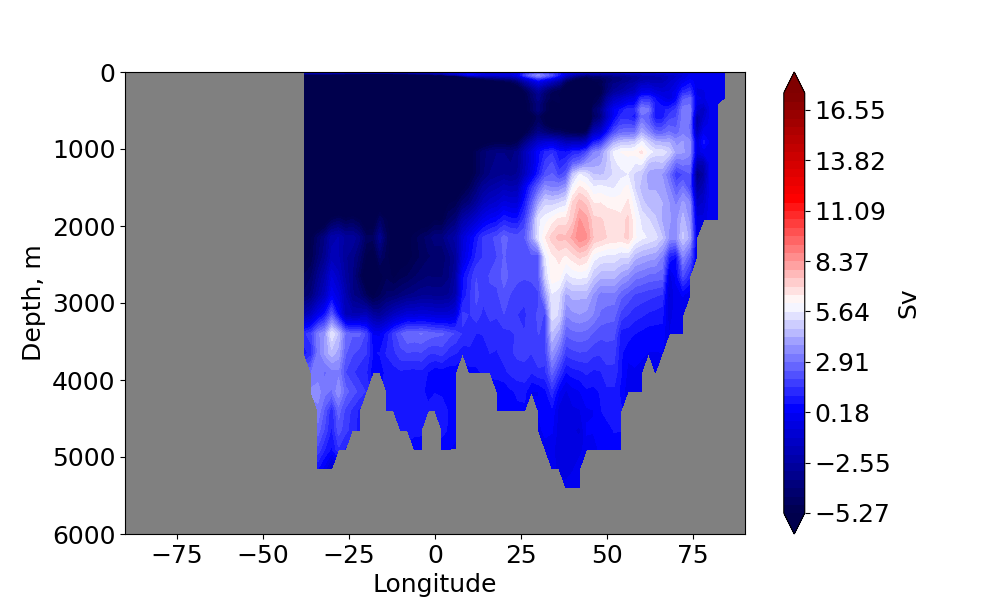}
			\caption{PI AMOC in 1948}
		\end{subfigure}
	\hfill
	\begin{subfigure}[t]{.49\textwidth}
			\includegraphics[width=\textwidth]{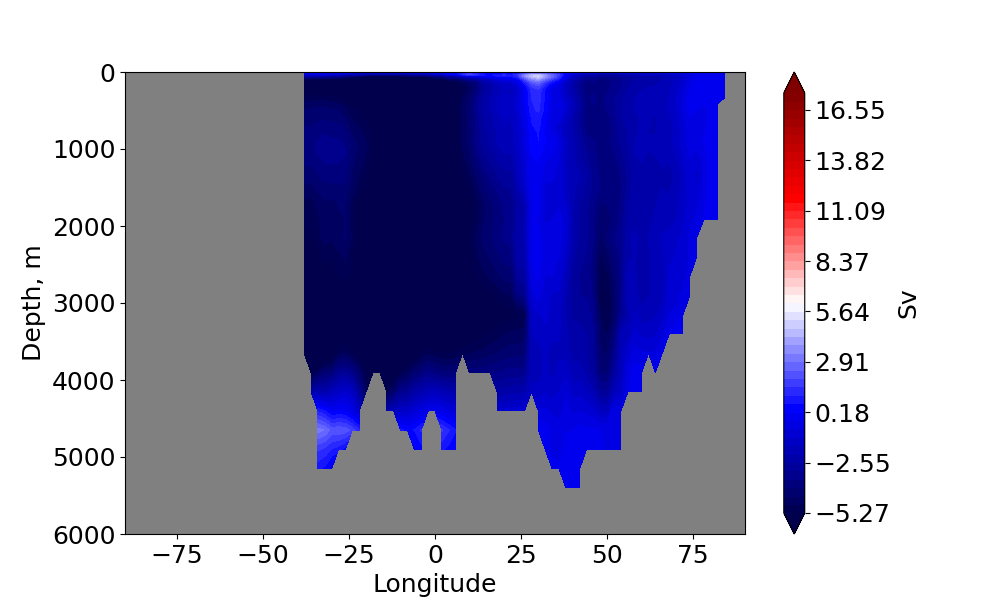}
			\caption{PI AMOC after spinup}
		\end{subfigure}
		\begin{subfigure}[t]{.49\textwidth}
			\includegraphics[width=\textwidth]{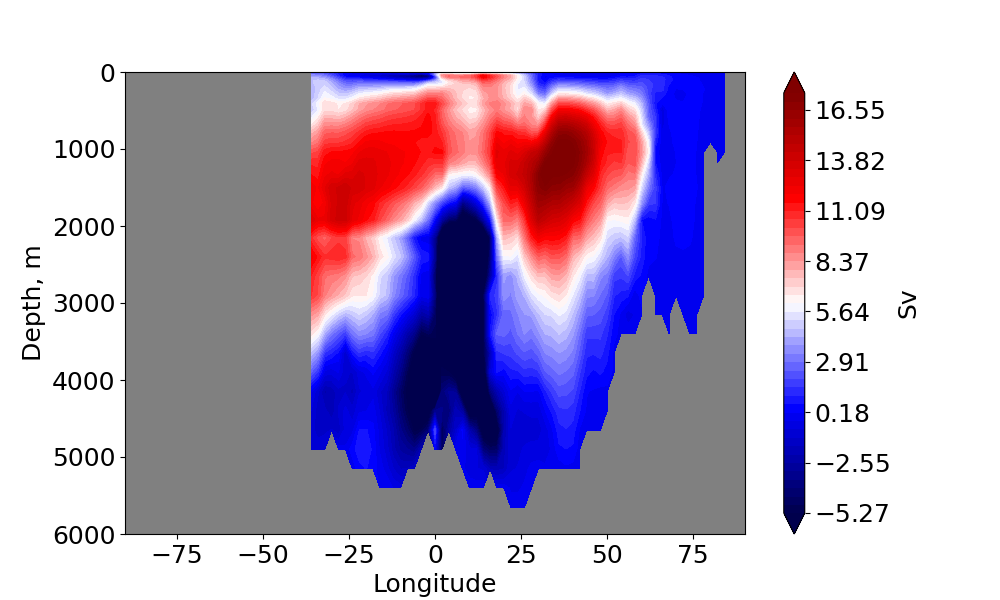}
			\caption{FPI AMOC in 1948}
		\end{subfigure}
	\hfill
	\begin{subfigure}[t]{.49\textwidth}
			\includegraphics[width=\textwidth]{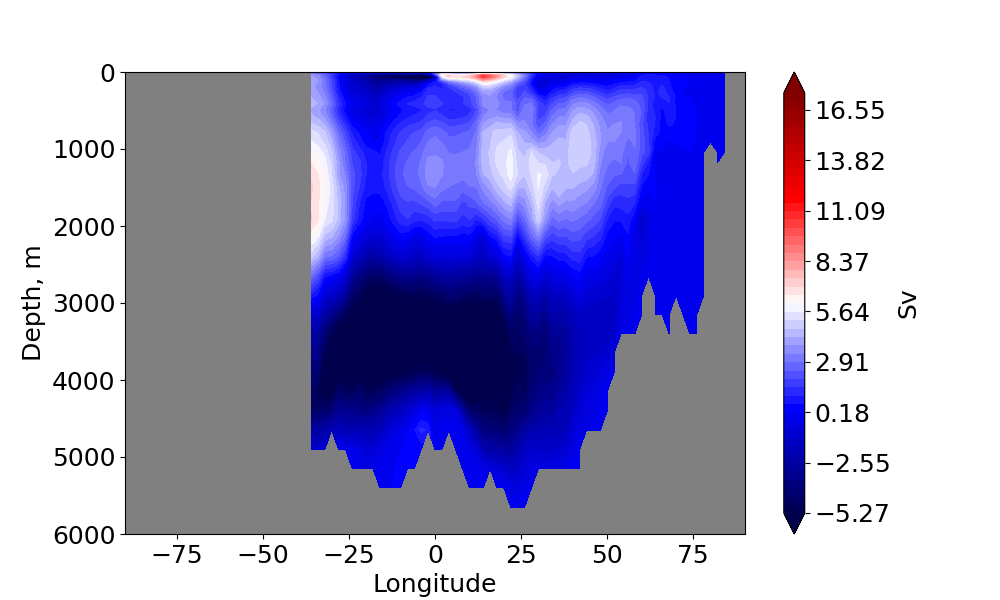}
			\caption{FPI AMOC after spinup}
		\end{subfigure}
	\begin{subfigure}[t]{.49\textwidth}
			\includegraphics[width=\textwidth]{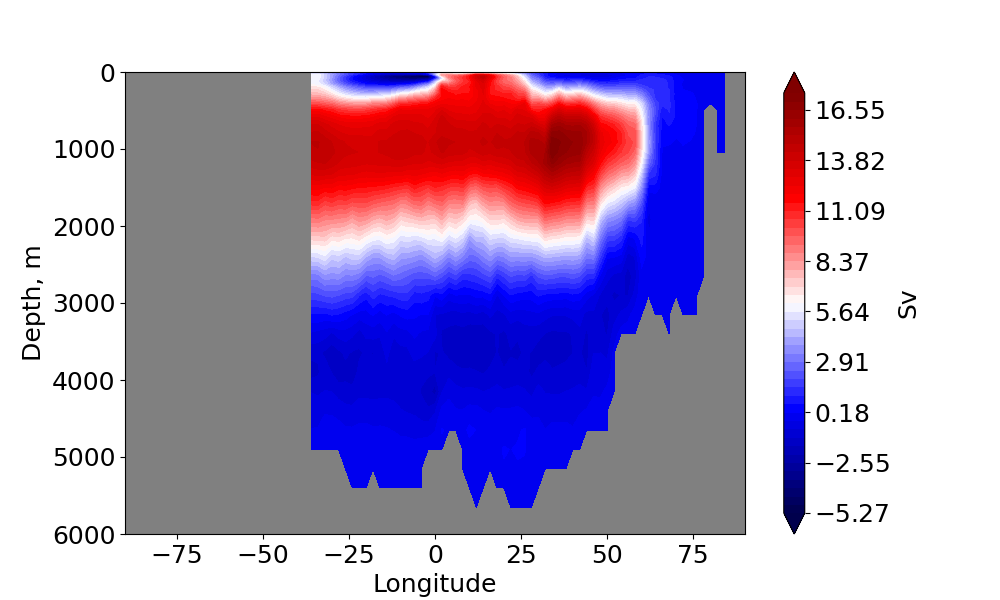}
			\caption{AMOC of the CORE run.}
		\end{subfigure}
	\caption{2D visualization of the AMOC between 40S to 80N. First row shows the AMOC from the PI mesh at the beginning (a) and end (b) of the spin-up. Second row represents results from the FPI setting and in (e) the CORE mesh is given as a reference.}
	\label{FIG:AMOC2D}
\end{figure}

In Fig.\ref{FIG:AMOC2D} a 2D-representation of the AMOC is given along 40S to 80N through the Atlantic Ocean. The AMOC was evaluated after the first year of the spin-up, where the results are close to the initial condition imposed by the climatology data set, and after 600 years of simulation time. In agreement with the results from the AMOC evaluation at 26.5N, we observe a breakdown of the ocean circulation for the PI mesh configuration.

In Fig.\ref{FIG:AMOC2D}(a) the loss in the circulation magnitude is already observed after one year. In Fig.\ref{FIG:AMOC2D}(b) the AMOC is close to zero Sverdrup throughout the evaluated part of the Atlantic Ocean. The PI mesh fails to provide any of the structure seen in the CORE mesh reference in Fig.\ref{FIG:AMOC2D}(e). For the FPI mesh the decline in the circulations magnitude is given in Figs.\ref{FIG:AMOC2D}(c) and \ref{FIG:AMOC2D}(d). Although the structure clearly is less developed the FPI configuration allows for reduced dynamics. \\

With the observed break-down (PI) or reduction (FPI) of the currents in the Atlantic Ocean, the results contribute to the understanding of the mean temperature profiles at surface level and 500m depth. Without transport in the layers below the sea surface temperature for the PI and FPI configuration is slightly higher than for the CORE mesh. Correspondingly, an cool down of the layers below due to the lack of currents is observed.

\subsubsection{Horizontal Mean Velocity}

\begin{figure}[H]
	\centering
	\begin{subfigure}[t]{.49\textwidth}
			\includegraphics[width=\textwidth]{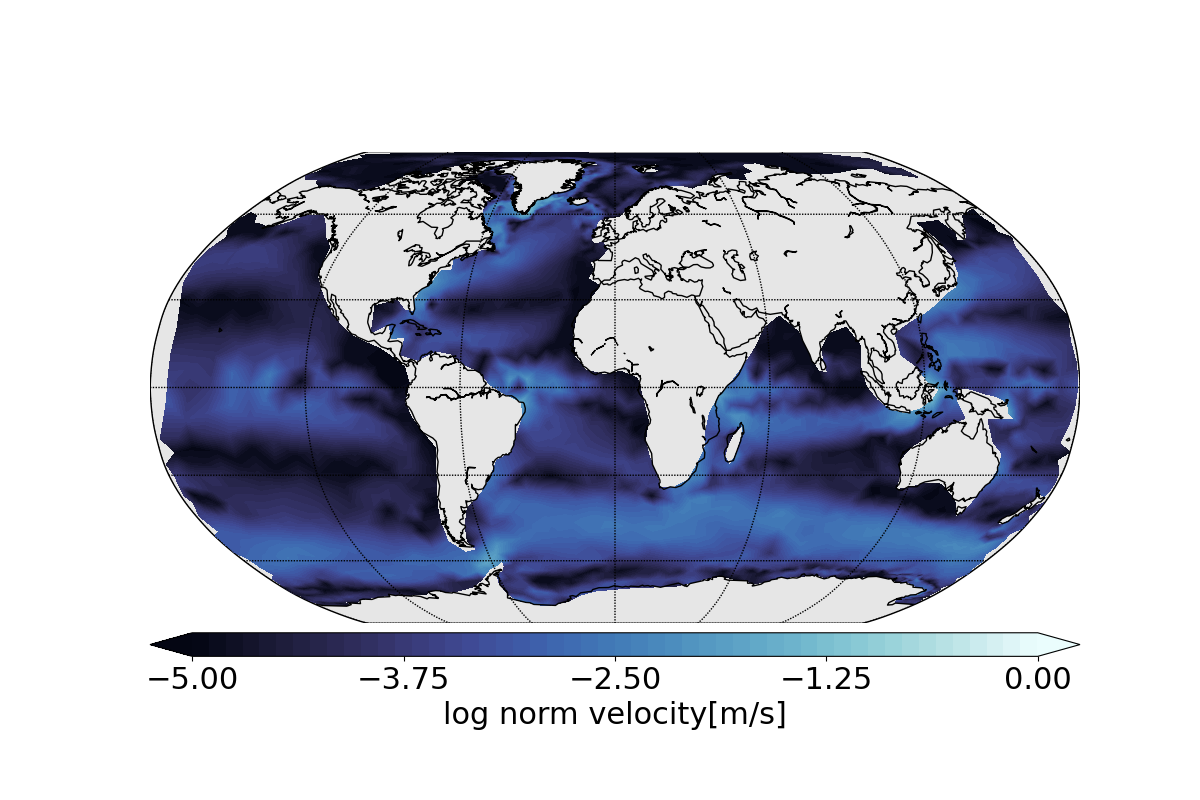}
			\caption{PI initial state at 1948}
		\end{subfigure}
	\hfill
	\begin{subfigure}[t]{.49\textwidth}
			\includegraphics[width=\textwidth]{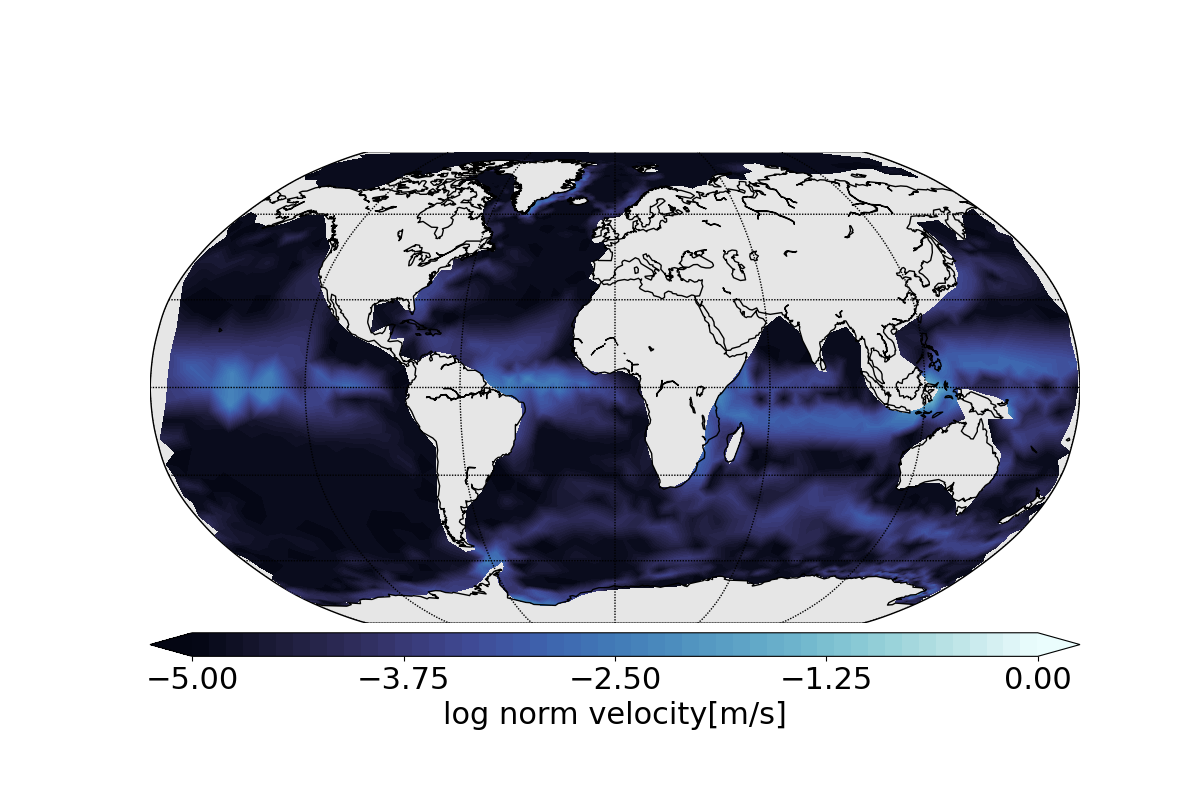}
			\caption{PI after spinup}
		\end{subfigure}
	\begin{subfigure}[t]{.49\textwidth}
			\includegraphics[width=\textwidth]{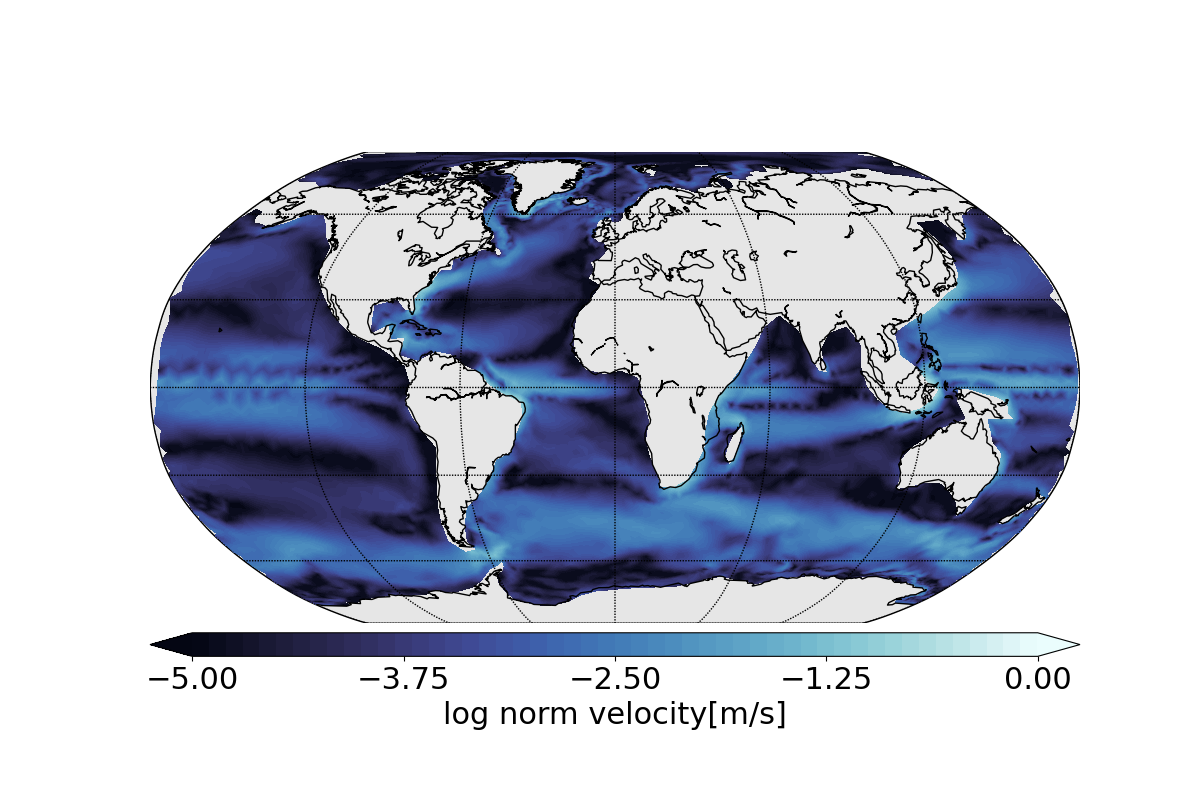}
			\caption{FPI initial state 1948}
		\end{subfigure}
	\hfill
	\begin{subfigure}[t]{.49\textwidth}
			\includegraphics[width=\textwidth]{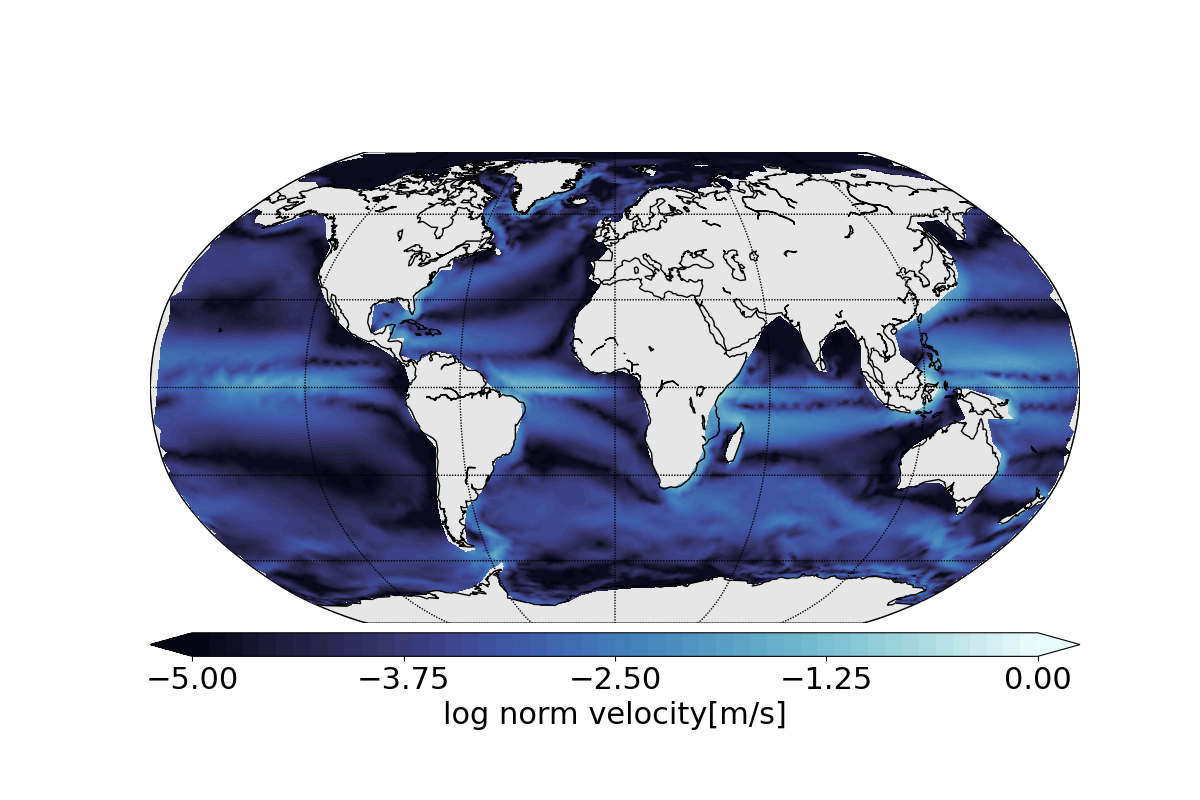}
			\caption{FPI after spinup}
		\end{subfigure}
	\begin{subfigure}[t]{.49\textwidth}
			\includegraphics[width=\textwidth]{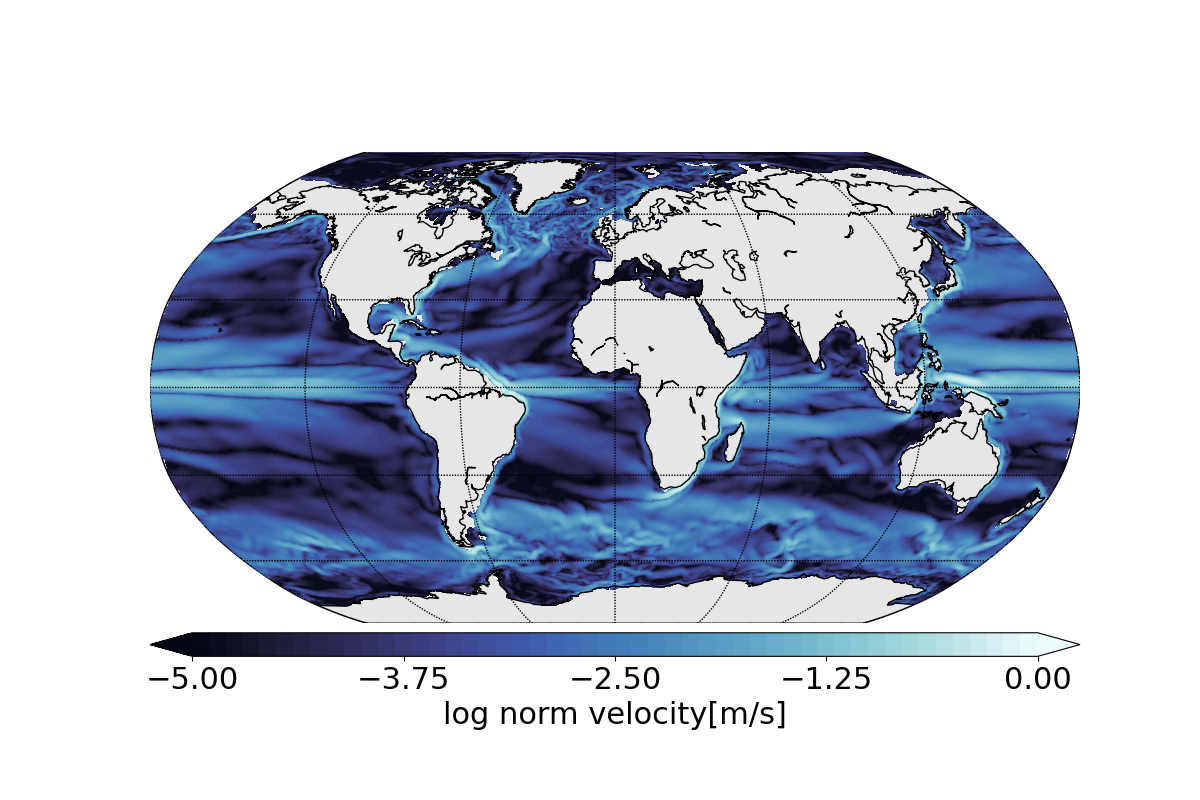}
			\caption{CORE reference run.}
		\end{subfigure}
	\caption{2D visualization of the mean horizontal velocities at the surface.}
	\label{FIG:MEANVEL}
\end{figure}

In Fig.\ref{FIG:MEANVEL} the horizontal mean velocities on the surface are compared. In the left row the velocity is given after the first year in 1948 and in the right at the end of the spin-up run. In agreement with the evaluation of the AMOC the PI mesh fails to preserve the ocean dynamics. After one year of simulation already, see Fig.\ref{FIG:MEANVEL}a), a decrease in magnitude can be observed. The stillstand mentioned in the AMOC evaluation is clearly visible in Fig.\ref{FIG:MEANVEL}b). The FPI mesh proofs better ability to preserve the ocean currents, but also suffers from losses in magnitude of the velocities after the spin-up. For comparison we attached the results of the CORE mesh, which were evaluated after 300 years of simulation. Due to its ability to preserve momentum the currents after one year of simulation are virtually indistinguishable by these visualizations.

\subsubsection{Salinity}

\begin{figure}[H]
	\centering
	\begin{subfigure}[t]{.49\textwidth}
			\includegraphics[width=\textwidth]{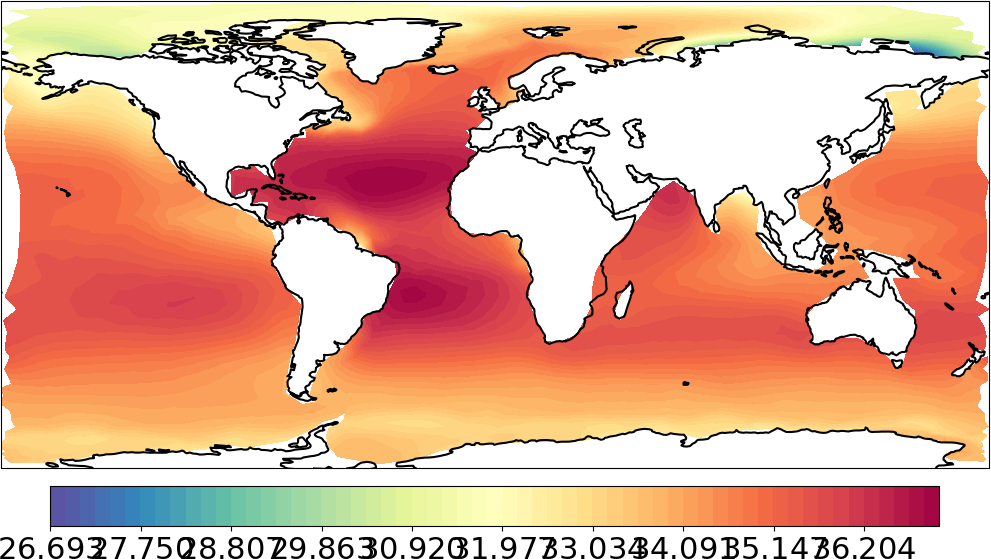}
			\caption{PI: Initial salinity in 1948}
		\end{subfigure}
	\hfill
	\begin{subfigure}[t]{.49\textwidth}
			\includegraphics[width=\textwidth]{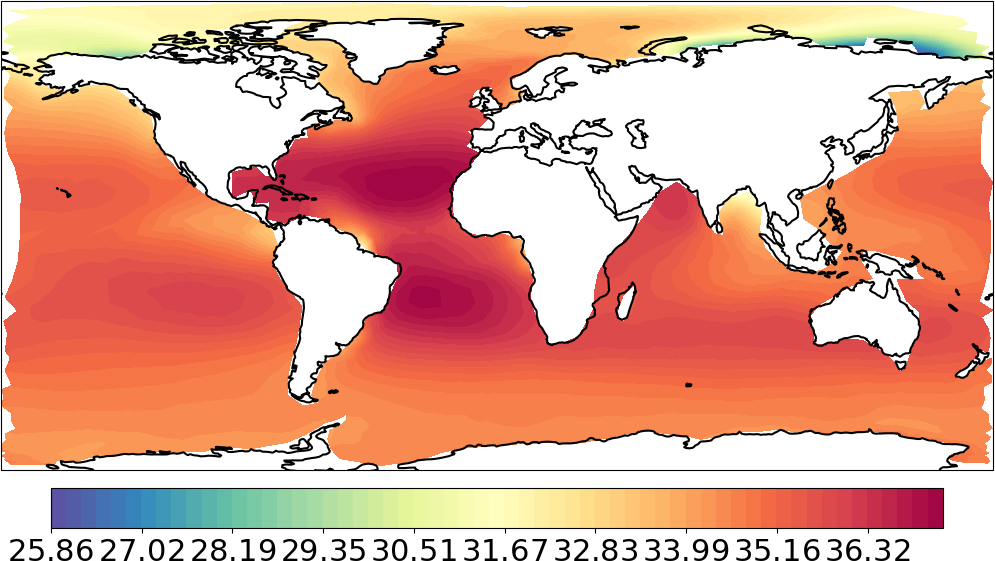}
			\caption{PI: Initial salinity after spinup}
		\end{subfigure}
		\begin{subfigure}[t]{.49\textwidth}
			\includegraphics[width=\textwidth]{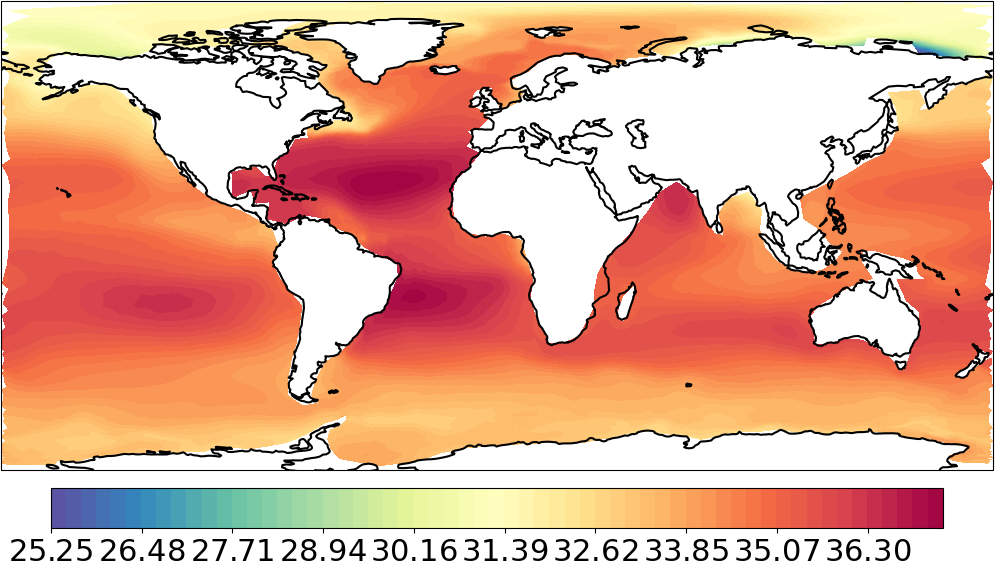}
			\caption{FPI: Initial salinity in 1948}
		\end{subfigure}
	\hfill
	\begin{subfigure}[t]{.49\textwidth}
			\includegraphics[width=\textwidth]{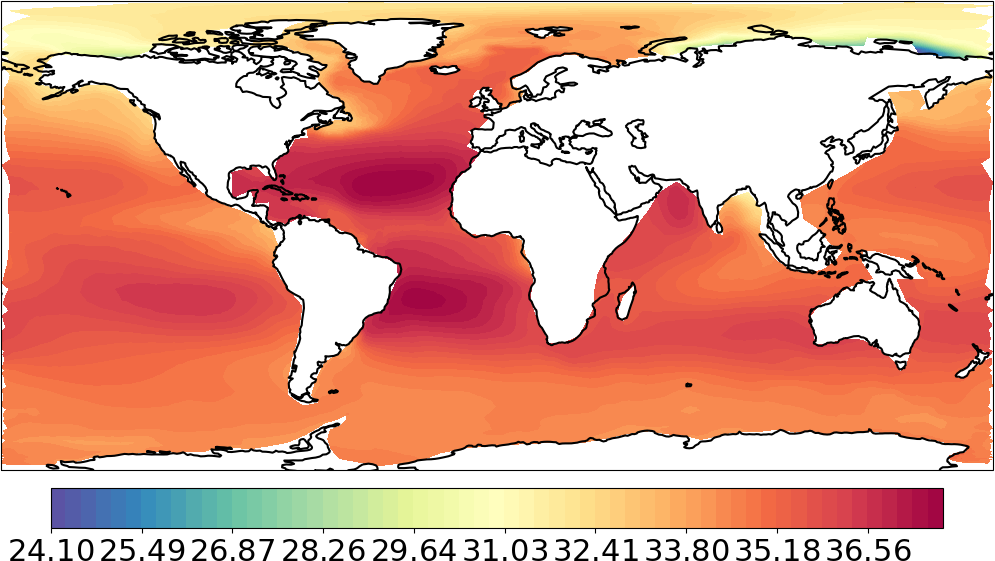}
			\caption{FPI: Initial salinity after spinup}
		\end{subfigure}
	\begin{subfigure}[t]{.49\textwidth}
			\includegraphics[width=\textwidth]{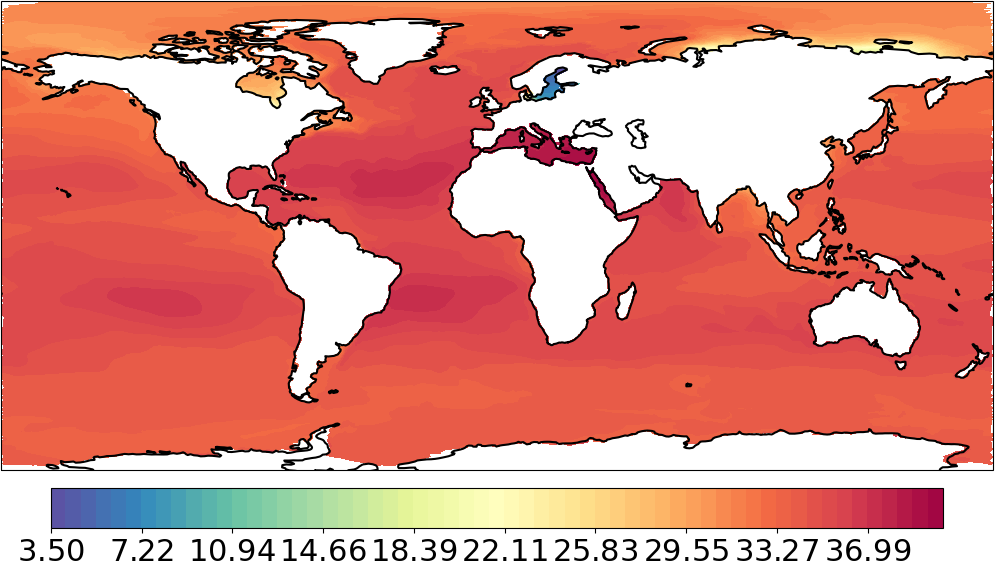}
			\caption{CORE: Initial salinity in 1948}
		\end{subfigure}
	\hfill
	\begin{subfigure}[t]{.49\textwidth}
			\includegraphics[width=\textwidth]{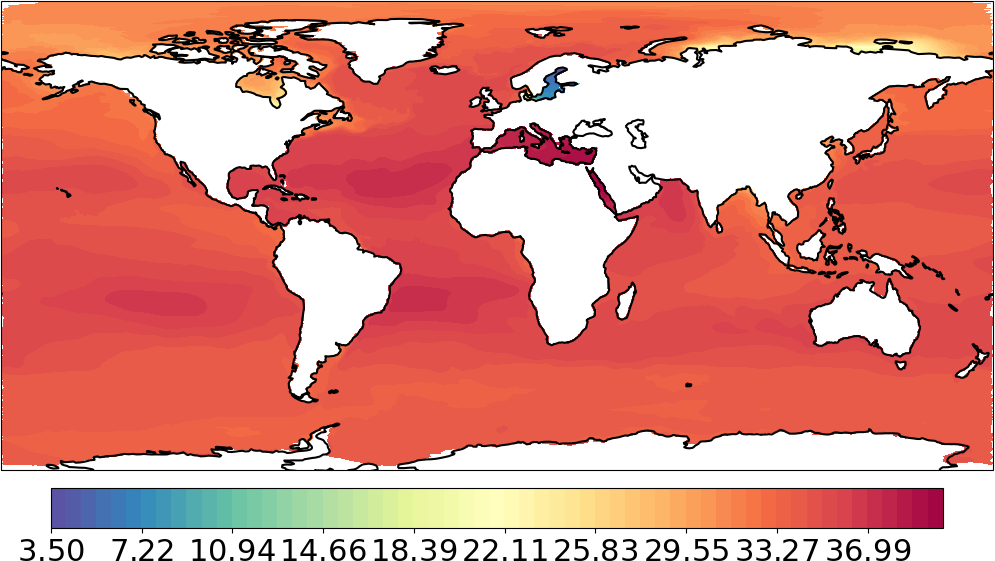}
			\caption{CORE: Initial salinity after spinup}
		\end{subfigure}
	\caption{Salinity comparison at surface level at the beginning of and after the spin-up.}
	\label{FIG:SALT}
\end{figure}

In Fig.\ref{FIG:SALT} we compare the salinity distribution at the ocean surface. The differences in the distributions at the beginning and end of the spin-up show minimal differences. The major difference is found on the CORE mesh for Baltic Sea, which is low in Salinity due to large contribution of fresh water fluxes. Apart from the outlier in the Baltic Sea, all meshes show comparable data ranges for salinity. We found this statement to hold true also for the lower levels of the ocean and for the sake of brevity decided against a visualization at this point.

\subsubsection{Impact of Forcing}

\begin{figure}[H]
	\centering
	\includegraphics{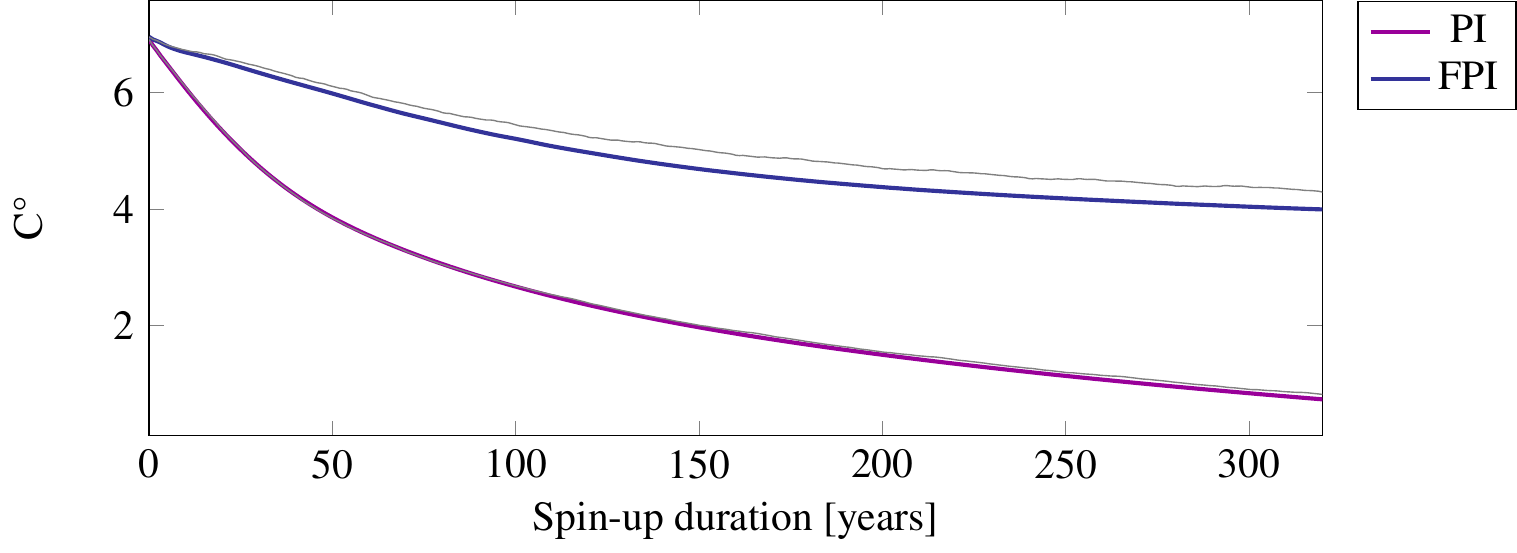}
	\caption{Global mean temperature at 500m depth with one year forcing. Gray lines indicate}
	\label{FIG:MeanTempNOF}
\end{figure}

\begin{figure}[H]
	\centering
	\includegraphics{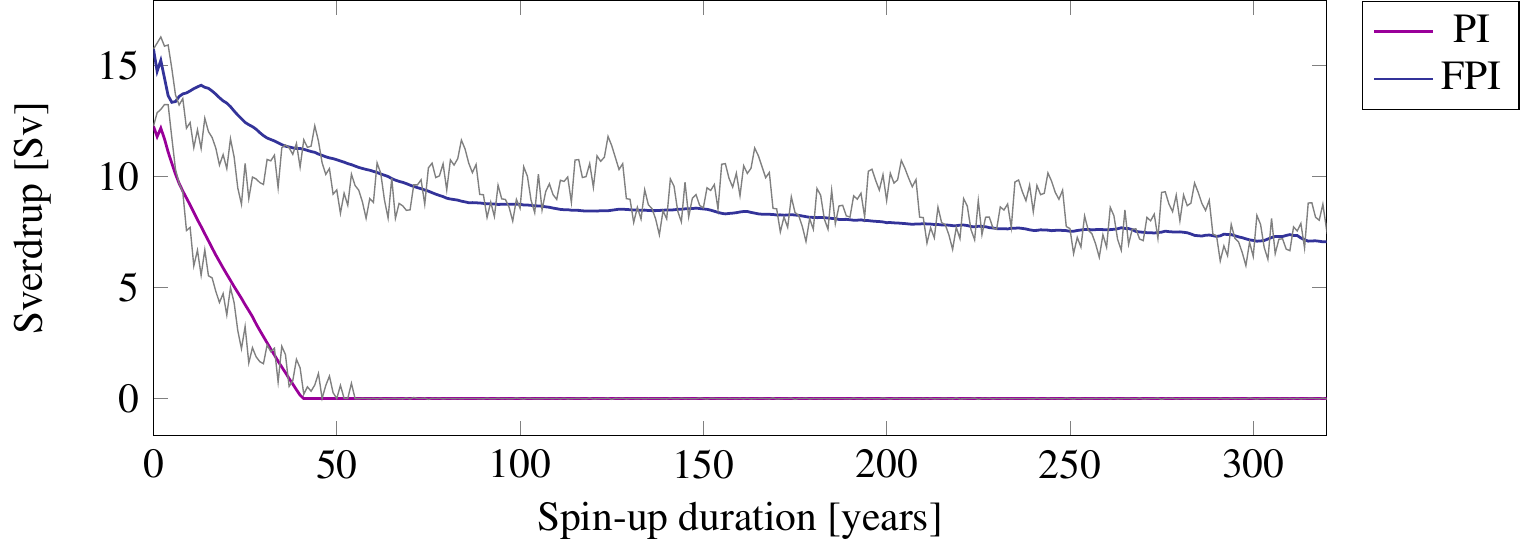}
	\caption{Comparison AMOC.}
	\label{FIG:AMOCNOF}
\end{figure}

With respect to achieving a converged state for the spin-up simulation we investigated the impact of the forcing data set on the rate of convergence. Since the ocean circulation and therefore FESOM2 relies on the input, we decided to run a spin-up over 300 years with the forcing data solely from 1948. By removing the variability in the results due to CORE2 forcing over 58 years, we intended to exclude possibly unwanted contributions to the converged state. The results of those runs are given for the mean temperature at 500m depth in Fig.\ref{FIG:MeanTempNOF} and for the AMOC at 26.5N in Fig.\ref{FIG:AMOCNOF}. The gray lines in both figures represent the respective run with full forcing cycles applied. \\
In Fig.\ref{FIG:MeanTempNOF} the mean temperature shows the same decline over time while converging slightly faster. For the FPI run the temperature difference is 0.3°C and 0.09°C for the PI setting. The AMOC evaluation at 26.5N in Fig.\ref{FIG:AMOCNOF} depicts a smoother profile for both runs, since the variability imposed by the entire forcing data set is not permitted. The general evolution of the diagnostic variables in time is not affected by different integration of the CORE2 forcing.

We conclude, that the impact of the CORE2 forcing cycle is of minor importance within the context of converging to a steady-state in spin-up simulations. Although the cyclic behavior in the variability due to repeated CORE2 spin-up simulations is self-evidently missing, there is no contribution to a mesh's ability to preserve momentum or energy. Hence, an appropriate representation of the oceans is determined primarily by the spatial resolution.

\subsubsection{Impact of Temporal Resolution on the Annual Averaged Temperature}

\begin{figure}[H]
	\centering
	\includegraphics{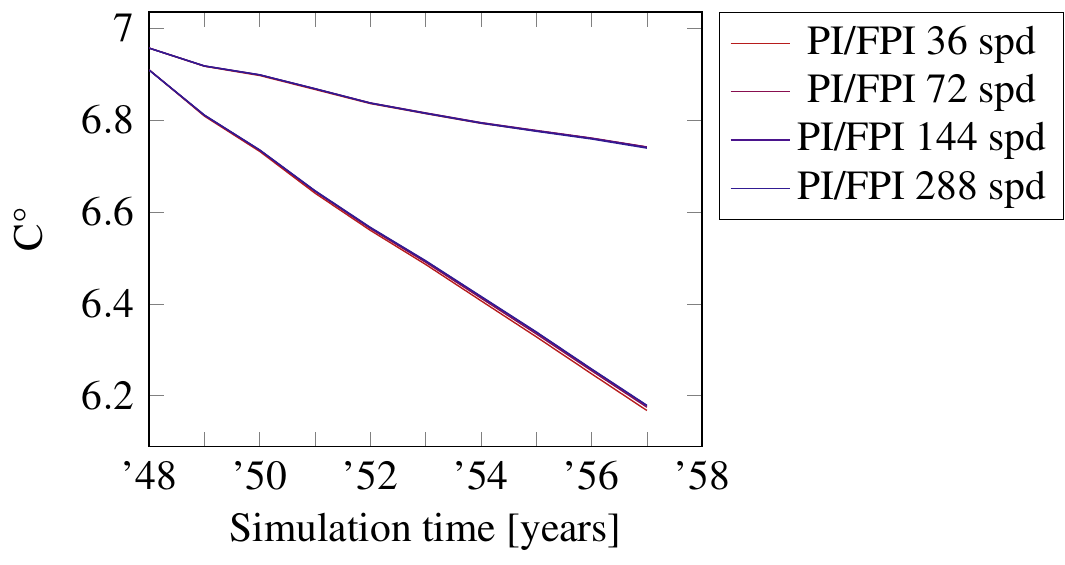}
	\caption{Comparison of the impact of the temporal resolution with regard to different resolutions of the PI and FPI mesh. With higher resolution the annual averaged temperature at 500m depth appears to be less affected.}
	\label{FIG:MeshTime}
\end{figure}

To conclude the mesh evaluation we investigated the impact of the temporal resolution on the mean temperature. This investigation is kept to a minimum, since it has been addressed thoroughly in \cite{philippi2022parareal}, where a Parareal implementation for the PI mesh with different time step sizes is given. In Fig.\ref{FIG:MeshTime} the mean temperature results at 500m depth of simulations over a period of 10 years with different time step sizes are shown. The largest time steps of 2400 seconds corresponds to 36 spd and the smallest of 300 seconds to 288 spd, respectively.  

Virtually there is no difference in the temperature profiles to be found. Hence, the variation of time step sizes fails to improve simulation results and therefore, the largest possible time step is the configuration of choice for serial computations. Though, in the documentation \url{https://fesom2.readthedocs.io/en/latest/index.html} it is recommended for higher resolution meshes to refine the time step size when encountering stability problems. In the context of FESOM2 the time step size can be considered as a possibility to maintain stability during simulations, which is exploited in one of the Parareal experiments later on. \\

With the closure of this section we will address the time step size of 36spd being applicable for PI, FPI and CORE mesh. Although all meshes investigated in the mesh evaluation section have different horizontal resolution, they share the same amount of vertical layers. The dimensions of the horizontal planes are significantly larger than of the ocean's depth. In all meshes used in this study the lowest depth is 6.25km. On the contrary, the Pacific ocean, as the largest one, is spanning a width up to 20.000km. Hence, the vertical resolution presents the restricting factor when it comes to the largest time step size, which is identical for the three meshes.

\section{Parareal on the DKRZ Cluster}

In this section the realization of the Micro-Macro-Parareal algorithm on the DKRZ cluster Levante is given. The management and job scheduling system provided on Levante is SLURM. All algorithm steps carried out are provided in a sbatch script file. The file contains the execution of FESOM2 as fine and coarse propagator, as well as the instructions for file manipulation and interpolation tasks. Will all instructions gathered within one script file only one job has to be submitted, where the total amount of nodes and CPUs required for the Parareal experiment have to be defined beforehand. 

\subsection{SLURM}

\begin{figure}[H]
	\centering
	\tikzsetnextfilename{FESOM2folder}
\begin{forest}
	for tree={
		font=\ttfamily,
		grow'=0,
		child anchor=west,
		parent anchor=south,
		anchor=west,
		calign=first,
		inner xsep=7pt,
		edge path={
			\noexpand\path [draw, \forestoption{edge}]
			(!u.south west) +(7.5pt,0) |- (.child anchor) pic {folder} \forestoption{edge label};
		},
		before typesetting nodes={
			if n=1
			{insert before={[,phantom]}}
			{}
		},
		fit=band,
		before computing xy={l=15pt},
		file/.style={edge path={\noexpand\path [draw, \forestoption{edge}]
				(!u.south west) +(7.5pt,0) |- (.child anchor) \forestoption{edge label};},
			inner xsep=2pt,font=\small\ttfamily
		}
	}  
	[\$PATH
	[G
	[slice x
	[namelist.config, file]
	[namelist.forcing, file]
	[namelist.oce, file]
	[namelist.ice, file]
	[fesom.x, file]
	[output
	[fesom.1948.restart.oce.nc, file]
	[fesom.1948.restart.ice.nc, file]
	[fesom.clock, file]
	]
	]
	]
	]
\end{forest}
	\caption{Folder structure for FESOM2 computations on Levante is given. The example shows the folder before the coarse propagator execution in an arbitrary slice x with initial values from the year 1948.}
	\label{FIG:FESOM2folder}
\end{figure}
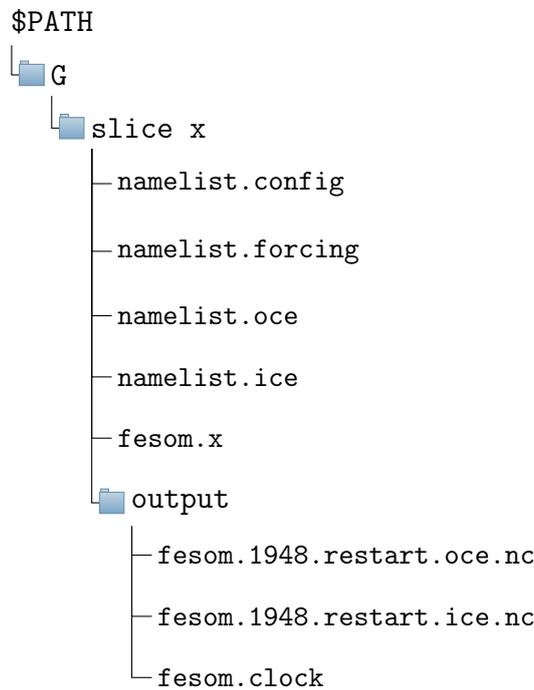

FESOM2 requires a pre-defined folder structure where the code can be executed on cluster. Hence, for every call of FESOM2 during the algorithm's process an own folder has to be provided. The general folder structure for each call of FESOM2 is given in Fig.\ref{FIG:FESOM2folder}.

In \verb|namelist.config| contains the general configuration settings, where the location of the mesh, amount of CPUs for spatial domain decomposition, time stepping and simulation interval are given. In the remaining \verb|namelist| files parameters for the physical models of ocean, ice and forcing are defined. We kept all parameters in these files as they were provided with downloading the standard basic configuration from the developers github. In \verb|namelist.oce| are around 50 parameters to be found, which allow for modification and adjustment of the ocean model. We considered the investigation of the impact of parameter settings on simulation results beyond the scope of this study. Even more so, after the evaluation of our test meshes PI and FPI has shown the limits of low resolution meshes. \\

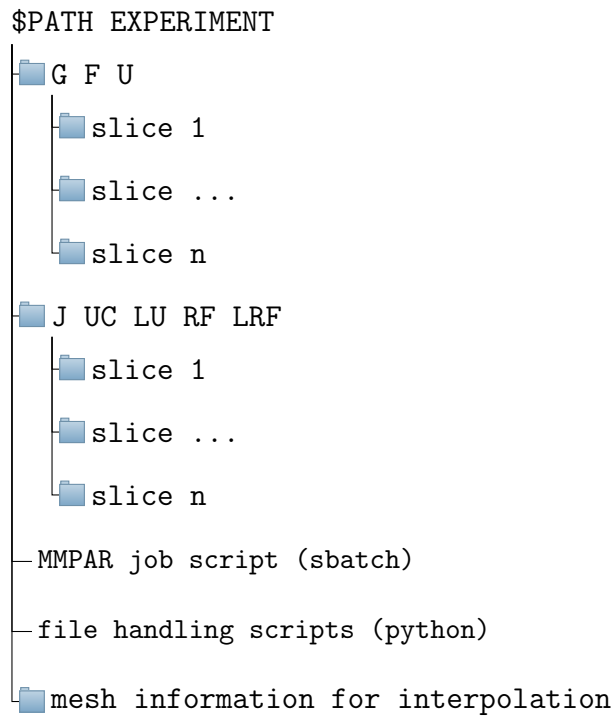
\begin{figure}[H]
	\centering
	\tikzsetnextfilename{MMPARfolder}
\begin{forest}
	for tree={
		font=\ttfamily,
		grow'=0,
		child anchor=west,
		parent anchor=south,
		anchor=west,
		calign=first,
		inner xsep=7pt,
		edge path={
			\noexpand\path [draw, \forestoption{edge}]
			(!u.south west) +(7.5pt,0) |- (.child anchor) pic {folder} \forestoption{edge label};
		},
		before typesetting nodes={
			if n=1
			{insert before={[,phantom]}}
			{}
		},
		fit=band,
		before computing xy={l=15pt},
		file/.style={edge path={\noexpand\path [draw, \forestoption{edge}]
				(!u.south west) +(7.5pt,0) |- (.child anchor) \forestoption{edge label};},
			inner xsep=2pt,font=\small\ttfamily
		}
	}  
	[\$PATH EXPERIMENT
	[G F U
	[slice 1]
	[slice \dots]
	[slice n]
	]
	[J  UC LU RF LRF
	[slice 1]
	[slice \dots]
	[slice n]
	]
	[MMPAR job script (sbatch), file]
	[file handling scripts (python), file]
	[mesh information for interpolation]
	]	
\end{forest}
	\caption{Folder structure for a Parareal run on SLURM. The G, F and U folders contain the results from FESOM2 simulations on a pre-defined amount of slides. J, UC, LU, RF and LRF are auxiliary folders for the temporary storage of the Parareal jumps and interpolations. }
	\label{FIG:MMPARfolder}
\end{figure}

For each fine and coarse call of FESOM2 a dedicated folder must be provided. For each Parareal test case a parent folder is created, containing the run folders for $F$ and $G$ as well as a storage folder for the iteratives $U$. For the sake of traceability and possible error analysis folders for the computation of the Parareal jumps and interpolated files were created. While these auxiliary folders help to keep the process of file handling manageable, there are obviously possibilities for optimization left. Additionally, the Parareal sbatch file and the python scripts for interpolation tasks are given. The parent folder structure is given in Fig.\ref{FIG:MMPARfolder}.

With the underlying folders structure being introduced excerpts from the sbatch job script are given to explain how serial and parallel executions of FESOM2 are carried out. In this section we will restrict the introduction of the sbatch script to the key instructions, since the entire script contains numerous tasks, e.g. the declaration of run-time variables. In Lst.\ref{lst:serialG} the serial progression of the computation of the initial coarse solution $U^0$ is given. At the end of each run within its respective folder the results are copied to the consecutive folder. 

\begin{lstlisting}[language=bash, caption=Extract from the Parareal SLURM Script: Serial execution of the coarse solver, label=lst:serialG]
for (( i=$TA ; i<=$TE ; i++ )); do
 
 ...
 
 # $TA and $TE define the first and last slice in the experiment
 # $pathg defines the absolute path to the coarse folder
 # $flgs_c defines the amount of processes for spatial parallelization
 
 srun $flgs_c --chdir=$pathg/slice$i ./fesom.x > \ 
                     "$pathg/slice$i/fesom2.0.out"
 ...
 
 # Copy results into the next year in G
 
 ((ip=$i+1))
 
 # $year denotes the respective simulation year in slice $i
 
 cp $pathg/slice$i/output/fesom.$year.ice.restart.nc $pathg/slice$ip/output/.
 cp $pathg/slice$i/output/fesom.$year.oce.restart.nc $pathg/slice$ip/output/.
 cp $pathg/slice$i/output/fesom.clock $pathg/slice$ip/output/.
 ...
done;
\end{lstlisting}

The parallel execution of the fine solver is given in Lst.\ref{lst:paraF}. To run the fine solver in parallel the ampersand is attached at the end of the \verb|srun| command. During the execution of the for-loop all FESOM2 simulations will be started without waiting for the completion of the previous task. The \verb|wait| command ensures that each parallel process is completed before starting the next instruction in the script.

\begin{lstlisting}[language=bash, caption=Extract from the Parareal SLURM Script: Parallel execution of the fine solver, label=lst:paraF]
for (( i=$TA ; i<=$TE ; i++ )); do
 ...
 
 # $TA and $TE define the first and last slice in the experiment
 # $pathf defines the absolute path to the fine folder
 # $flgs_f defines the amount of processes for spatial parallelization
 
 srun $flgs_f --chdir=$pathf/slice$i ./fesom.x > "$pathf/slice$i/fesom2.0.out" &
  
 ...
done;

wait
\end{lstlisting}

\subsection{Initial Value Handling}

Although the convergence of diagnostic variables is defined as goal of this study, we have to emphasize the particular importance of prognostic variables stored in restart files. The states local in time and space are sensitive to both spatial and temporal resolution. During the computation of the Parareal updates we found the jumps to exceed the predefined bounds for temperature and salinity, which will lead to model instabilities and eventually blow-ups. The predefined bounds are set by the FESOM2 developers to prevent the model to run into physically incorrect states, e.g. temperature below -4 degrees Celsius or salinity below 0 ppt. Due to the local differences of variables in the coarse and fine configurations or even by interpolation errors, we found these bounds violated during computation of new iterative solutions in restart files. The model is forced to terminate, even if just a single value in the entire ocean array is affected. To keep the simulations and consequently Parareal running, we have decided to correct such outliers and keep them within the defined limits. In Lst.\ref{lst:tempcheck} the function call for the correction in temperature arrays is given.

\begin{lstlisting}[language=Python,label=lst:tempcheck, caption=Extract of tempcheck.py: conservative interpolation function for element interpolation from fine to coarse mesh.]

def check_temp(var,ii):
	# ii :: horizontal layer index
			
	# DATA is read in from the fesom.yyyy.oce.restart.nc file
	x = DATA[var][0,:,ii]
				
	# For temperature -2.0 degrees Celsius 
	# is the set as the lower bound for all experiments.
	# all values smaller than temp_min are replaced.
	x[x < temp_min] = temp_min
			
	DATA[var][0,:,ii] = x
\end{lstlisting}

\section{Numerical Experiments}

In this section settings for the numerical experiments are introduced. The overall goal of this study was to investigate appropriate numbers of time slices that allow for parallelization in time by Parareal. We estimated the wall-clock time of each propagator on the Levante cluster in order to assess the theoretical speedups beforehand. At the time when the test cases were carried out Levante was only recently put into operation and thus, we observed considerable variability in the run-times for the coarse and fine propagators. On average, the FESOM2 requires 120 seconds to compute one year on the PI mesh and 430 seconds on the FPI mesh. With these estimates the run-time ratio is $m=3.6$ and would allow for two iterations at most in order to generate the minimal speedup of $S_2 = 1.2$. For the assessment of the convergence behavior we chose the maximum error norm $||\cdot||_\infty$ and relative error norm $||\cdot||_\infty / ||REF||_\infty$ of the approximations by Parareal with respect to the serial reference solution.

\subsection{Test Cases}

We distinguished the test cases by the their respective initial conditions into spin-up and restarted simulations. As listed in Tab.\ref{TAB:cases} we chose a time slice size of $\Delta T = 1$ year for all simulations. The time step size of both propagators is 36 spd since temporal refinement does not improve the approximation quality. Nevertheless, for the last test case we refined the time step size to investigate the impact on the stability of the Parareal algorithm. The simulations are distinguished into spin-up runs of 10 (case 1) and 20 years (case 2) simulation time. Further, we restarted computations from the Parareal solutions in Case 2 after 10 years to investigate the instabilities during the execution of FESOM2 in the second half of the 20 years time interval.

\begin{table}[H]
	\centering
	\caption{Numerical experiment settings for Parareal. The fine propagator is executed with $\Delta t_F =$ 36 spd for all configurations.}
	\begin{tabular}{|c|c|c|c|c|c|}
		\hline
		Case No. & $\Delta T$ & $T$ & Climatology IV & Restart IV & $\Delta t_C$ \\
		\hline
		1 & 1 year & 10 years & $\times$  & {} & 36 \\
		2 & 1 year & 20 years & $\times$  & {} & 36 \\
		\hline
		1.1 & 1 year & 10 years & {} & $\times$  & 36 \\
		1.2 & 1 year & 10 years & {} & $\times$  & 72 \\
		\hline
	\end{tabular}
	\label{TAB:cases}
\end{table}

\subsection{Results}

The convergence results presented in this section are restricted to a selection of diagnostic variables. FESOM2 simulations provide a vast amount of diagnostic variables and an investigation of those would extend the scope of this study. Similar to the procedure for mesh evaluation, we compiled a selection of diagnostics to analyze we considered appropriate in assessing convergence of the algorithm with respect to the ocean dynamics. The first test case will consider the annual mean temperature at various depths, AMOC at 26.5N and from south to north as well as an investigation of the errors in annual Salinity and Temperature fields before post-processing. \\

Before presenting the results it should be defined when they can be considered as converged, and accordingly what this is measured by. Inter-model comparisons of ocean circulation models result in differences in the simulation outcomes for the same time periods up to several degrees Celsius in annual mean temperature and Sverdrup in AMOC estimations. Hence, the approach to achieve convergence to double or even single precision was considered not suitable by the authors. Furthermore, when considering the sensibility of the FESOM2 solver regarding initial values, it was not possible to obtain convergence to machine precision in the first place. We decided that the absolute error between the reference and the iterated solution should be around $10^{-2}$ in order to consider the the respective test case as converged.

\subsubsection{Experiment 1}

\begin{figure}[H]
	\centering
	\includegraphics{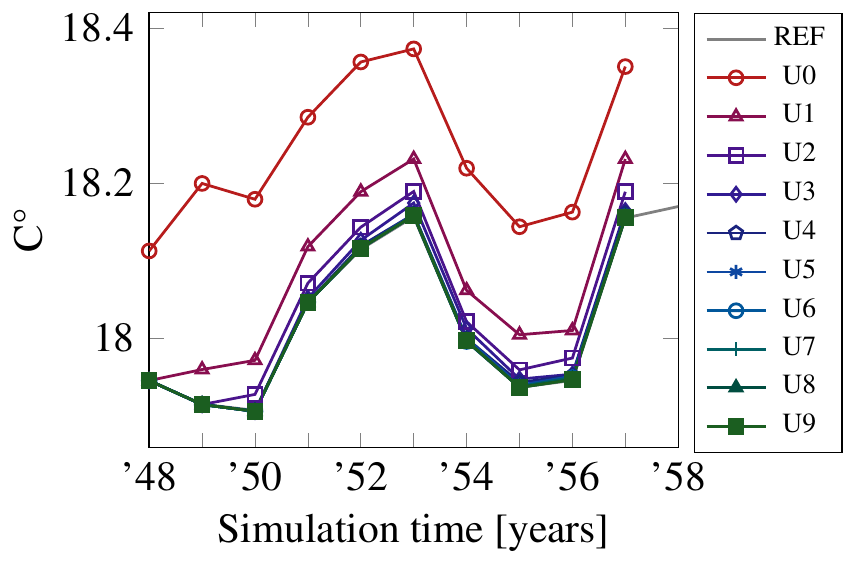}
	\caption{Annual sea surface temperature from 1948 to 1957 for the serial reference run and all iterations.}
	\label{FIG:SSTEXP1}
\end{figure}

In Fig.\ref{FIG:SSTEXP1} the annual sea surface temperature for the Parareal iterations and the serial reference run are given. At the surface an acceptable representation of the reference solution is observed with the second iteration already. 

\begin{figure}[H]
	\begin{subfigure}[t]{.48\textwidth}
			\includegraphics{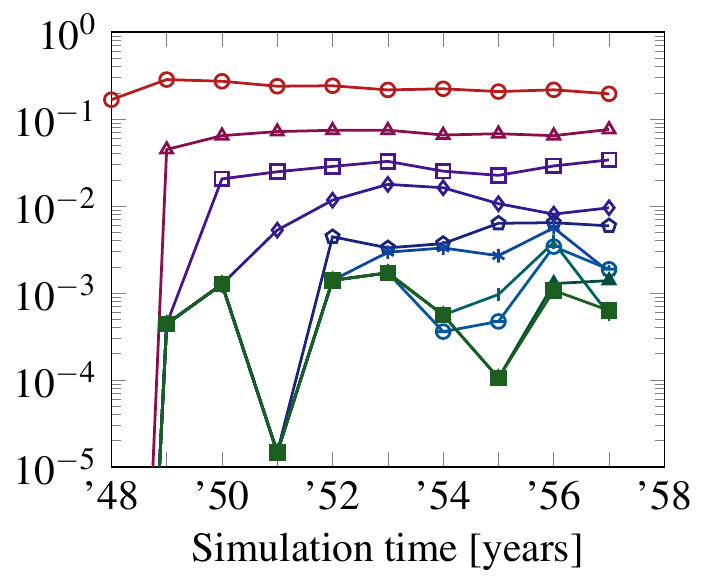}
		\end{subfigure}
	\begin{subfigure}[t]{.48\textwidth}
			\includegraphics{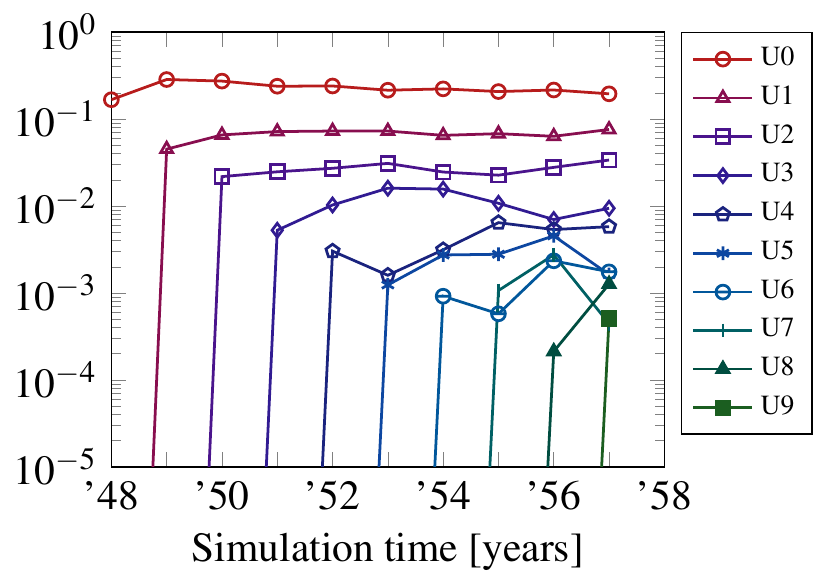}
		\end{subfigure}
	\caption{Absolute error of the annual mean surface temperature with respect to the reference solution (left) and the yearly restarted reference (right).}
	\label{FIG:SSTEXP1ERR}
\end{figure}

In the Fig.\ref{FIG:SSTEXP1ERR} two absolute error estimations are given. In the left panel the error with respect to the serial reference run is given. If FESOM2 is stopped after one year of computation and the then repeatedly restarted we found the solutions to slightly deviate between the serial runs. Since this is happening during the iteration procedure of Parareal, we took a repeated serial solution as a reference for the absolute error in the right panel of Fig.\ref{FIG:SSTEXP1ERR}. Since FESOM2 is a deterministic solver we find for increasing iteration number $k$ the error to be $10^{-16}$ in the time slices up to $k-1$. In the left panel the temperature difference vanishes in the first time slice of year 1948, since all computations in the first slice start from the same initial conditions. In both error estimations we observe the convergence to stagnate around $10^{-3}$. Taking into account the rather generous convergence criterion from the beginning of this section, the algorithm is assumed done after iterations. 

\begin{figure}[H]
	\centering
	\includegraphics{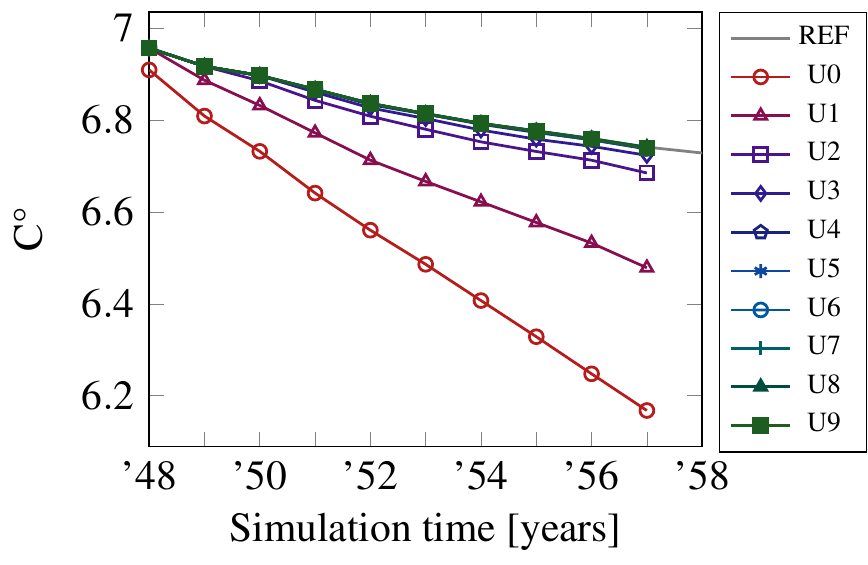}
	\caption{Annual mean temperature at 500m depth from 1948 to 1957.}
	\label{FIG:500TEMPEXP1}
\end{figure}

In Fig.\ref{FIG:500TEMPEXP1} the convergence behavior at 500m depth is given. In contrast to the surface temperature the difference between initial iteration $U_0$ and reference solution is increasing over time peaking in 1957 with 0.6°C. Although, the deviation between coarse and fine solution is larger, we find Parareal to reach convergence after 2 iterations. 

\begin{figure}[H]
	\begin{subfigure}[t]{.48\textwidth}
		\includegraphics{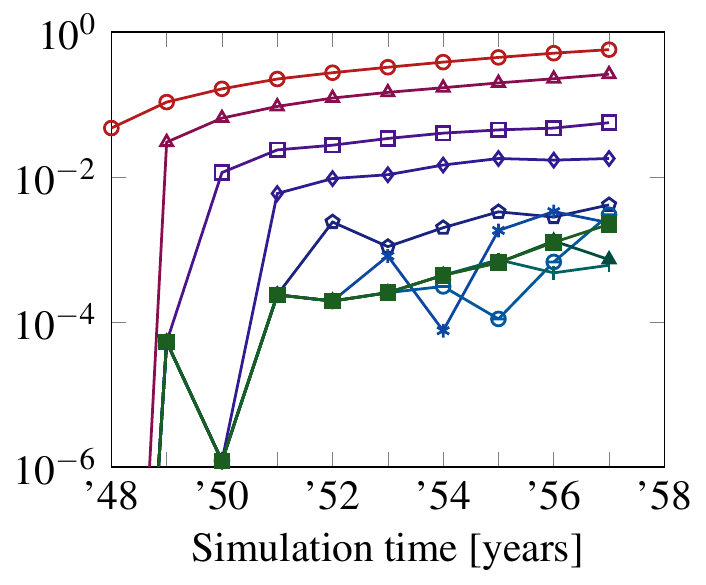}
	\end{subfigure}
	\begin{subfigure}[t]{.48\textwidth}
		\includegraphics{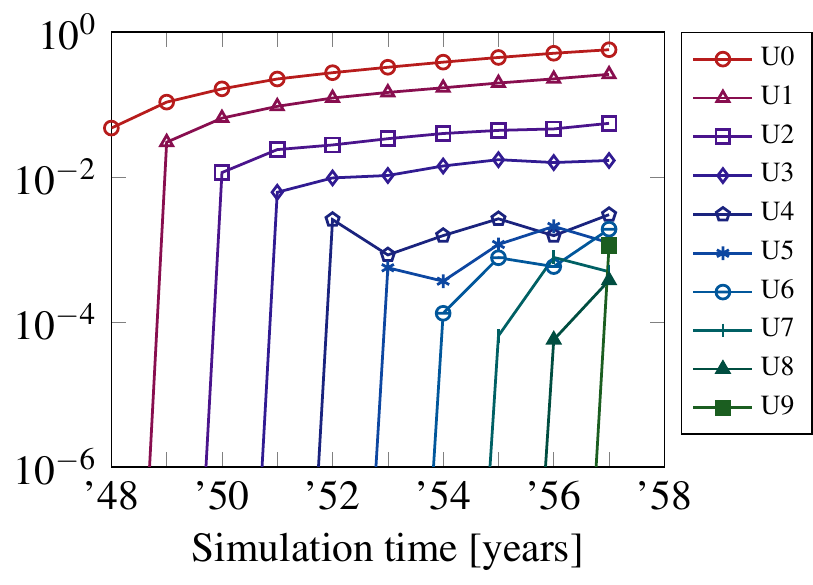}
	\end{subfigure}
	\caption{Absolute error of the annual mean temperature in 500m depth with respect to the reference solution (left) and the restarted fine solution (right).}
	\label{FIG:500TEMPEXP1ERR}
\end{figure}

In Fig.\ref{FIG:500TEMPEXP1ERR} an error reduction for all time slices can be observed up to fourth iteration. Beyond this point the convergence stagnates at $10^{-3}$, as observed for the sea surface temperature errors.

\begin{figure}[H]
	\centering
	\includegraphics{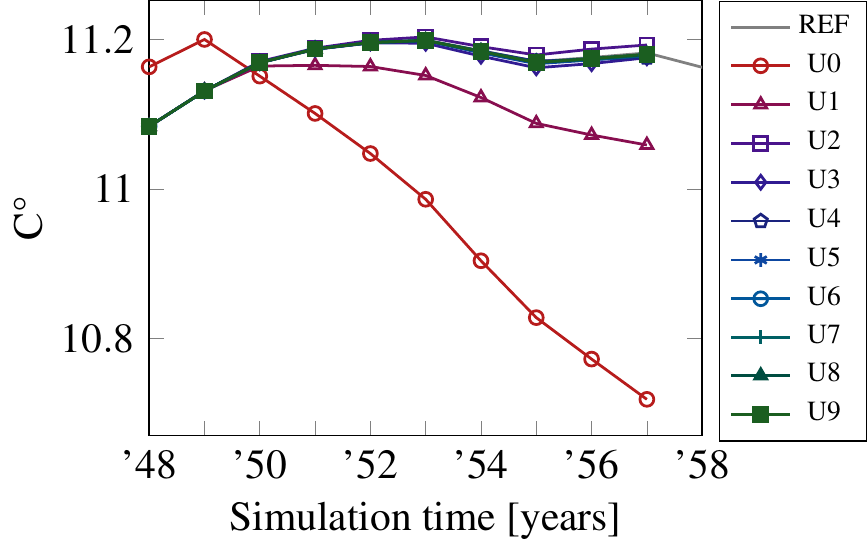}
	\caption{The annual mean temperature of the oceans averaged between surface to 500m depth.}
	\label{FIG:INTEMPEXP1}
\end{figure}

For the sake completeness of the investigation of the ocean's surface layers the spatially averaged temperature between sea surface and 500m depth is given in Fig.\ref{FIG:INTEMPEXP1}. A slightly better convergence is found with the computation of the mean over several layers. 

\begin{figure}[H]
	\begin{subfigure}[t]{.48\textwidth}
			\includegraphics{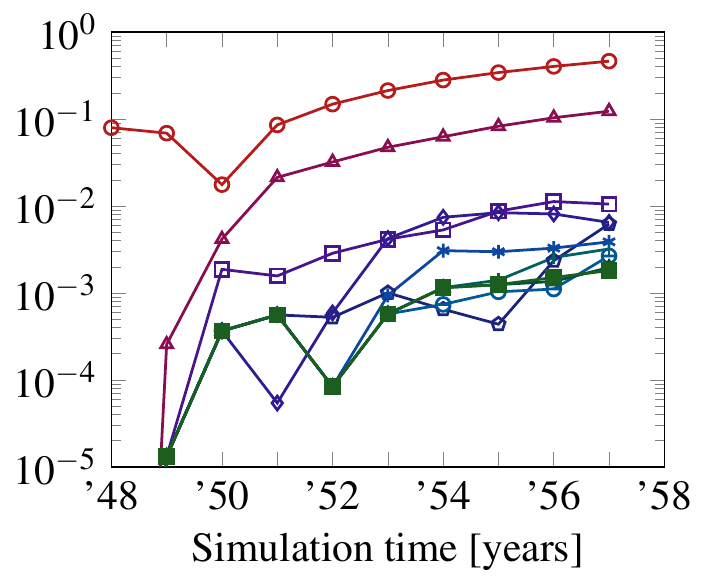}
		\end{subfigure}
	\begin{subfigure}[t]{.48\textwidth}
			\includegraphics{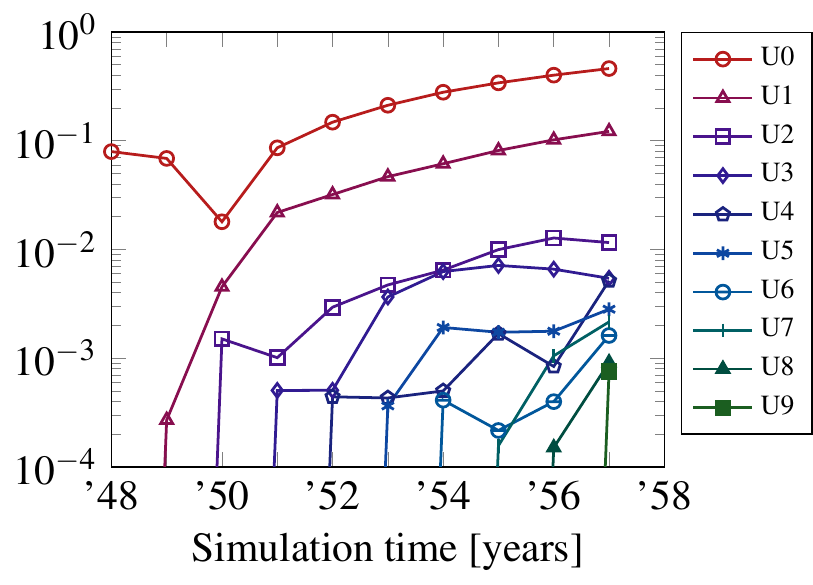}
		\end{subfigure}
	\caption{Absolute error of the annual mean temperature averaged between surface and 500m depth with respect to the reference solution (left) and the restarted fine solution (right).}
	\label{FIG:INTEMPEXP1ERR}
\end{figure}

In contrast to the former evaluations we find an improved error reduction in the first three iterations in Fig.\ref{FIG:INTEMPEXP1ERR} until the stagnation point is reached. It shows that with additional averaging in vertical direction that the overall convergence of the time interval improves. 

We conclude for the annual mean temperature that an error threshold around $10^{-3}$ in the absolute error cannot be overcome. The sensitivity to initial values in FESOM2 simulations, the impact of the interpolation methods and the manipulation of restart files during the update procedure of Parareal seem to form an obstacle that cannot be overcome. \\

\begin{figure}[H]
	\centering
	\includegraphics{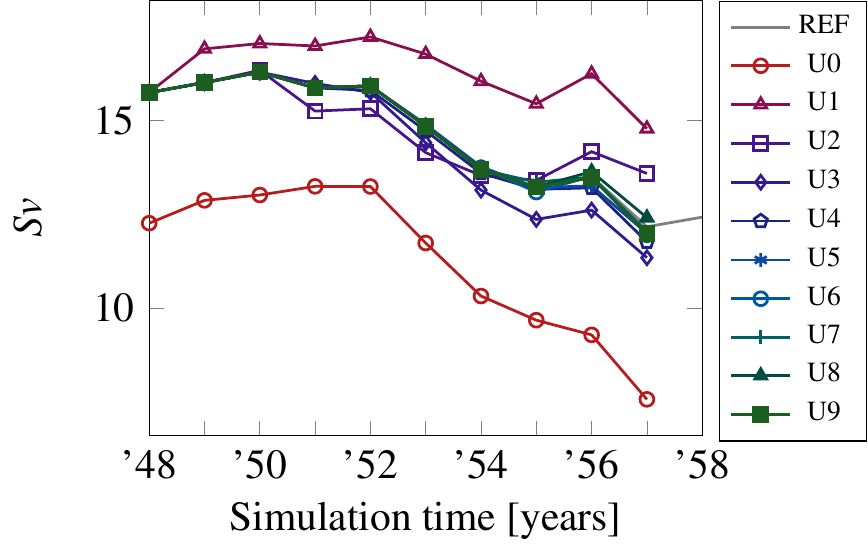}
	\caption{Parareal iterations of the AMOC at 26.5N.}
	\label{FIG:AMOCEXP1}
\end{figure}

In Fig.\ref{FIG:AMOCEXP1} the convergence to the AMOC at 26.5N is depicted. The initial approximation by the first iteration $U_0$ is already around 3 Sverdrup below the reference representation of the AMOC. In accordance with the convergence for temperature we find the second iteration to be sufficiently accurate.

\begin{figure}[H]
	\begin{subfigure}[t]{.48\textwidth}
			\includegraphics{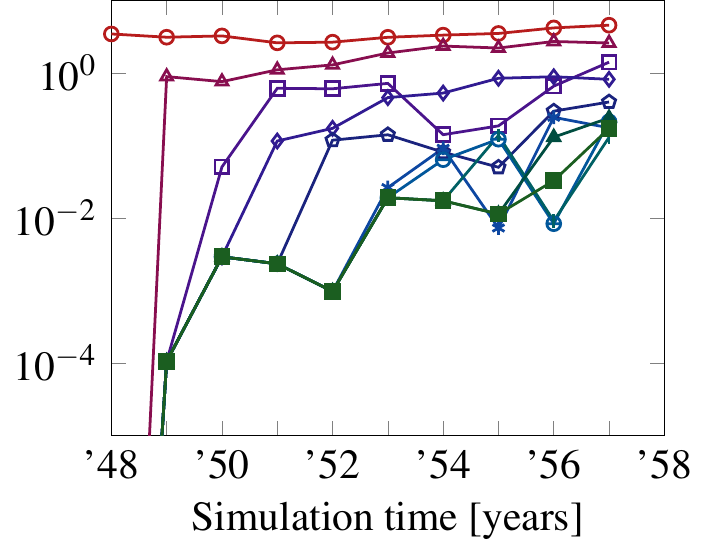}
		\end{subfigure}
	\begin{subfigure}[t]{.48\textwidth}
			\includegraphics{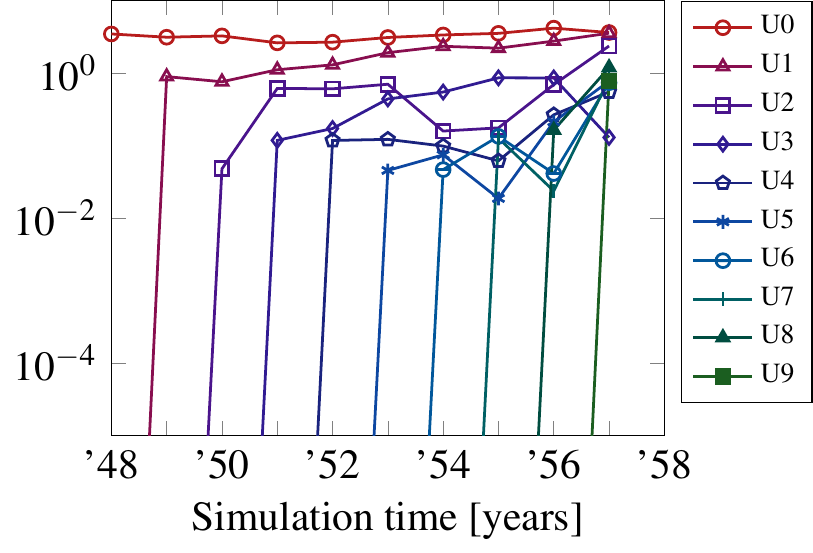}
		\end{subfigure}
	\caption{Absolute error of iterative solutions to the reference AMOC at 26.5N (left) and the repeated serial computation (right).}
	\label{FIG:AMOCEXP1ERR}
\end{figure}

In the absolute error evaluations in Fig.\ref{FIG:AMOCEXP1ERR} the distinct error reduction between consecutive iterations breaks at the second Parareal approximation $U_2$. In agreement with the previous results an error stagnation is found, for the AMOC around $10^{-1}$. Again a rather acceptable convergence behavior in the profile in Fig.\ref{FIG:AMOCEXP1} is contrasted by the inability to reduce errors to double or even single precision, shown in Fig.\ref{FIG:AMOCEXP1ERR}.  

From the point of view of the FESOM2 users and the expectations regarding the accuracy of the simulation results, the application for the AMOC can be considered a success. If, on the other hand, the algorithm is evaluated on the basis of numerical accuracy, it is difficult to assume a successful implementation. \\

\begin{figure}[H]
	\centering
	\begin{subfigure}[t]{.49\textwidth}
			\includegraphics[width=\textwidth]{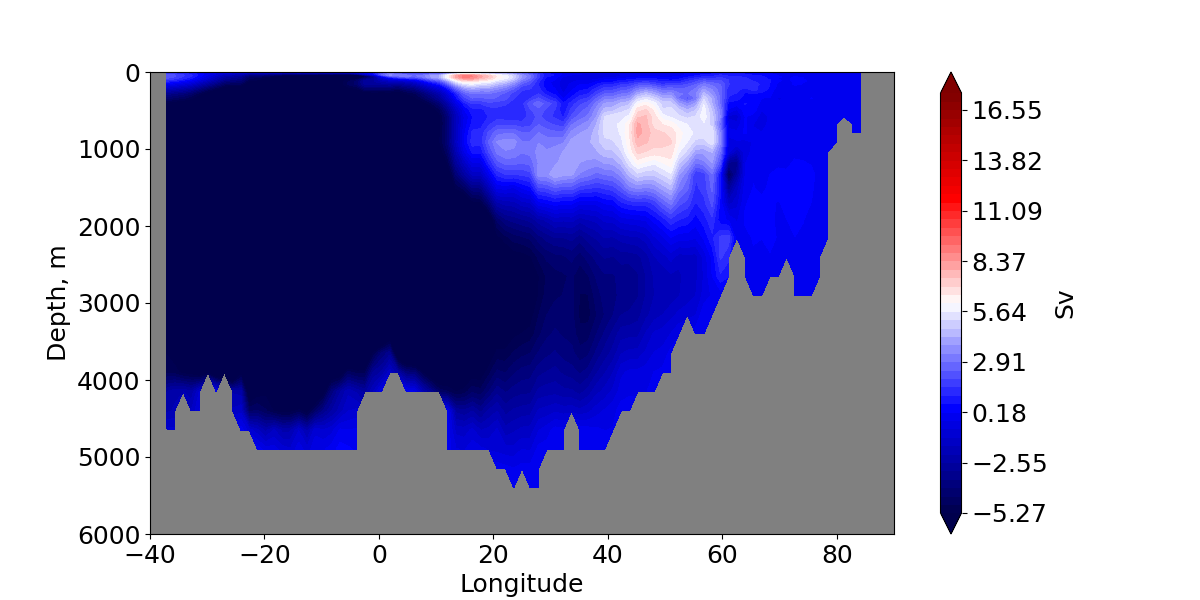}
			\caption{Iteration 0}
		\end{subfigure}
	\hfill
	\begin{subfigure}[t]{.49\textwidth}
			\includegraphics[width=\textwidth]{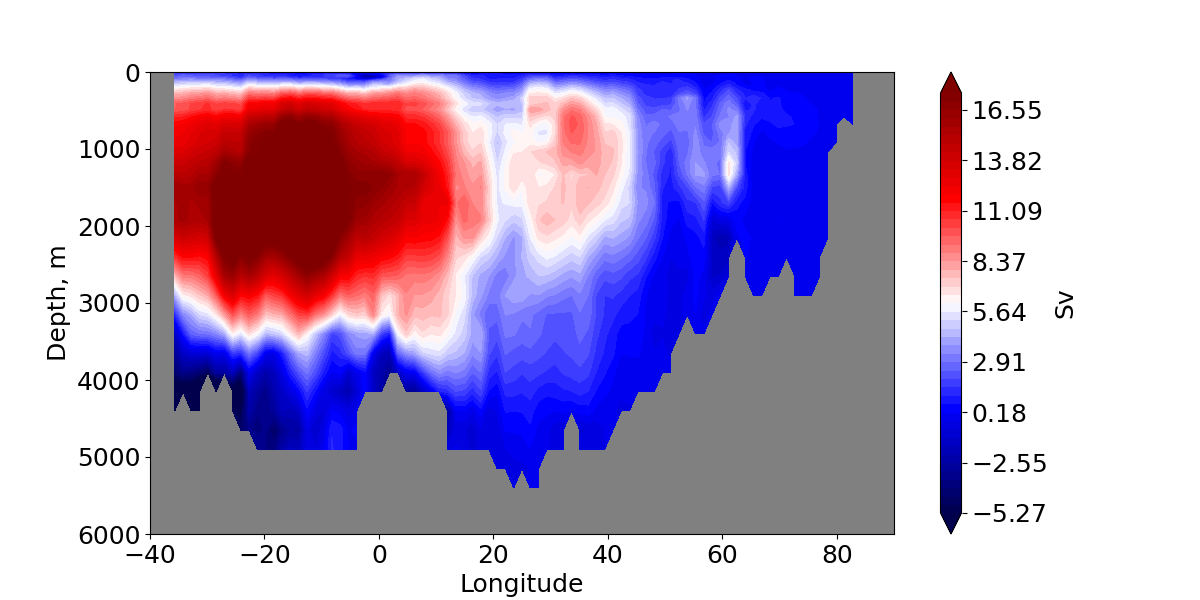}
			\caption{Difference to reference.}
		\end{subfigure}
	\begin{subfigure}[t]{.49\textwidth}
			\includegraphics[width=\textwidth]{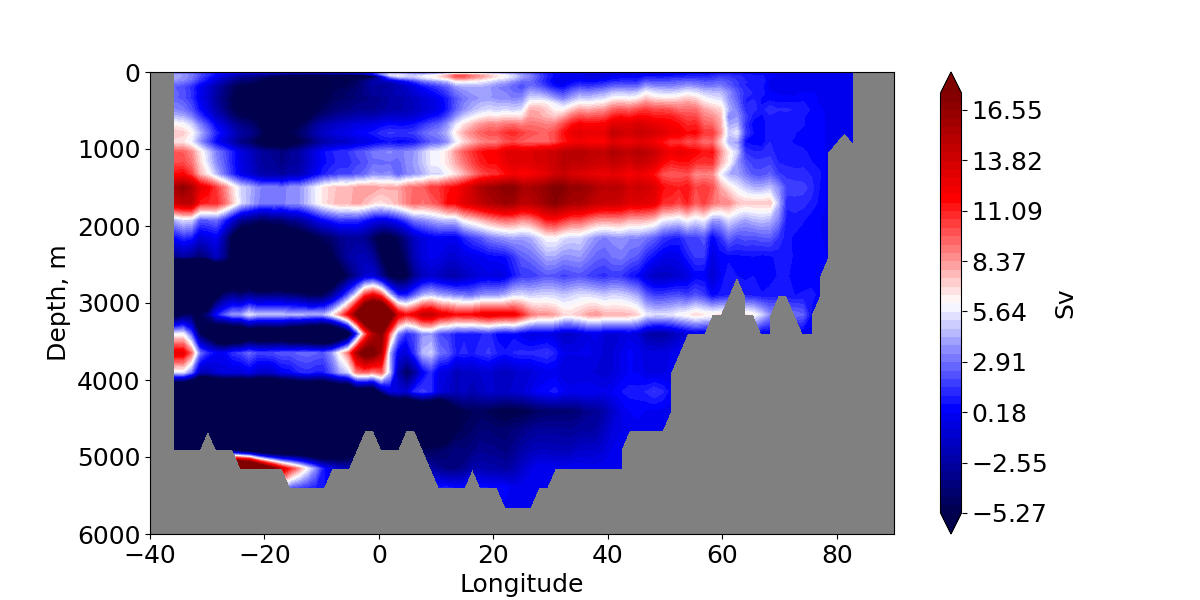}
			\caption{Iteration 1}
		\end{subfigure}
	\hfill
	\begin{subfigure}[t]{.49\textwidth}
			\includegraphics[width=\textwidth]{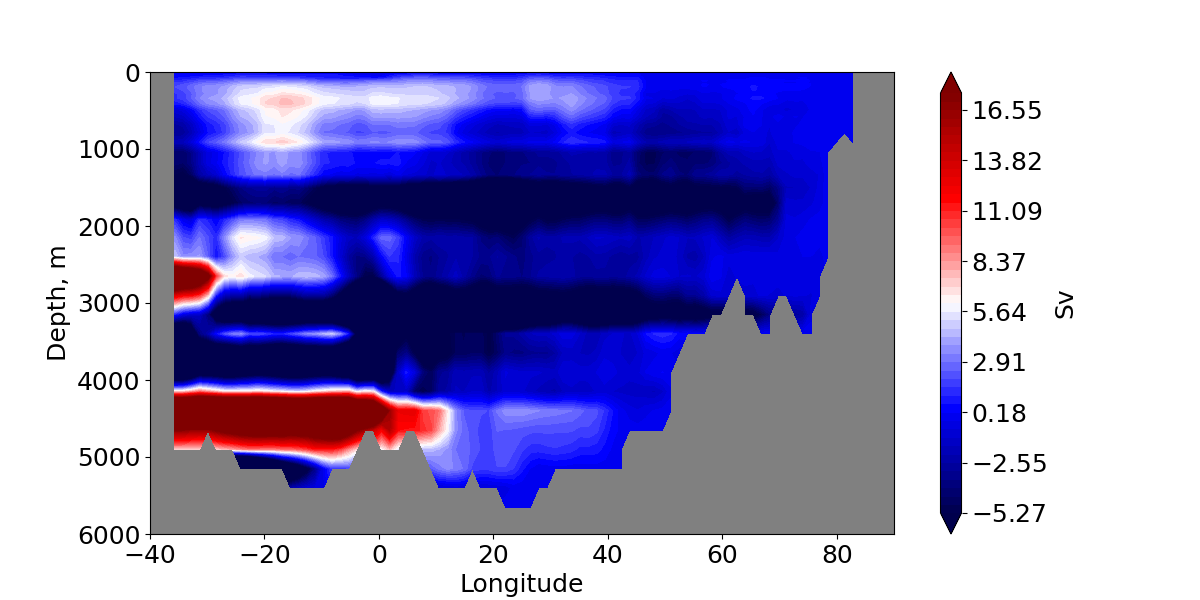}
			\caption{Difference to reference.}
		\end{subfigure}
	\begin{subfigure}[t]{.49\textwidth}
			\includegraphics[width=\textwidth]{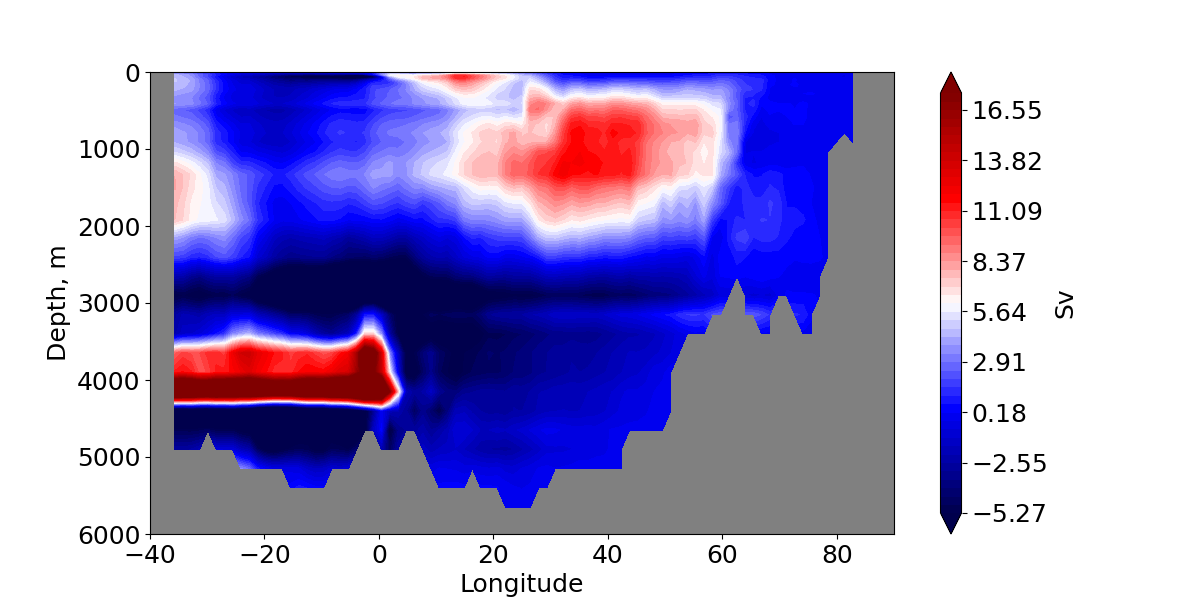}
			\caption{Iteration 3}
		\end{subfigure}
	\hfill
	\begin{subfigure}[t]{.49\textwidth}
			\includegraphics[width=\textwidth]{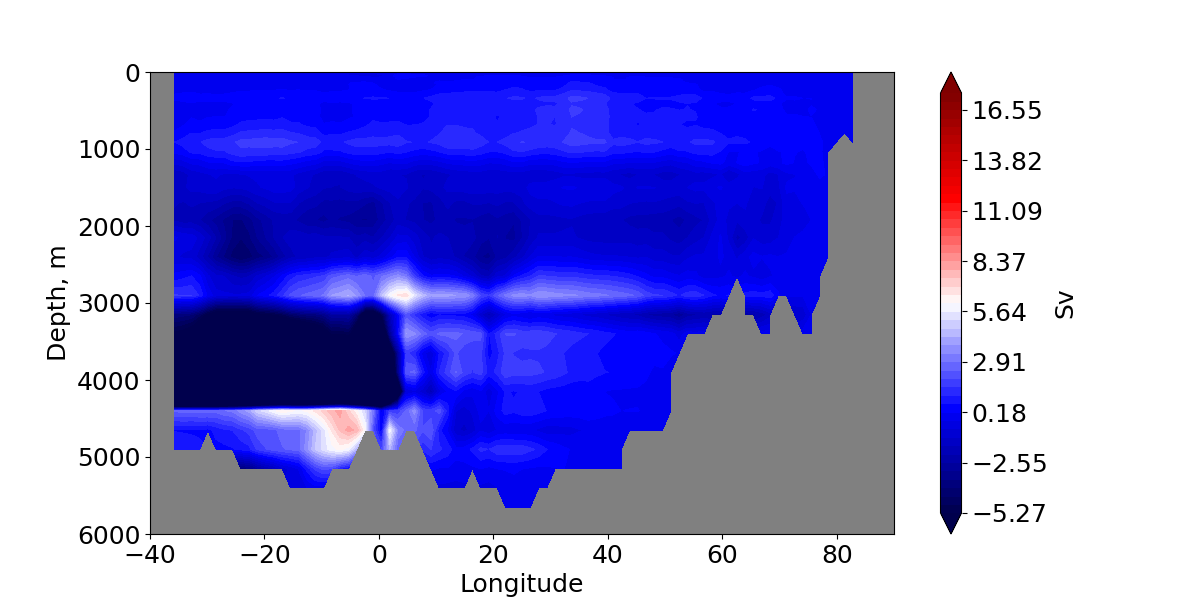}
			\caption{Difference to reference.}
		\end{subfigure}
	\begin{subfigure}[t]{.49\textwidth}
			\includegraphics[width=\textwidth]{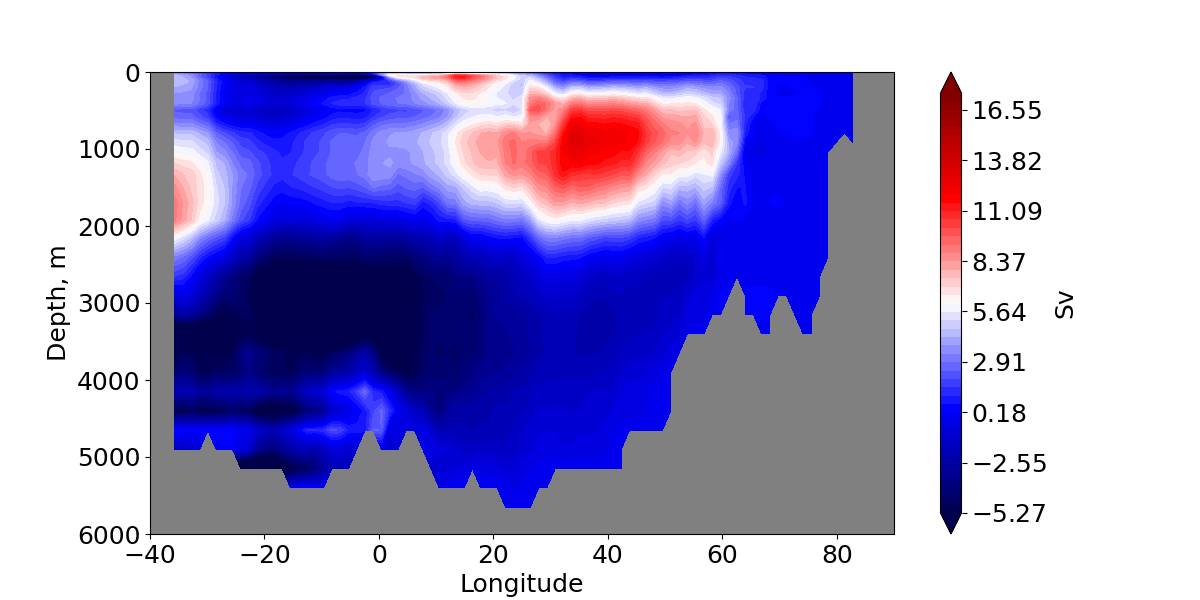}
			\caption{Iteration 5.}
		\end{subfigure}
	\hfill
	\begin{subfigure}[t]{.49\textwidth}
			\includegraphics[width=\textwidth]{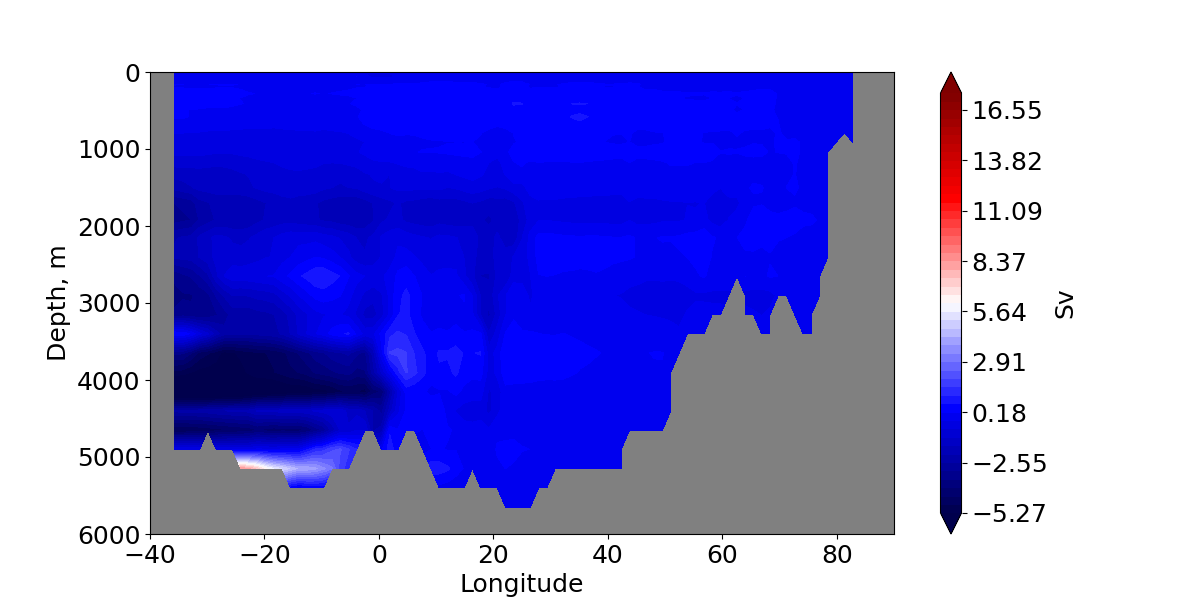}
			\caption{Difference to reference.}
		\end{subfigure}
	\begin{subfigure}[t]{.49\textwidth}
			\includegraphics[width=\textwidth]{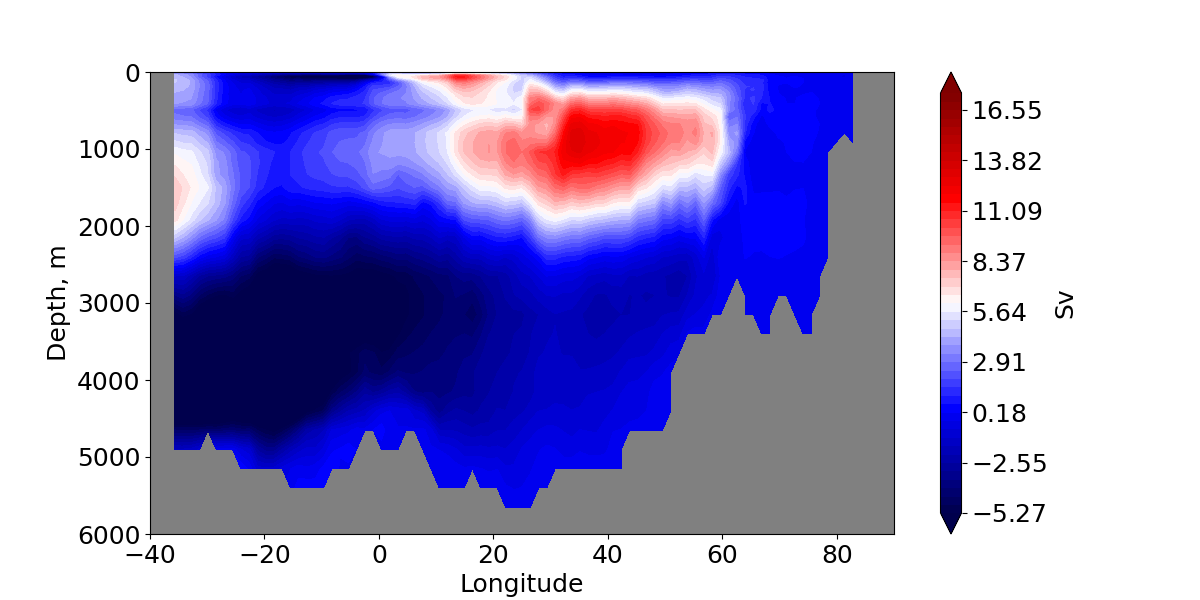}
			\caption{Iteration 7.}
		\end{subfigure}
	\hfill
	\begin{subfigure}[t]{.49\textwidth}
			\includegraphics[width=\textwidth]{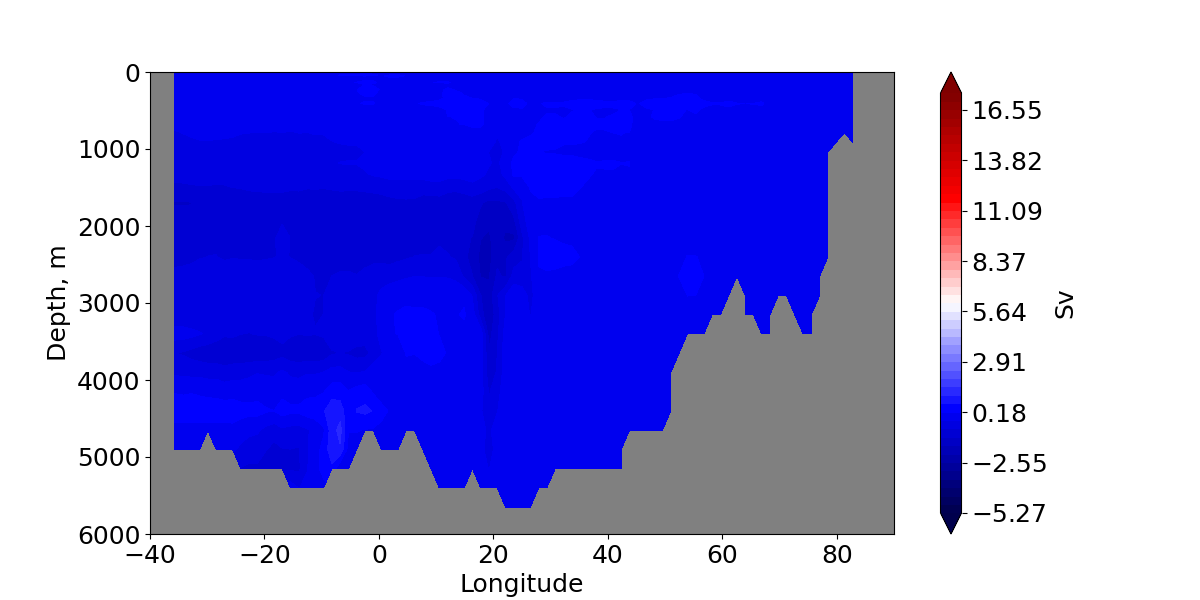}
			\caption{Difference to reference.}
		\end{subfigure}
	\caption{Parareal iterations of the 2-D AMOC visualization evaluated from 40S to 80N.}
	\label{FIG:AMOC2DEXP1}
\end{figure}

The second investigation of the AMOC is evaluated from 40S to 80N and in depth resulting in a 2D representation in Fig.\ref{FIG:AMOC2DEXP1}. We observe the initial approximation $U_0$ to fail entirely in solving the overturn circulation. As mentioned in the mesh evaluation part earlier on, the PI mesh's coarse spatial resolution reveals its inability to simulate the AMOC properly. In fact, the total break down of the ocean circulation in the ocean leads to an approximation close to zero Sverdrup in magnitude. It takes at least 5 iterations to recover the AMOC profile of the reference simulation. \\

Investigating the North-South evaluation of the AMOC underlines the necessity to consider as many aspects of the Atlantic as possible. The convergence in the AMOC at 26.5N in Fig.\ref{FIG:AMOC2DEXP1} alone would have been misleading in the assessment of Parareal. Considering 5 iterations as sufficiently accurate would still exceed the allowed number of iterations as formulated a priori for the theoretical speed up. 

\begin{figure}[H]
	\begin{subfigure}[t]{.45\textwidth}
		\includegraphics{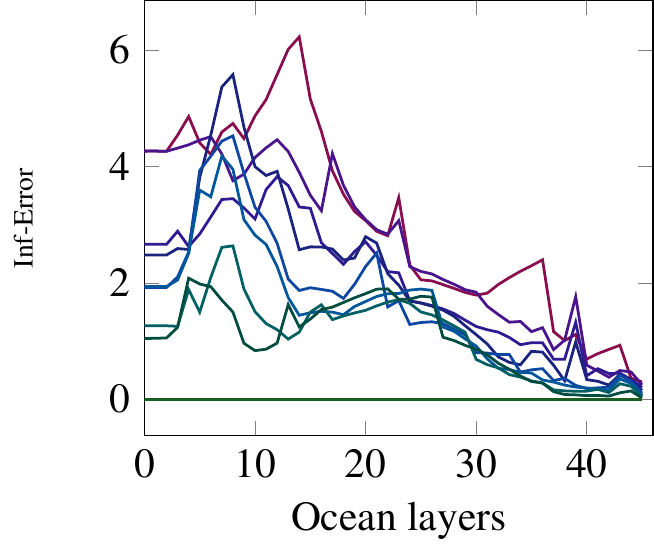}
	\end{subfigure}
	\begin{subfigure}[t]{.45\textwidth}
		\includegraphics{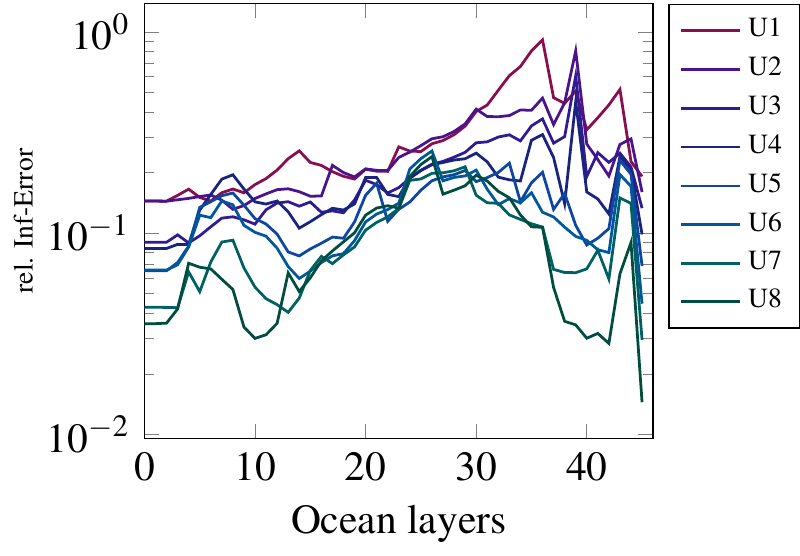}
	\end{subfigure}
	\caption{Parareal error reduction in the time-averaged horizontal temperature field over the ocean layers at the last slice. Layer 0 corresponds to the ocean surface.}
	\label{FIG:TEMPFIELDEXP1}
\end{figure}

To conclude the evaluation of the first experiment, the fields of temperature and salinity are considered before post-processing. Both data sets are time-averaged, and therefore diagnostic, but local in space. In Fig.\ref{FIG:TEMPFIELDEXP1} the absolute and relative maximum error of the temperature for each horizontal layer in the ocean is given. The fields have been evaluated in the last slice in 1957. We excluded the initial iteration from the investigation since it was computed solely on the PI mesh, and therefore  additional errors were introduced by lifting to the FPI mesh. The maximum error in the upper is generally found to be significantly larger than towards the bottom of the ocean. Due to slow dynamics in the so-called deep ocean changes take place over long time scales. Hence, the error is smaller and reaches absolute values close to zero while at the layers close to the surface we find differences up to 6°C. For the relative error, shown in the right panel of Fig.\ref{FIG:TEMPFIELDEXP1}, the deep layers show a slightly higher deviation for the first iterations. In every horizontal level the relative error was computed with respect to its maximum value its respective horizontal array. While temperatures on the surface can reach up to 30°C, the temperature in the deep ocean is only around 4 degrees. Accordingly, the relative error for all levels is about $10^{-1}$. 

By post-processing the time-averaged temperature fields and obtaining the annual mean temperature we seem to filter out the maximum errors in the respective horizontal levels. The error at surface level is for $U_1$ around 4 °C, while for the post-processed temperature in Fig.\ref{FIG:SSTEXP1} we see an reduction to 0.2°C. \\

\begin{figure}[H]
	\begin{subfigure}[t]{.45\textwidth}
		\includegraphics{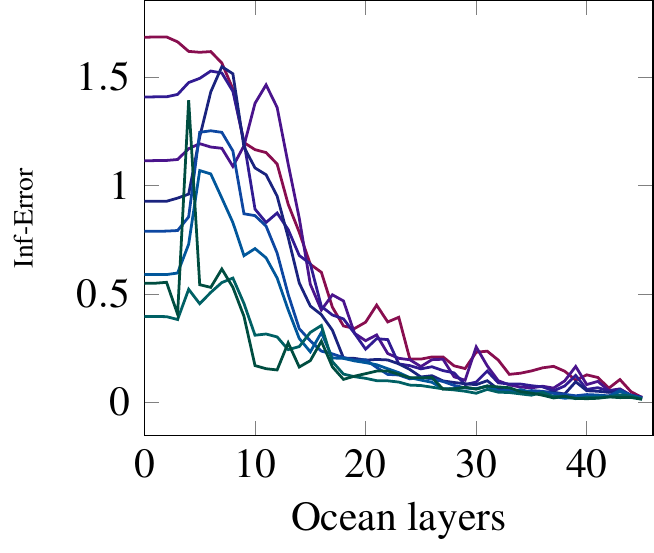}
	\end{subfigure}
	\begin{subfigure}[t]{.45\textwidth}
		\includegraphics{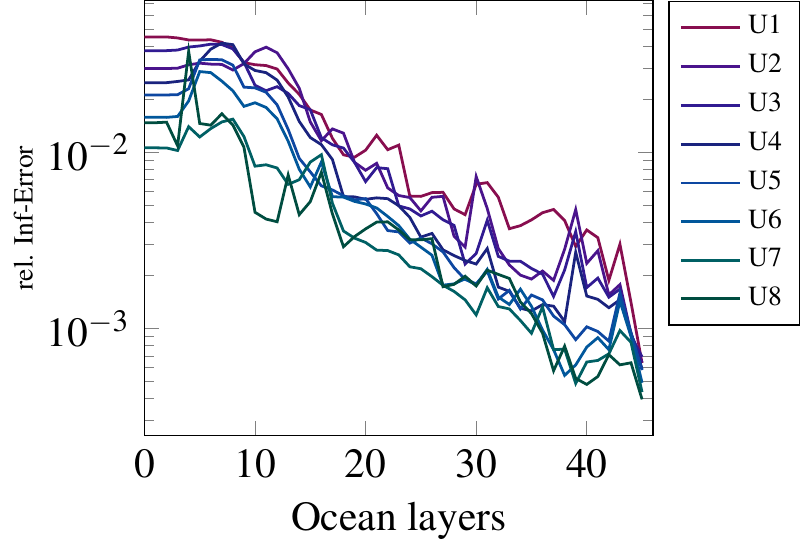}
	\end{subfigure}
	\caption{Parareal error reduction in the time-averaged horizontal salinity fields over the ocean layers at the last slice. Layer 0 corresponds to the ocean surface.}
	\label{FIG:SALTFIELDEXP1}
\end{figure}

For the salinity we restricted the convergence results to the annual averaged salinity local in space in the same manner we investigated the temperature fields before. In Fig.\ref{FIG:SALTFIELDEXP1} the maximum and relative maximum error is depicted. For both errors over the ocean layers the tendency to decrease with depth is given. The error reduction by Parareal takes mainly part in the first twenty levels. In the deep ocean the deviations from iteration 1 onwards are already very small.

\subsubsection{Experiment 2}

The long-term experiment over 20 years simulation time failed with the second iteration. In the results given one can find two different types of failures in the algorithm. Occurring blow-ups during the execution of FESOM2 will lead to the problem, that no output files are generated. If an output file is missing for a specific operation in the algorithm, the SLURM script will omit the respective step. Accordingly, it makes a difference if the blow-up occurs during the execution of the fine or coarse run. In case of a failure during the fine execution, no update can be performed and this operation will be skipped. If the coarse run is incomplete no restart files are generated and FESOM2 will generate a initial condition from climatology for the consecutive time slice. Both failure types can be found during the long-term run. \\

\begin{figure}[H]
	\centering
	\includegraphics{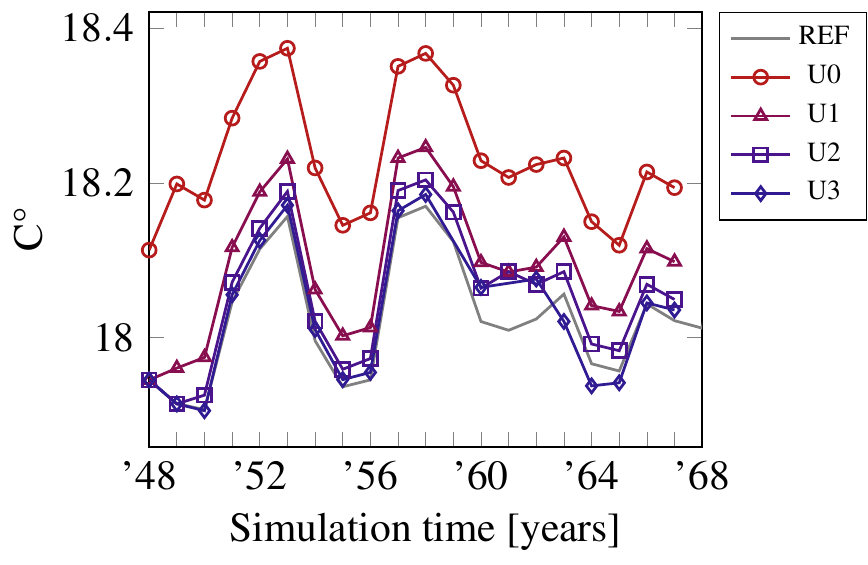}
	\caption{Annual sea surface temperature from 1948 to 1967 for the serial reference run and all iterations.}
	\label{FIG:SSTEXP2}
\end{figure}

In Fig.\ref{FIG:SSTEXP2} the annual mean sea surface temperature is shown over 20 years of simulation from 1948 to 1967. During the computation of the second iteration $U_2$ the fine propagator suffered from instabilities during the year 1961. Thus, no output files were generated and the update procedure could not be executed and accordingly no fine contribution was added. Therefore the temperature in 1961 is the same as for the first iteration $U_1$.
In iteration $U_3$ values are missing for the years 1959 and 1961, where the coarse propagator failed to complete the computation without error. The updated iteration is written into the lifted files of the coarse run. With this file missing, no results for the respective year can be processed. 

Apart from the instabilities during the execution, the convergence tendency during the first two iterations seems promising. Nevertheless, we must consider the application of Parareal to longer simulation intervals as a failure.

\begin{figure}[H]
	\centering
	\includegraphics{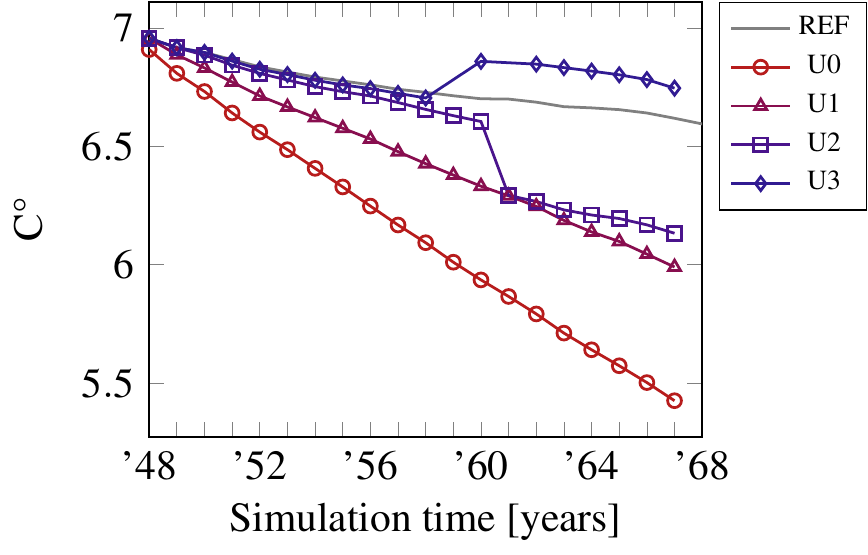}
	\caption{The annual mean temperature of the oceans averaged between surface to 500m depth.}
	\label{FIG:500TEMPEXP2}
\end{figure}

\begin{figure}[H]
	\centering
	\includegraphics{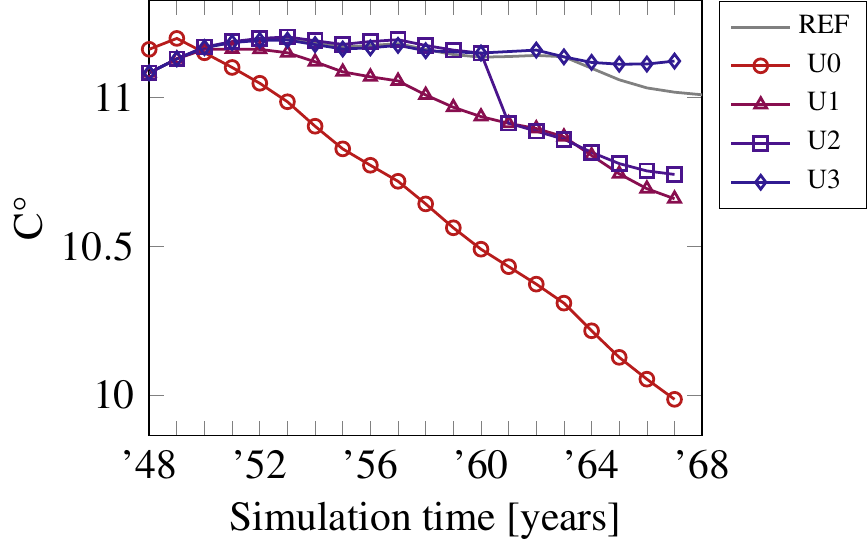}
	\caption{Annual mean temperature at 500m depth from 1948 to 1967.}
	\label{FIG:INTEMPEXP2}
\end{figure}

Consequently, the temperature profiles at 500m depth in Fig.\ref{FIG:500TEMPEXP2} and averaged between surface and 500m depth in Fig.\ref{FIG:INTEMPEXP2} suffer from the blow-ups arising during the execution of FESOM2 on both meshes. In Fig.\ref{FIG:500TEMPEXP2} one can observe best how the failure in the respective propagators affects the evolution of the mean temperature. In absence of results by the fine solver the algorithm proceeds from the previous iteration. In case of a coarse propagator failure no output files are provided and FESOM2 will start by deriving initial states from climatalogy data sets, shifting the profile towards an observed (by experiment) state.

\begin{figure}[H]
	\centering
	\includegraphics{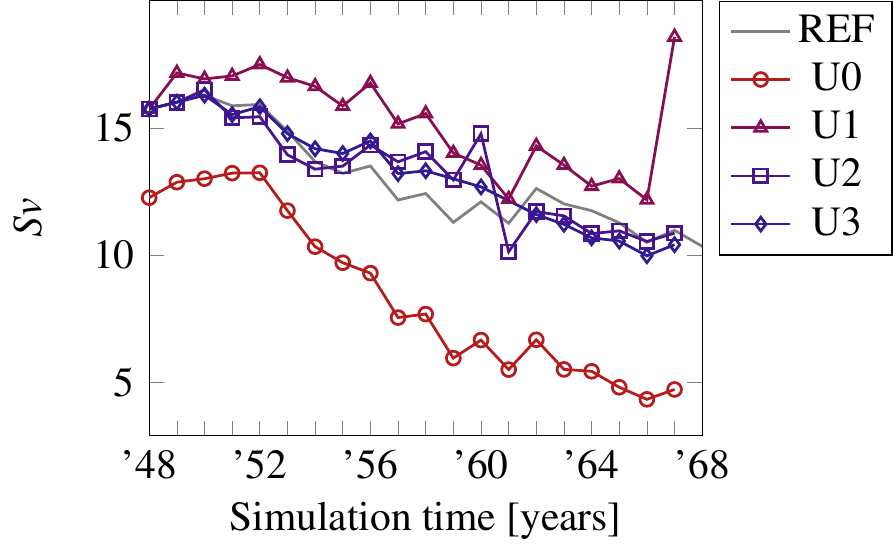}
	\caption{First three iterations for the AMOC at 26.5N.}
	\label{FIG:AMOCEXP2}
\end{figure}

Since the iteration is incomplete, we consider an evaluation of the absolute and relative errors to be of little use. In agreement with the first iterations for the annual mean temperature the AMOC at 26.5N is given in Fig.\ref{FIG:AMOCEXP2}. 

\begin{figure}[H]
	\begin{subfigure}[t]{.45\textwidth}
			\includegraphics{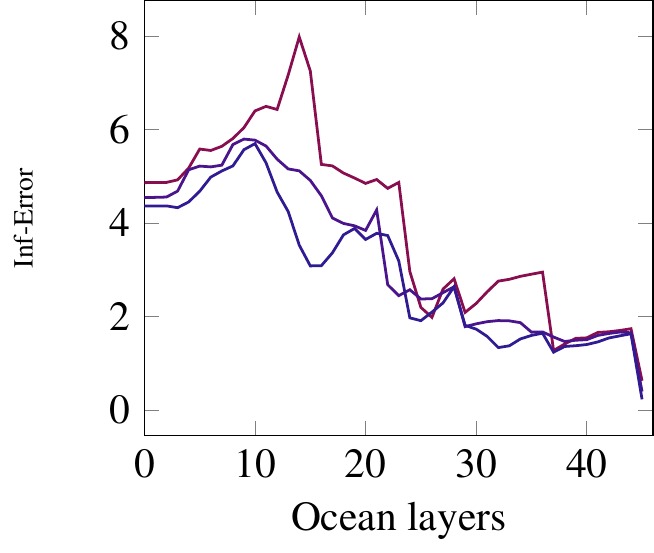}	
		\end{subfigure}
	\begin{subfigure}[t]{.45\textwidth}
			\includegraphics{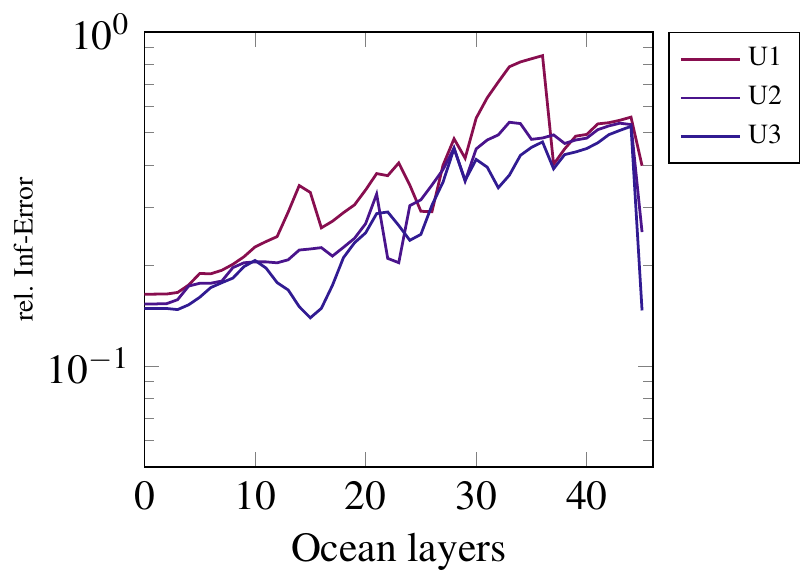}
		\end{subfigure}
	\caption{Parareal error reduction in the time-averaged horizontal temperature field over the ocean layers at the last slice. Layer 0 corresponds to the ocean surface.}
	\label{FIG:TEMPFIELDEXP2}
\end{figure}

\begin{figure}[H]
	\begin{subfigure}[t]{.45\textwidth}
			\includegraphics{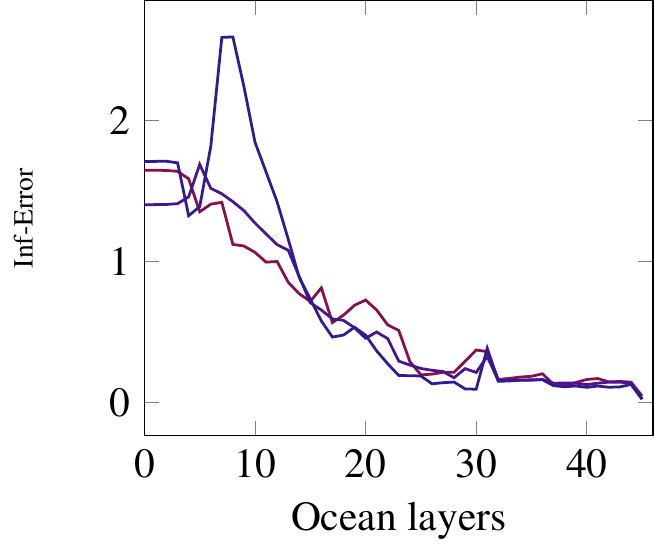}
		\end{subfigure}
	\begin{subfigure}[t]{.45\textwidth}
			\includegraphics{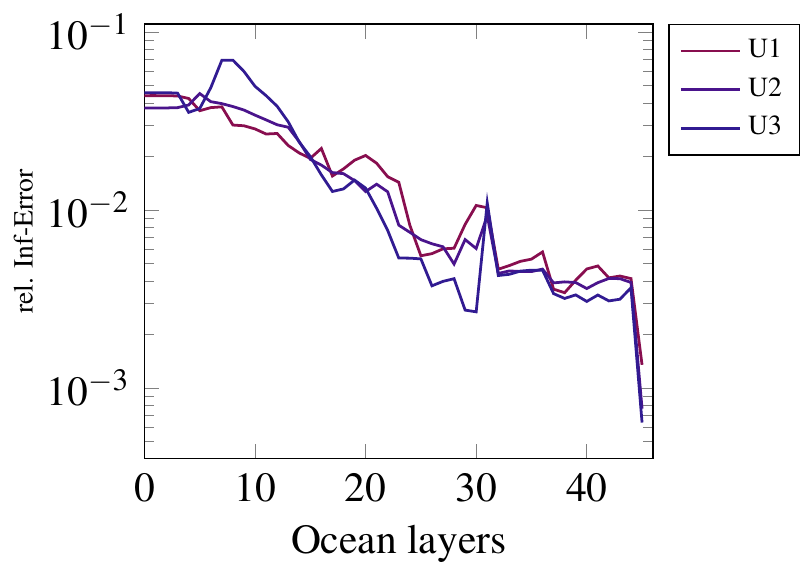}
		\end{subfigure}
	\caption{Parareal error reduction in the time-averaged horizontal salinity fields over the ocean layers at the last slice. Layer 0 corresponds to the ocean surface.}
	\label{FIG:SALTFIELDEXP2}
\end{figure}

In Figs.\ref{FIG:TEMPFIELDEXP2} and \ref{FIG:SALTFIELDEXP2} the absolute and relative errors of the horizontal layers are given for temperature and salinity. In comparison to the first test case the error build up in the last slice is significantly larger. For time-averaged temperature the maximum error reaches up to 8°C and for salinity 2.5g/kg, respectively. During the analysis of error log files by FESOM2 we found local violations of the vertical CFL number and values outside of predefined ranges for temperature and salinity. The update steps of Parareal algorithm were implemented in a way that would conserve quantities, e.g. energy. Hence, the operations carried out during the update procedure lead to an error build up that eventually causes FESOM2 simulations to diverge. 

\subsubsection{Experiment 3}

\begin{figure}[H]
	\centering
	\includegraphics{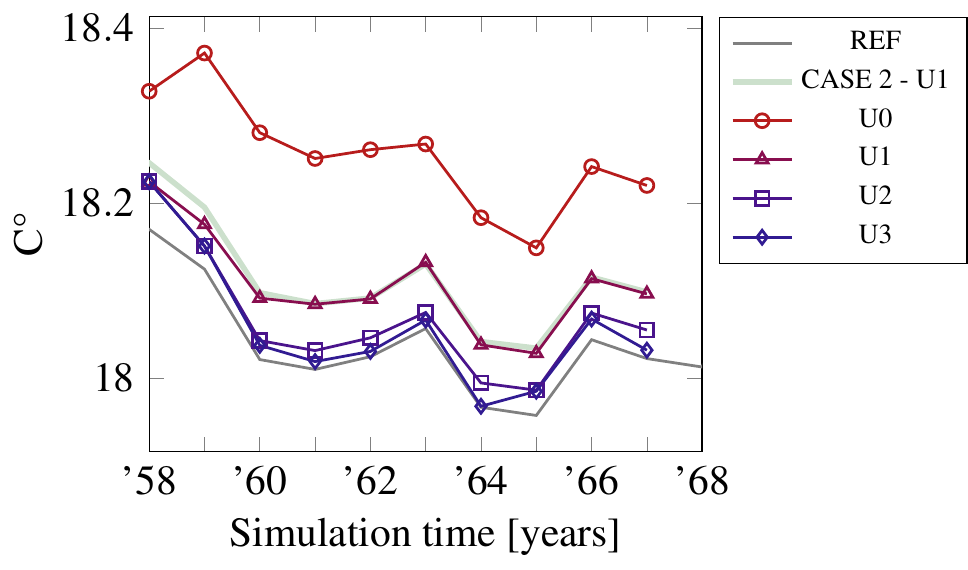}
	\caption{Sea surface temperature from 1958 to 1967. The results from 1957 in the first iteration $U_1$ from experiment 2 were used as initial value in slice 1.}
	\label{FIG:SSTEXP3}
\end{figure}

Since the application of Parareal fails to converge for 20 years simulation a restart of the algorithm is performed after the completion of first 10 years in test case 2. Unfortunately, we needed the algorithm to start from the result of iteration 1. For all other Parareal approximations of higher iteration number in case 2 the execution will suffer from blow-ups right away in the first iteration. We were able to obtain three iterations until the error build-up became a problem and caused instabilities, again. Therefore, the results are restricted to the annual mean temperature at the surface and in 500m depth.
With the choice of the initial value being the first iteration the best approximation obtainable appears to be limited. Exception is the sea surface temperature in Fig.\ref{FIG:SSTEXP3}, where the interface to the forcing data set is located, providing accurate input, e.g. from atmosphere, fresh water fluxes, etc. \\

In Fig.\ref{FIG:SSTEXP3} the convergence in the first three iterations for the sea surface temperature is shown. The first iteration $U_1$ recovers the first iteration of the previous experiment. With $U_2$ and $U_3$ the algorithm converges towards the serial reference solution. 

\begin{figure}[H]
	\centering
	\includegraphics{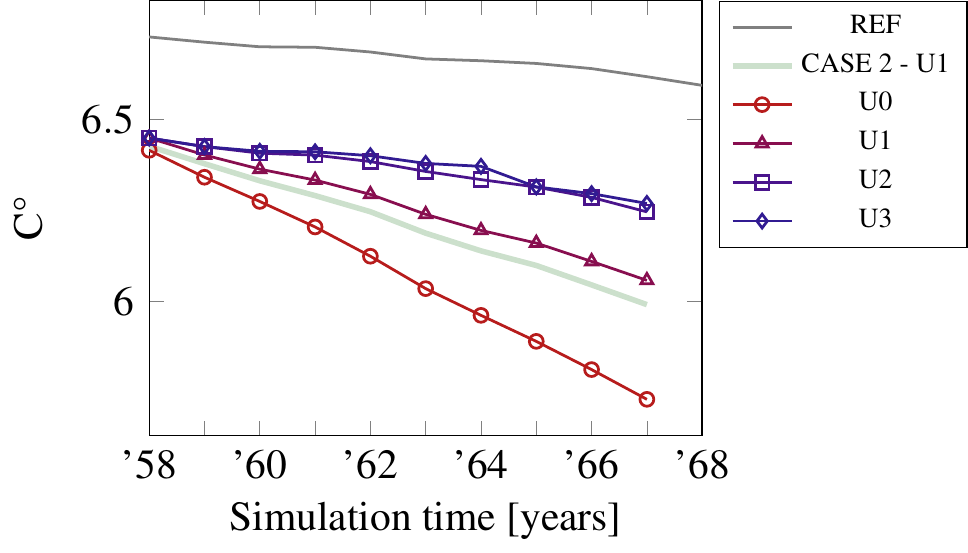}
	\caption{Annual mean temperature from 1958 to 1967 at 500m depth.}
	\label{FIG:500TEMPEXP3}
\end{figure}

In contrast to the sea surface temperature we observe the convergence in 500m depth stagnating with the second iteration $U_2$, compare Fig.\ref{FIG:500TEMPEXP3}. Without the immediate interface to the forcing data set the outcome is dominantly dependent on the initial values provided by the previous test case. Hence, the approximation cannot approach the serial reference in the way observed in the other experiments. 

\subsubsection{Experiment 4}

In an attempt to tackle the stability problems arising in the first test case for restarted Parareal the time step size of the coarse solver was increased from 36 to 72 spd. According to the sensibility of diagnostics towards time step size (ref to section) and the FESOM2 documentation (ref) we assumed that increasing the time step size improves the stability during the execution of FESOM2, but has no impact to time-averaged quantities. Hence, we expected the Parareal iteration to become more stable and the results presented in this section seem to back up this approach. We chose iteration 3 of case 2 to restart Parareal in the year 1958 and were able to compute nine iterations. \\

\begin{figure}[H]
	\centering
	\includegraphics{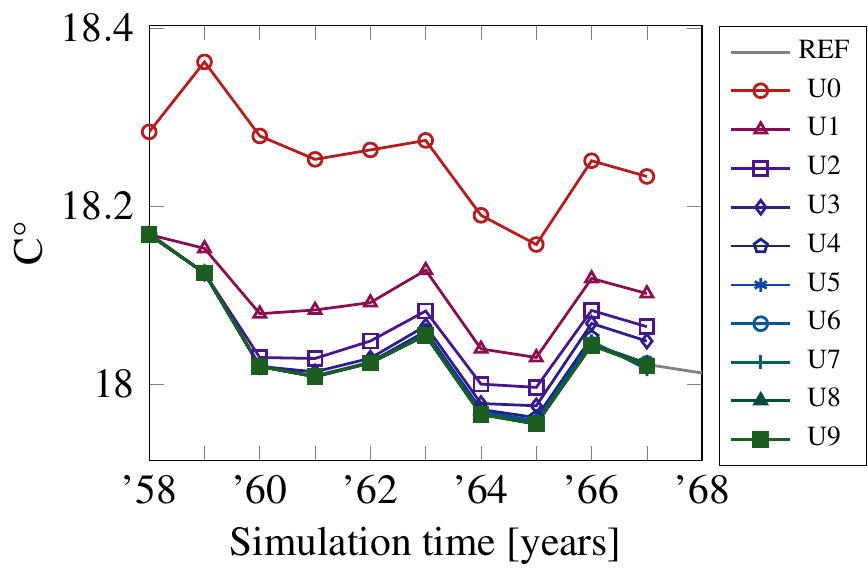}
	\caption{Sea surface temperature from 1958 to 1967 after restarting Parareal from iteration $U_3$ of test case 2.}
	\label{FIG:SSTEXP4}
\end{figure}

In Fig.\ref{FIG:SSTEXP4} the annual sea surface temperature from 1958 to 1967 is given. Although all iterations were computed, the last time slice in iteration $U_4$ suffered a failure in the coarse propagation. Further, the convergence is slightly worse than in the first experiment. With $U_4$ being incomplete it takes 5 iterations to recover the reference solution. 

\begin{figure}[H]
	\centering
	\begin{subfigure}[t]{.45\textwidth}
		\includegraphics{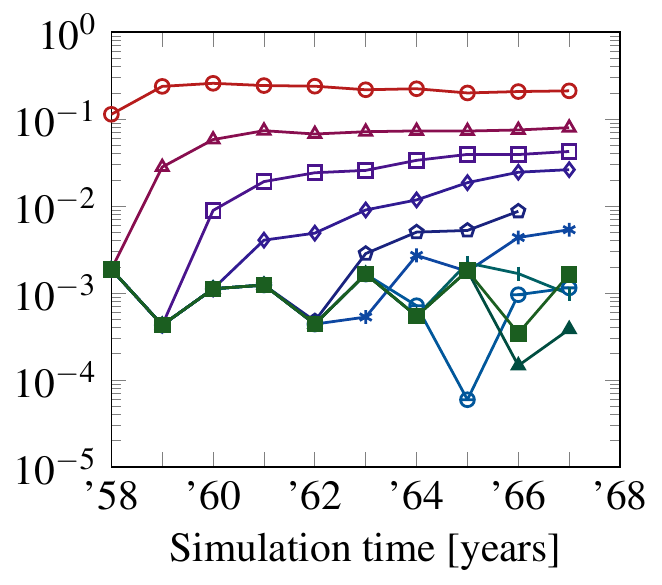}
	\end{subfigure}
	\begin{subfigure}[t]{.45\textwidth}
		\includegraphics{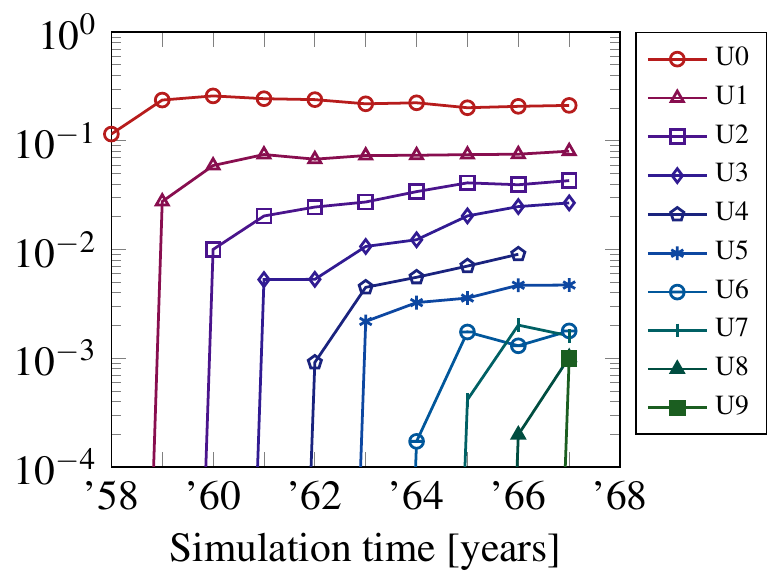}
	\end{subfigure}	
	\caption{Maximum error with respect to the reference solution startet in 1948 (left) and the repeated reference computed from the initial value given in 1958 (right).}
	\label{FIG:SSTEXP4ERR}
\end{figure}

\begin{figure}[H]
	\centering
	\includegraphics{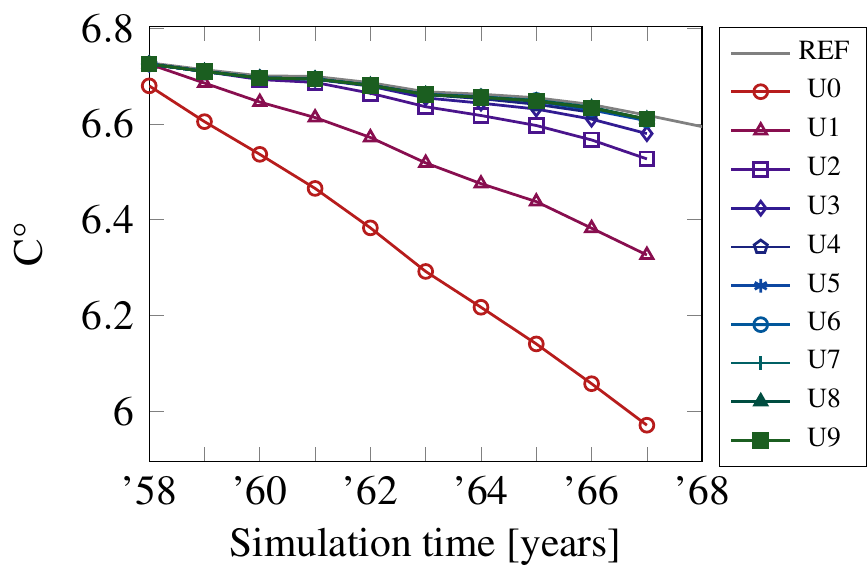}
	\caption{Annual mean temperature at 500m depth from 1958 to 1967 after restarting Parareal from iteration $U_3$ of test case 2.}
	\label{FIG:500TEMPEXP4}
\end{figure}

\begin{figure}[H]
	\centering
	\begin{subfigure}[t]{.45\textwidth}
		\includegraphics{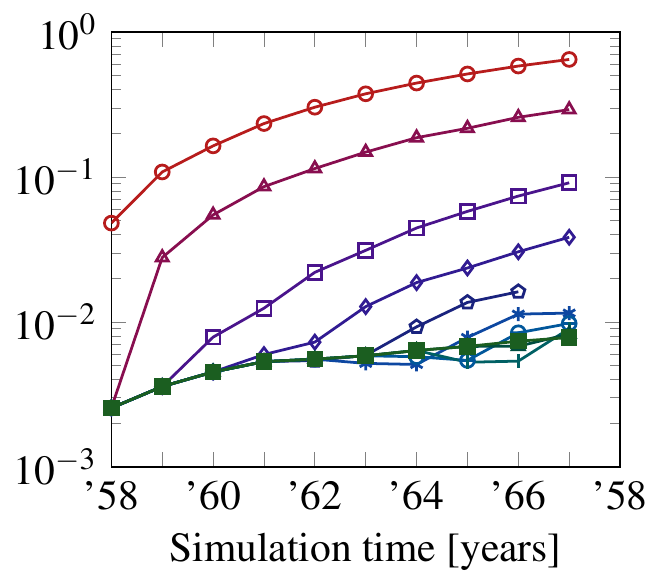}
	\end{subfigure}
	\begin{subfigure}[t]{.45\textwidth}
		\includegraphics{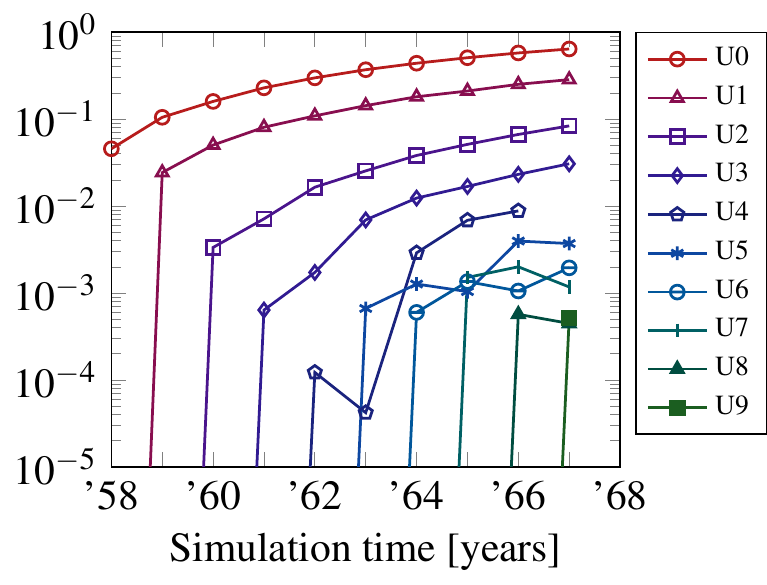}
	\end{subfigure}	
	\caption{Maximum error for annual mean temperature in 500m depth with respect to reference solution started in 1948 (left) and the repeated reference computed from the initial value given in 1958 (right).}
	\label{FIG:500TEMPEXP4ERR}
\end{figure}

The error analysis given in Fig.\ref{FIG:SSTEXP4ERR} reflect the convergence stagnation around $10^{-3}$ for both estimates as observed in the previous test cases. The only difference is the error in the first timeslice of the left panel. Since the reference solution was started in 1948 and Parareal starts with the third iteration from case 2 in 1958, the initial error depends on the approximation error of $U_3$ to reference solution. \\
For the annual mean temperature at 500m depth, see Fig.\ref{FIG:500TEMPEXP4}, we were seemingly able to reproduce the convergence behavior of the first experiment from 1948 to 1957. In Fig.\ref{FIG:500TEMPEXP4ERR} the error estimates for each iteration demonstrate is deviating from the clear error reduction until the stagnation point is reached. From 1960 to 1962 the second iteration is slightly better than $U_3$. Further, the error stagnation for the estimate with respect to the reference started in 1948 is shifted towards $10^{-2}$.

\begin{figure}[H]
	\centering
	\includegraphics{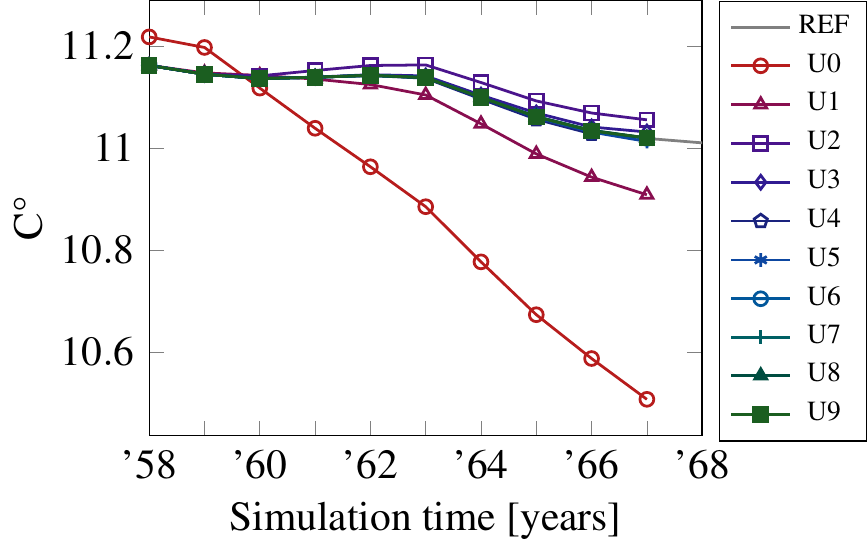}
	\caption{Annual mean temperature averaged between surface and 500m depth from 1958 to 1967 after restarting Parareal from iteration $U_3$ of test case 2.}
	\label{FIG:0500TEMPEXP4}
\end{figure}

\begin{figure}[H]
	\centering
	\begin{subfigure}[t]{.45\textwidth}
		\includegraphics{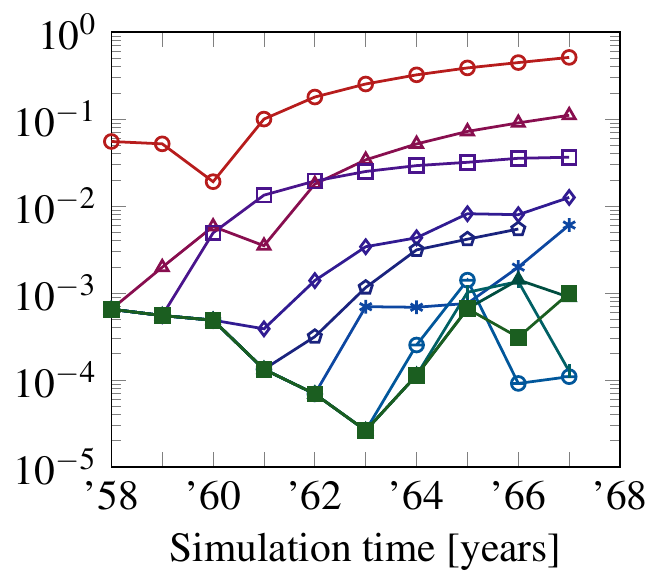}
	\end{subfigure}
	\begin{subfigure}[t]{.45\textwidth}
		\includegraphics{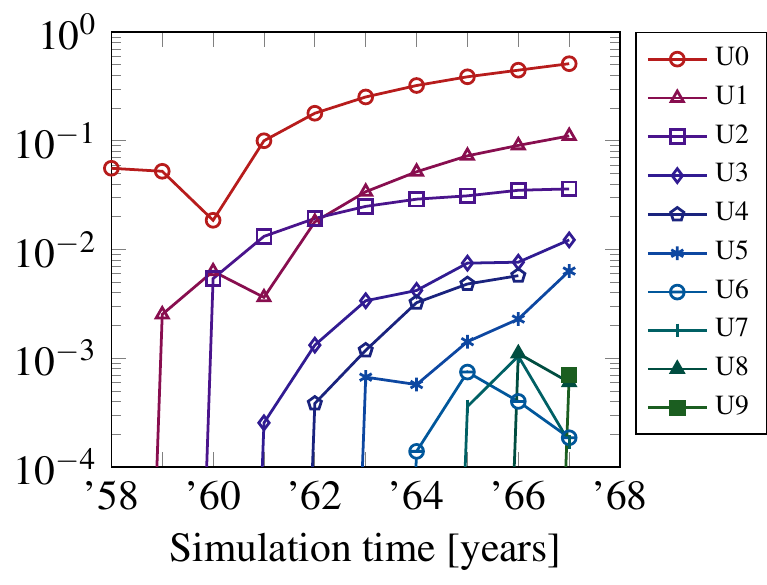}
	\end{subfigure}	
	\caption{Maximum error for annual mean temperature averaged between surface and 500m depth with respect to reference solution started in 1948 (left) and the repeated reference computed from the initial value given in 1958 (right).}
	\label{FIG:0500TEMPEXP4ERR}
\end{figure}

The averaged annual mean temperature of the layer between surface and 500m depth are in agreement with the previous results, compare Fig.\ref{FIG:0500TEMPEXP4}. The stagnation in error reduction with respect to reference solution, given in the right panel of Fig.\ref{FIG:0500TEMPEXP4ERR}, was recovered. Again, the loss in accuracy between 1960 and 1962 in $U_2$ can be observed. 

\begin{figure}[H]
	\centering
	\includegraphics{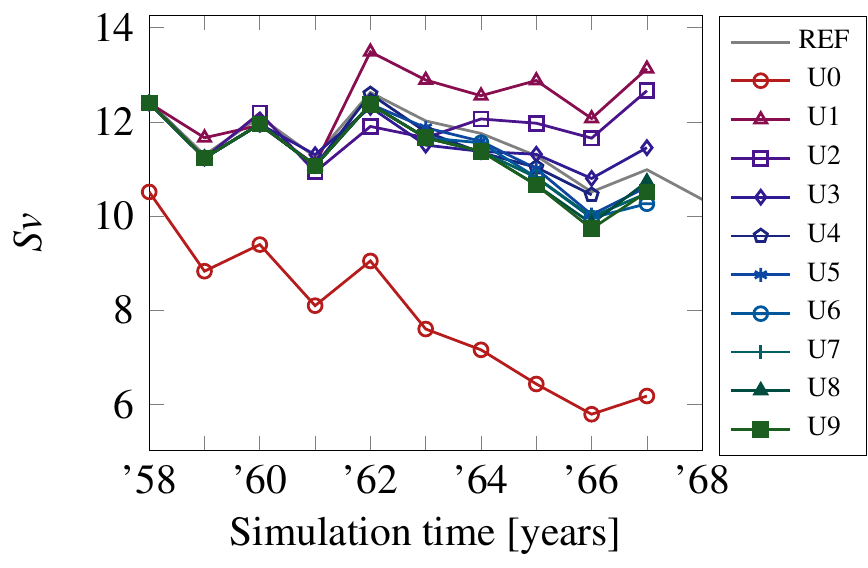}
	\caption{AMOC at 26.5N from 1958 to 1967 after restarting Parareal from iteration $U_3$ of test case 2.}
	\label{FIG:AMOCEXP4}
\end{figure}

\begin{figure}[H]
	\centering
	\begin{subfigure}[t]{.45\textwidth}
		\includegraphics{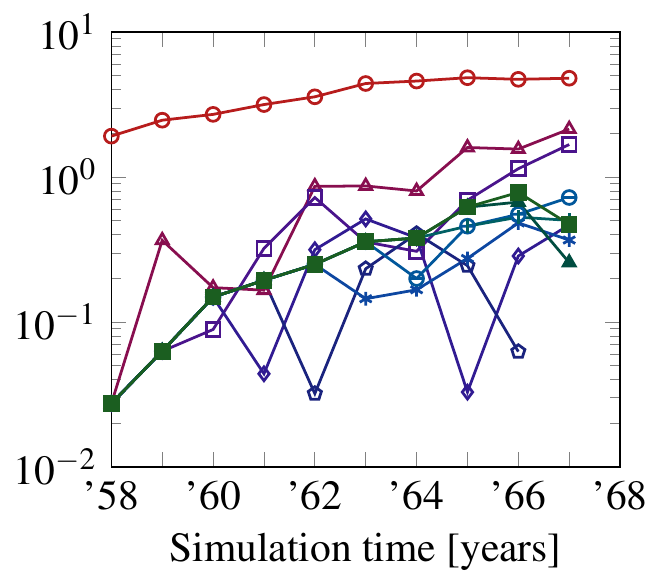}
	\end{subfigure}
	\begin{subfigure}[t]{.45\textwidth}
		\includegraphics{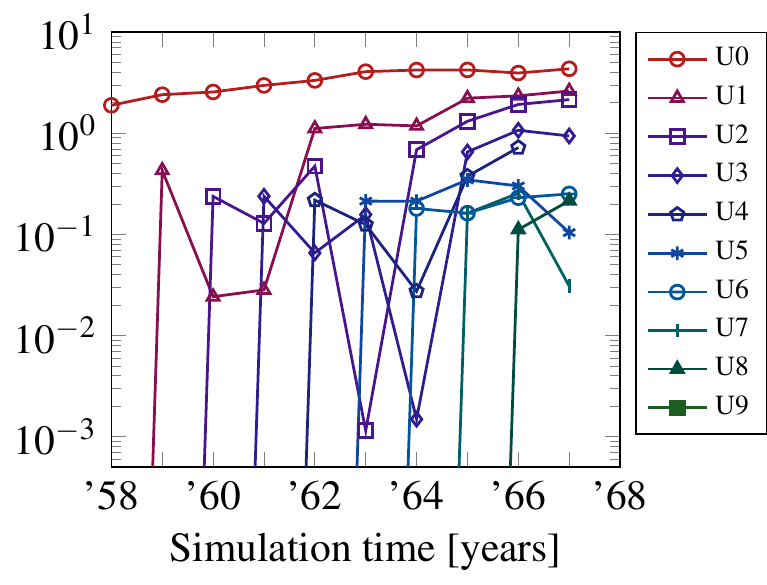}
	\end{subfigure}	
	\caption{Maximum error of approximations to the AMOC at 26.5N with respect to reference solution started in 1948 (left) and the repeated reference computed from the initial value given in 1958 (right).}
	\label{FIG:AMOCEXP4ERR}
\end{figure}

In Fig.\ref{FIG:AMOCEXP4} we find Parareal converging to a solution that is slightly off the AMOC reference at 26.5N. Aside the jump in reduction from the initial iteration $U_0$, shows the error reduction in Fig.\ref{FIG:AMOCEXP4ERR} an oscillatory behavior around the $10^{-1}$ mark. The analysis of the velocities shows most clearly that perturbation in the starting values lead to deviations in the diagnostic variables. In comparison to the convergence results for the annual mean temperature we find the post-processed velocities to be more sensitive to the initial values. By using the third iteration as an initial value for the current Parareal test case, we were not able to approximate to the desired solution. As mentioned in the introduction to this section, one has to decide with respect to his use case whether or not the solutions provided by Parareal can be considered sufficient. 

\begin{figure}[H]
	\begin{subfigure}[t]{.45\textwidth}
			\includegraphics{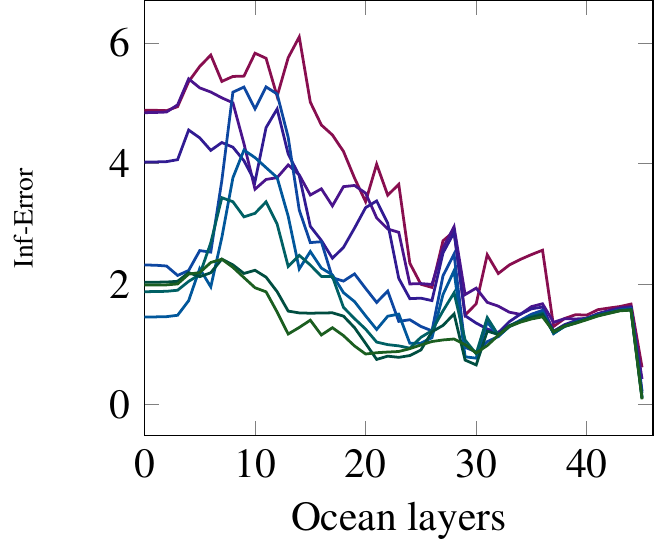}
		\end{subfigure}
	\begin{subfigure}[t]{.45\textwidth}
			\includegraphics{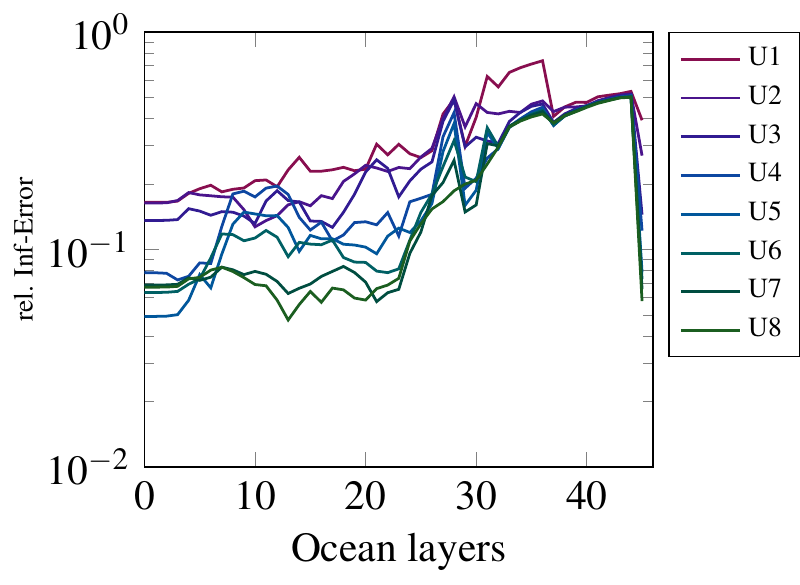}
		\end{subfigure}
	\caption{Parareal error reduction in the time-averaged horizontal temperature field over the ocean layers at the last slice. Layer 0 corresponds to the ocean surface.}
	\label{FIG:TEMPFIELDEXP4}
\end{figure}

With the amount of time steps increased to 72 spd, were able to stabilize the Parareal algorithm. Considering the error estimate for the time-averaged temperature fields in Fig.\ref{FIG:TEMPFIELDEXP4}, an reduction in the maximum error is observed to 6°C in the upper layers. In comparison to the instable iteration in experiment 2, where the maximum was 8°C, a positive impact on the stability of Parareal can be confirmed. 

\begin{figure}[H]
	\begin{subfigure}[t]{.45\textwidth}
			\includegraphics{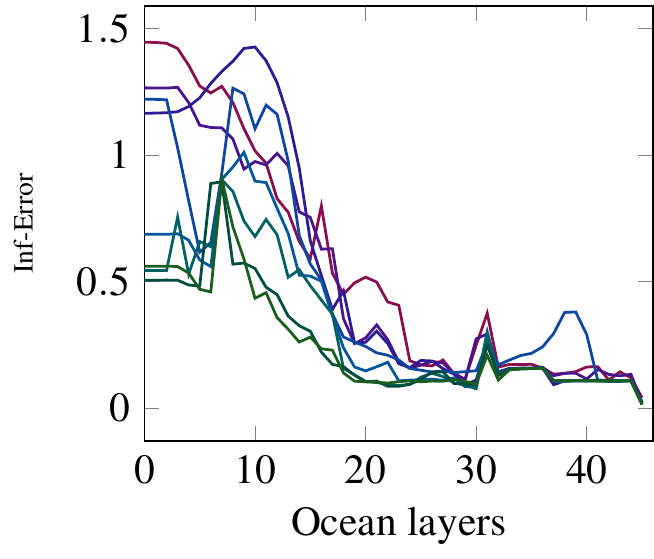}
		\end{subfigure}
	\begin{subfigure}[t]{.45\textwidth}
			\includegraphics{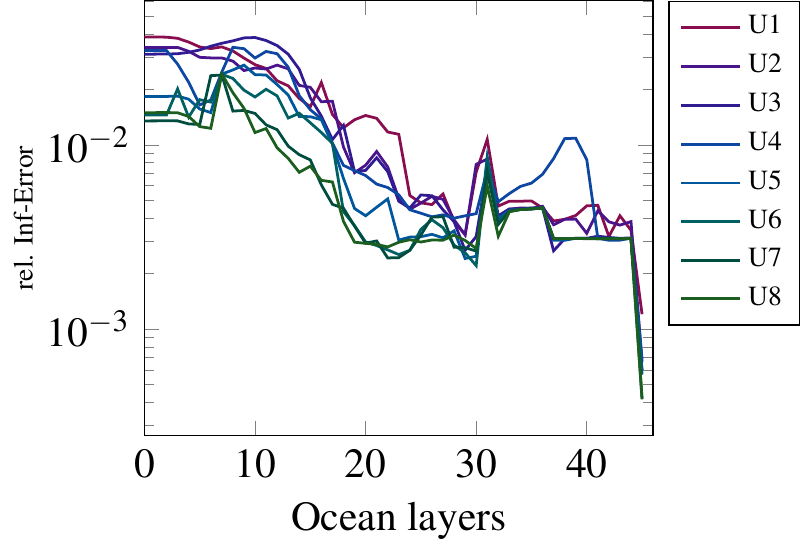}
		\end{subfigure}
	\caption{Parareal error reduction in the time-averaged horizontal salinity fields over the ocean layers at the last slice. Layer 0 corresponds to the ocean surface.}
	\label{FIG:SALTFIELDEXP4}
\end{figure}

Equally, an error reduction in the time-averaged salinity fields was observed, as depicted in Fig.\ref{FIG:SALTFIELDEXP4}. We conclude, that refining the temporal resolution for the coarse solver proves beneficial for stability of Parareal. On the other hand, it increases the run-time of FESOM2. In this case, where the time step size is refined by a factor of 2 the wall-clock time for the coarse propagator is doubled. Consequently, the a priori speedup estimate would be cut in half.

\newpage

\section{Discussion}

\begin{figure}[H]
	\centering
	\begin{subfigure}[t]{.45\textwidth}
			\includegraphics[width=\textwidth]{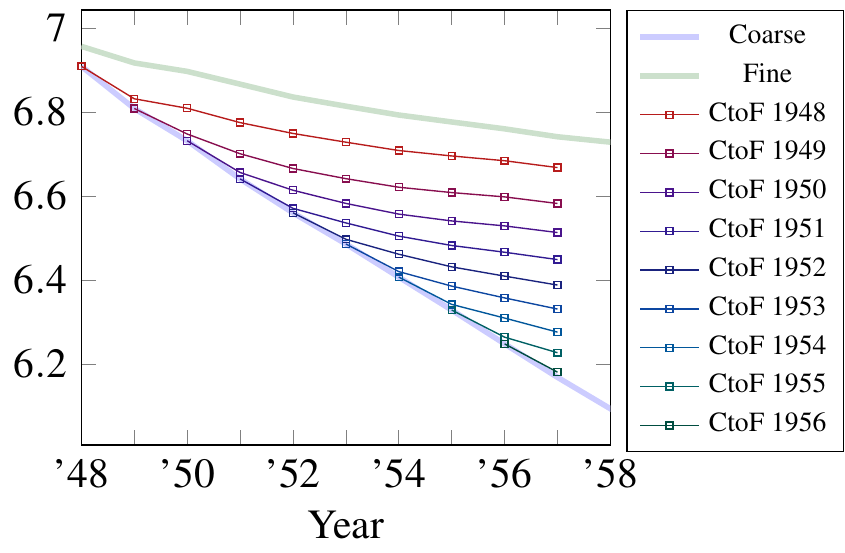}
		\end{subfigure}
	\hfill
	\begin{subfigure}[t]{.45\textwidth}
			\includegraphics[width=\textwidth]{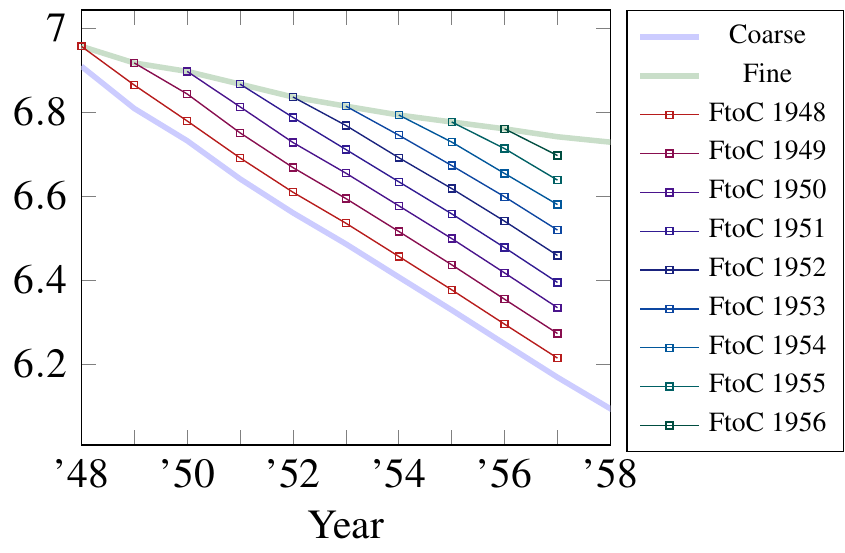}
		\end{subfigure}
	\caption{Annual mean temperature at 500m depth from 1948 to 1957 on the PI mesh (green) and FPI mesh(blue). After each year the results were interpolated from coarse to fine (left) and fine to coarse (right) mesh.}
	\label{FIG:CROSSRUNS}
\end{figure}

With the convergence results presented and the necessary steps introduced to allow for the application of the Parareal algorithm to FESOM2 we are confronted with the question of whether the first attempt of time-parallelism is a success. For this purpose, we decided to give this part the following structure:
\vspace{.25cm}
\begin{itemize}
	\item[1.] Convergence and speedup.
	\item[2.] Limits to the setting chosen in this paper.
	\item[3.] FESOM2 and mesh design in the context of compatibility to Parareal.
\end{itemize}
\vspace{.25cm}
As stated before, the question for convergence of Parareal throughout this study turns out to be twofold. There are no analytical solutions known to the problem of ocean circulation as it is treated numerically by FESOM2. Additionally, we observe the computed solutions to be highly sensitive to initial conditions, a feature that causes problems for the Parareal algorithm. The usual way to verify numerical experiments is the comparison to observational data. As both methods are prone to errors and uncertainties, the aim is to assess general geophysical trends in diagnostic variables, for which an assessment by known metrics, like the maximums norm and others, is rarely used. Evaluating the results on the basis of convergence to machine precision may seem reasonable from a numerical point of view, but it ignores the way FESOM2 and other models are used. The algorithm itself failed to reduce errors to machine precision in diagnostic variables for all experiments, except it is $K=N_t$ iterations. In this respect, it is impossible to speak of a successful application. In the context of climate research, Parareal was able to reconstruct the ocean dynamics computed on the FPI mesh. It must be taken into account that the coarse propagator used throughout this study does not offer a method of low order, but a significantly different representation of the physics by using the low resolution PI mesh. There is no simulation setting available for the PI mesh that could compute a result comparable to simulations carried out on the FPI or even CORE mesh, although the same equations are solved. In Fig.\ref{FIG:CROSSRUNS} the annual temperature at 500m depth computed on the PI and FPI mesh from 1948 to 1957 are shown. After each year we interpolated the results from coarse to fine (left) and fine to coarse (right) mesh. With the  interpolated initial values we restarted FESOM2 on the respective other mesh. These test simulations emphasize the impact of mesh resolutions and initial values on the outcome. None of the cases in Fig.\ref{FIG:CROSSRUNS} is able to recover the coarse or fine reference.
Despite the strong dependence on spatial resolution, it is possible to approximate the state of the oceans on a target grid, even if the demand for accuracy has to be weakened. Including the additional variations in diagnostics by different use of FESOM2 for the identical simulation interval (in one consecutive run or restarted after each year), we consider the application of Parareal in this regard as a success. \\
We introduced an a priori measure for the expected speedup and the maximum amount of iterations allowed. The measure does not take communication on the cluster and wall-clock time necessary for interpolation or file manipulation into account. Our wall-clock time measurements for both propagators showed a run-time ratio of $m=\tau_F / \tau_C = 3.6$, resulting in a theoretical speedup of $S_1=m/(1+1)=1.8$ for one iteration and $S_2=m/(2+1)=1.2$ for two iterations, respectively. From the results we decided that two iterations are required to adequately approximate the reference diagnostics. Consequently, the algorithm is only minimally faster. Unfortunately, the interpolation of the restart files and diagnostic output takes a considerable amount of time, even though the respective python scripts have been parallelised and are executed on the cluster for each netcdf file in parallel. The execution of the coarse solver over one year requires 120 seconds on average on Levante. Lifting the restart files from the coarse to the fine mesh takes on average the same amount of wall-clock time. The interpolation scripts were designed in a way, that allows for different interpolation types and the aspect of efficiency was therefore secondary. Any implementation, no matter how efficient, would require additional time for the interpolation, and therefore the algorithm is not able in any case to generate speedups. To obtain actual wall-clock time reductions, the FPI mesh would need to be further refined to allow for more overhead. \\

Although, Parareal was able to recover the fine solution in test case 1 and 4, stability problems occurred during the experiments 2 and 3. The error log files of FESOM2, when an abort has been forced, have indicated the blowup due to a violation of the vertical CFL number. In the introduction to our interpolation schemes, we emphasized the importance of the conservative interpolation for horizontal velocities, which are stored at the cell centers. The vertical velocity on the other hand is located midlevel at the nodes, for which the bilinear scheme was used. We suppose that, the repetitive mapping between coarse and fine meshes throughout the algorithm to be a major cause for the instabilities. We figured, that a conservative mapping must be applied for all velocities to solve this problem. For the prototype we developed for this study, the chosen scheme is sufficient, but for further investigation of Parareal for FESOM2 the velocity interpolation method needs improvement. A second contribution to the instability is assumed by discontinuities occurring during the computation of the jumps during the Parareal update step. \\
The most important limit for this configuration is certainly the selection of the meshes. They are well suited for this first attempt in parallelizing state-of-the-art software in time, but they still fall short of the expectations one would have of the results in climate research. The spatial resolution simply is too coarse. The only option would be to coarsen the CORE mesh and adjust the interpolation methods to the new configuration. Furthermore, all simulations in this study are hindcast experiments and rely on data sets from other climate models, e.g. atmosphere. Therefore, the Parareal application is bound to simulations in the past, for which the necessary observational data can be provided. If the goal would be set to simulate without data sets and predict future climate change, additional models are required to be coupled to FESOM2. This would further increase the already considerable complexity and we firmly assume that coupled simulations extend over the limit for the application of Parareal. \\

Lastly, we want to address the design of FESOM2 and the meshes used in simulations. In the context of Parareal we used FESOM2 as an black-box that provides output files for prognostic and diagnostic variables, on which manipulation and interpolation procedures were performed. The process of writing to and reading from the respective locations on the server is anything but desirable for an efficient Parareal implementation. One has to keep in mind, that FESOM2 clearly was not designed for such applications like parallelization in time. If an efficient implementation of Parareal is the goal, the entire code is needed to be rewritten, such that FESOM2 could be called by a single function and in parallel as it is needed. Furthermore, writing results to output file could be restricted, such that the overall output during the execution of Parareal is reduced to a minimum. Naturally, this would interfere with ongoing improvements to the solver by the AWI and would require additional resources. Particularly the latter does not seem feasible. \\
Although unstructured grids bring their advantages for the complex geometry of the oceans, they pose problems when in it comes to altering the mesh resolution and available interpolation schemes. We demonstrated how an existing mesh can be refined and which steps were taken to obtain the FPI mesh. The spatial refinement of a mesh is straightforward to accomplish, when the procedure proposed in this study is chosen. However, the coarsening of a grid is a significantly more complex matter, if coast lines and bathymetry are to be preserved. Which is an key requirement for the feasibility of interpolation schemes and ultimately the Parareal algorithm. The available meshes differ not only in their resolution, but also in which physical phenomenon they are designed to capture, like El Nino in the Pacific or Sea-ice around the North Pole. Coarsening these meshes makes great demands on the preservation of these very purposes. In the end, they have not been generated for being used in Micro-Macro Parareal applications and yet to do it requires a considerable effort. \\
 
In conclusion, this study demonstrates that Parareal can be successfully applied to diagnostic variables in state-of-the-art climate research. As has become clear in the course of this report, a realisation of the algorithm means considerable effort already for the prototype presented here. Without being able to provide wall-clock time reduction in the execution, a way to disclose and assess the obstacles in tackling these complex problems has been laid out.

\section{Conclusion}

This study has been conducted to demonstrate that it is possible to apply the Parareal to state-of-the-art climate models with a micro-macro variant utilizing two meshes of different spatial resolution. Based on the results of a former study on Parareal for FESOM2 \cite{philippi2022parareal} we focused on the convergence to diagnostic variables in this paper. We have given an extensive introduction on the preceding amount of work involved to make the algorithm work. The conclusions drawn from this study are given by:
\vspace{.25cm}
\begin{itemize}
	\setlength\itemsep{1em}
	\item Executing FESOM2 on two different meshes to obtain coarse and fine propagators results in significantly deviating results regarding the diagnostic variables. By applying micro-macro Parareal we were able to recover the reference diagnostics on the fine mesh.
	\item Reaching convergence of the diagnostic variables to floating point precision is possible with maximum iterations $K=N_t$, only. If convergence is evaluated from an application point of view, a sufficient maximum error around $10^{-2}$ can be achieved within $K=2$ iterations. We observed an overall stagnation in the error reduction at $10^{-3}$ beginning with iteration 4. 
	\item The vast amount of output files and interpolations performed on those files create an overhead that prevents the algorithm to reduce wall-clock times. 
	\item Generating new meshes and providing interpolation methods for FESOM2 that can be used in a Parareal configuration is a time consuming process that should be carried out with care. 
	\item For the longer simulation interval $T=20$ years we observed instabilities occurring during the FESOM2 execution on both meshes. By refining the time step size of the coarse propagator from 36 to 72 spd we have been able to restore stability, but had to accept a run-time twice as long on the coarse mesh.
\end{itemize}
\vspace{.25cm}
To our knowledge the attempt of implementing Parareal to FESOM2, or comparable climate research software, has not been made so far. Furthermore, it was not possible for us to assess in advance whether there would be any convergence in diagnostics at all. With the new FPI mesh generated and interpolation methods chosen for this study we conducted numerical experiments over $T=10$ and $T=20$ years of simulation time with time slice length $\Delta t=1$ year. Using the proposed pragmatic stopping criterion of $10^{-2}$ allows for convergence in two iterations. Encountered stability problems could be solved by increasing the amount of time steps per day by factor 2 in the coarse solver. Unfortunately, the coarse propagator requires twice as much wall-clock time and thereby destroys any prospect of run-time reduction. We assume the origin of the stability problems in the interpolation method for the vertical velocity component. In order to avoid changing temporal resolution we propose to include the vertical velocity into the conservative interpolation scheme to preserve the divergence free flow, and ultimately overcome stability issues. \\
In the context of a first assessment of micro-macro Parareal in climate research we consider this study a success. For further investigations of the algorithm applied to FESOM2 we see three major steps that have to be taken:
\vspace{.25cm}
\begin{itemize}
	\setlength\itemsep{1em}
	\item Optimize the interpolation scripts to reduce overhead for Parareal applications.
	\item Extend the conservative interpolation method to as many variables as possible, but at least for the vertical velocity component.
	\item Using CORE as the fine mesh to approximate results from current research. Sufficient coarsening of the mesh is mandatory to allow for substantial speedups.
\end{itemize}
\vspace{.25cm}
The first suggestion might be the easiest to accomplish. During the course of this study we tried several interpolation methods of which the bilinear case was the only one applicable. By hard coding this procedure the execution would be more efficient, but mesh dependent. The extension of the conservative interpolation method to vertical velocity between FESOM2 meshes is a challenge that we have not yet been able to overcome. Although it goes beyond the scope of this study, we consider its realization in the future to be indispensable. Ultimately, the micro-macro algorithm has to be tested with the CORE mesh setting to prove its applicability in standard test cases in climate research. We expect the coarsening of the CORE mesh to be subject to an optimization problem in order to preserve numerical stability. This would also mean that the previous two steps have to be adapted for a CORE configuration. \\
In conclusion, with investing the necessary time and resources to realize the proposed next steps, we expect a successful implementation of the micro-macro Parareal algorithm to FESOM2 as achievable in the future.


\printbibliography	

\end{document}